\documentclass[sn-mathphys,Numbered]{sn-jnl}

\usepackage{graphicx}%
\usepackage{multirow}%
\usepackage{amsmath,amssymb,amsfonts}%
\usepackage{amsthm}%
\usepackage{mathrsfs}%
\usepackage[title]{appendix}%
\usepackage{xcolor}%
\usepackage{textcomp}%
\usepackage{manyfoot}%
\usepackage{booktabs}%
\usepackage{algorithm}%
\usepackage{algorithmicx}%
\usepackage{algpseudocode}%
\usepackage{listings}%
\usepackage{mathtools}%
\usepackage{booktabs}%
\usepackage{longtable}%
\usepackage{empheq}%
\usepackage[capitalize, noabbrev]{cleveref}%
\crefname{subsection}{Subsection}{Subsections}%
\usepackage[binary-units]{siunitx}%
\DeclareSIUnit\crate{C}%
\usepackage[labelformat=simple]{subcaption}%
\usepackage{graphicx}%
\usepackage{tikz}%
\usepackage{caption}%
\usepackage{ulem}
\usepackage{hyperref}

\let\oldbibitem\bibitem
\def\bibitem{\vfill\oldbibitem}

\def\figpath{figures}

\usepackage{amsmath,amssymb}
\usepackage{dsfont}

\makeatletter
\newcommand{\@defs@vec}[1]{\boldsymbol{#1}}
\newcommand{\@defs@tens}[1]{\mathbf{#1}}

\DeclareMathOperator{\RelTol}{RelTol}
\DeclareMathOperator{\AbsTol}{AbsTol}
\DeclareMathOperator{\tr}{tr \!}

\newcommand{\norm}[3]{|\!| #1 |\!|_{#2}^{#3}}
\newcommand{\Norm}[3]{\big|\!\big| #1 \big|\!\big|_{#2}^{#3}}
\newcommand{\NOrm}[3]{\Big|\!\Big| #1 \Big|\!\Big|_{#2}^{#3}}
\newcommand{\NORm}[3]{\bigg|\!\bigg| #1 \bigg|\!\bigg|_{#2}^{#3}}

\newcommand{\de}{\,\mathrm{d}}
\newcommand{\ptl}{\partial}

\newcommand{\IR}{\mathbb{R}}
\newcommand{\IRp}{\mathbb{R}_{\geq0}}
\newcommand{\IN}{\mathbb{N}}

\newcommand{\II}{\mathbb{I}}
\newcommand{\Ctensor}{\mathds{C}}
\newcommand{\tensorIpi}{\mathds{I}_\Pi}

\newcommand{\vect}[1]{\@defs@vec{#1}{}}
\newcommand{\tens}[1]{\@defs@tens{#1}{}}

\newcommand{\trp}{\textsf{T}}

\newcommand{\grad}{\boldsymbol{\nabla}}

\newcommand{\tfinal}{t_\text{end}}
\newcommand{\npo}{{n+1}}
\newcommand{\halb}{\frac{1}{2}}
\newcommand{\drittel}{\frac{1}{3}}

\newcommand{\dualreduction}{\!:\!}

\newcommand{\tin}{\text{in }}
\newcommand{\ton}{\text{on }}

\newcommand{\ch}{\text{ch}}
\newcommand{\el}{\text{el}}
\newcommand{\pl}{\text{pl}}
\newcommand{\gradL}{\grad_0}
\newcommand{\divgL}{\grad_0 \!\cdot\!}
\newcommand{\something}{\boldsymbol{\cdot}}
\newcommand{\minus}{\scalebox{0.75}[1.0]{$-$}}
\newcommand{\minusone}{{\minus 1}}

\newcommand{\normalL}{\vect{n}_0}

\newcommand{\Fp}[2]
{\frac{\partial #1}{\partial #2}}

\newcommand{\ext}{\text{ext}}
\newcommand{\devi}{\text{dev}}
\newcommand{\iso}{\text{iso}}
\newcommand{\trial}{\text{tri}}
\newcommand{\myln}{\ln \!}

\newcommand{\q}{\quad}
\newcommand{\qq}{\qquad}
\newcommand{\fdg}{\,:\,}

\newcommand{\con}{c}
\newcommand{\gradu}{\gradL \vect{u}}
\newcommand{\Fpl}{\tens{F}_\text{pl}}
\newcommand{\conb}{\OB{c}}
\newcommand{\OB}[1]{\mkern 1.5mu\overline{\mkern-1.5mu#1\mkern-1.5mu}\mkern
  1.5mu}

\makeatother

\theoremstyle{thmstyleone}%
\theoremstyle{thmstyletwo}%

\theoremstyle{thmstylethree}%

\begin{document}

  \title[Modeling and Simulation]{Modeling and Simulation of
    Chemo-Elasto-\linebreak Plastically Coupled Battery Active Particles}


  \author*[1]{\fnm{Raphael} \sur{Schoof}}\email{raphael.schoof@kit.edu}

  \author[2]{\fnm{Johannes}
    \sur{Niermann}}\email{niermann.johannes@gmx.de}

  \author[2]{\fnm{Alexander} \sur{Dyck}}\email{alexander.dyck@kit.edu}

  \author[2]{\fnm{Thomas} \sur{B\"ohlke}}\email{thomas.boehlke@kit.edu}

  \author[1]{\fnm{Willy} \sur{D\"orfler}}\email{willy.doerfler@kit.edu}

  \affil*[1]{\orgdiv{Institute for Applied and Numerical Mathematics (IANM)},
    \orgname{Karlsruhe Institute of Technology (KIT)},
    \orgaddress{\street{Englerstr. 2}, \city{Karlsruhe}, \postcode{76131},
      \country{Germany}}}

  \affil[2]{\orgdiv{Institute of Engineering Mechanics (ITM), Chair for
      Continuum Mechanics}, \orgname{Karlsruhe Institute of Technology (KIT)},
    \orgaddress{\street{Kaiserstr. 10}, \city{Karlsruhe}, \postcode{76131},
      \country{Germany}}}





\abstract{

  As an anode material for lithium-ion batteries, amorphous silicon offers a
  significantly higher energy density than the graphite anodes currently used.
  Alloying reactions of lithium and silicon,
  however,
  induce large
  deformation and lead to volume changes up to 300\%.
  We formulate a thermodynamically consistent continuum model for the
  chemo-elasto-plastic diffusion-deformation
  behavior of amorphous silicon and it's alloy with lithium
  based on finite deformations.
  In this paper,
  two plasticity theories, i.e. a rate-independent theory with linear
  isotropic hardening and a rate-dependent one,
  are formulated
  to allow
  the evolution of plastic deformations and reduce
  occurring stresses.
  Using
  modern numerical techniques, such as higher order finite element
  methods as well as efficient space and time adaptive solution algorithms,
	the diffusion-deformation behavior resulting from both theories
	is compared. In order to further increase the computational efficiency,
	an automatic differentiation scheme is used, allowing for a significant
	speed up in assembling time as compared to an algorithmic linearization
	for the global finite element Newton scheme.
  Both plastic approaches lead to a
  more heterogeneous
  concentration
  distribution and
  to
  a change to tensile tangential Cauchy stresses at the
  particle surface at the end of one charging cycle.
  Different parameter studies show how an amplification of the plastic
  deformation is affected.
  Interestingly,
  an elliptical particle shows only
  plastic deformation at the
  smaller half axis.
  With the demonstrated efficiency of the applied methods, results after five
  charging cycles are also discussed and can provide indications for the
  performance of lithium-ion batteries in long term use.
}

  \keywords{lithium-ion battery, finite deformation, (visco-)plasticity, finite
    elements, numerical simulation, automatic differentiation}



  \pacs[MSC Classification]{74C15, 74C20, 74S05, 65M22, 90C33}

  \maketitle



\section{Introduction}
\label{sec:introduction}

Lithium (Li)-ion batteries gained an enormous amount of research interest in
the past two decades \cite{zhao2019review}, as a mean of storing electric energy
and propelling electro-mobility~\cite{tomaszewska2019lithium-ion,
  vasconcelos2022chemomechanics}.
However, due to the complex electro-chemo-mechanically coupled
processes
occurring during charging and discharging of Li-ion batteries, ongoing research
still aims at improving battery lifetime, reducing costs and increasing capacity
by, e.g.,
varying the materials composing the battery~\cite{zhao2019review,
  vasconcelos2022chemomechanics}.
State of the art is the usage of graphite as anode
material~\cite{zhao2019review}.
A promising candidate to be used as anode material in Li-ion batteries is
amorphous silicon~(aSi), due to its large capacity and capability to form an
amorphous $\text{Li}_x\text{Si}$
alloy~\cite{uxa2019on}
with the diffusing Li-ions, increasing battery
capacity~\cite{tomaszewska2019lithium-ion}.
A disadvantage is the large volume increase aSi particles
undergo during
alloying,
which can reach up to
$300\,\%$~\cite{zhang2011review}.
Numerous simulative studies have shown
that these large deformations are
accompanied by plastic deformations of aSi which are inherently linked to
battery lifetime and capacity, see e.g.~\cite{zhao2011large,
  di-leo2015diffusion-deformation, poluektov2018modelling,
  vasconcelos2022chemomechanics, vadhva2022towards}.
With the goal of using aSi as anode material, it is therefore imperative to
study plastic deformation mechanisms at the particle level of aSi anodes and
their interplay with battery performance during charging and discharging using
physical models and computational investigations.

To this end, geometrically and physically nonlinear chemo-mechanically coupled
continuum theories have proven to be a valuable tool, see
e.g.~\cite{zhao2011lithium-assisted, di-leo2015diffusion-deformation,
  castelli2021efficient, schoof2023simulation}.
For the mechanical part of the model, most works rely on a multiplicative split
of the deformation gradient
into a chemical, an elastic and a plastic part
using finite deformation, see e.g.
\cite[Section~10.4]{bertram2021elasticity} and
\cite[Section~8.2.2]{lubliner2006plasticity},
\cite{kolzenberg2022chemo-mechanical,
  di-leo2015diffusion-deformation}.
Discrepancies in modeling
strategies occur in the nonlinear strain measure used, ranging from the
Green--Lagrange strain tensor~\cite{castelli2021efficient, schoof2023simulation}
to the Hencky strain tensor~\cite{di-leo2015diffusion-deformation}.
In addition, several models consider plastic deformation to be
rate-dependent~\cite{di-leo2015diffusion-deformation, poluektov2018modelling},
while others rely on a rate-independent plasticity
theory~\cite{zhao2011lithium-assisted, gritton2017using}.
Unfortunately, neither the atomic-level structural evolution, nor the
mechanical behavior of the aSi during lithiation and delithiation cycles
is 	well understood~\cite{basu2019structural}.
This also holds for the detailed mechanism of plastic deformation.
However, several studies concluded that plasticity does occur during charging
and discharging.
In experimental studies, c.f.~\cite{pharr2014variation},
a rate-dependent
plastic behavior is considered to explain the observed behavior.
In contrast, a numerical study conducted on a molecular level
in~\cite{sitinamaluwa2017deformation} seems to indicate rate-independent
plasticity.
The chemical part of the models, describing diffusion of Li-ions during
charging and discharging, is based on a diffusion equation relating changes in
concentration to the gradient of the species' chemical
potential and the species' mobility~\cite{di-leo2015diffusion-deformation}.
Models differ in their approach to define the chemical contribution to the
Helmholtz free energy, where approaches either rely on
open-circuit voltage~(OCV) curves~\cite{schoof2023simulation}
or assumptions for the entropy of mixing~\cite{di-leo2015diffusion-deformation}.
In addition, the mobility is defined to be either derived as the change of the
chemical part of the chemical potential with respect to the
concentration~\cite{di-leo2015diffusion-deformation}
or the entire chemical potential~\cite{kolzenberg2022chemo-mechanical}.
The coupling of deformations and diffusion arises due to the strains induced by
Li-ions as well as the influence of mechanical stresses on the chemical
potential.

Both finite difference~\cite{zhao2011large, kolzenberg2022chemo-mechanical}
and finite element~\cite{di-leo2015diffusion-deformation, castelli2021efficient}
schemes have been proposed to discretize the resulting equations, where the
latter
have been predominantly used lately, due to their superior applicability to
complex geometries.
Solving discretized non-linear coupled systems of equations is time consuming
and
expensive in terms of computational resources, due to small mesh sizes and
small time step sizes
required to resolve all mechanisms.
In contrast to a finite difference
discretization~\cite{kolzenberg2022chemo-mechanical}
or constant spatial discretization with
uniform temporal discretization~\cite{di-leo2015chemo-mechanics}
space and time adaptive solution algorithms, such as the one proposed
in~\cite{castelli2021efficient}, allow drastic reduction in computational
resources.
In addition, parallelization schemes \cite{schoof2022parallelization}
and
automatic differentiation (AD) can reduce simulation times considerably.
The AD concept is presented for the first time in combination with the
variable-step, variable-order time integration scheme and the inelastic
constitutive theory for a battery active particle.
Introducing plastic deformation is another
challenge, as the additional variables are either
considered as degrees of freedom~\cite{kolzenberg2022chemo-mechanical}
or static condensation
is
used~\cite{wilson1974static, di-pietro2014arbitrary-order,
  di-pietro2015hybrid, oday2005superposition}
to arrive at a primal
formulation~\cite{di-leo2015diffusion-deformation,
  frohne2016efficient},
where the variables are only computed at integration point level.
A final remark concerns the choice of the chemical degree of freedom in
chemo-mechanically coupled models. In non-phase separating models, commonly,
only a single chemical degree of freedom is introduced, which is either
the concentration or the chemical potential of the diffusing
Li~\cite{barrera2016modelling, di-leo2015diffusion-deformation}.
However,
when phase separation occurs, often modeled by a Cahn--Hilliard type
theory, splitting methods~\cite{castelli2021efficient} or micromorphic
approaches~\cite{di-leo2014cahn-hilliard-type} introduce a second chemical
degree of freedom in order to circumvent the discretization of a partial
differential equation of order four and stabilize the numerical scheme.

The goal of this work is to investigate predictions of a chemo-mechanically
coupled model for large chemo-elasto-plastic deformation processes in aSi anode
particles using modern numerical techniques including space and time adaptivity,
parallelization schemes and AD.
Regarding the model, we take into account
plastic deformation of the aSi particles, where
we consider the initial yield stress to be a function of lithium
concentration~\cite{sitinamaluwa2017deformation}.
As no consensus exists in the experimental literature regarding the mechanisms
of plastic deformation in aSi, we formulate both a rate-dependent
viscoplasticity,
as well as a rate-independent plasticity theory and discuss the different
implications on particle behavior~\cite{basu2019structural, pharr2014variation,
  sitinamaluwa2017deformation}.
We use static condensation to arrive at a primal formulation for the mechanical
equations and consider plastic deformation at integration point level within
each finite element~\cite{frohne2016efficient}.
We explicitly derive a projector onto the admissible stresses, relying on the
classical return mapping method~\cite[Chapter~3]{simo1998computational}
and~\cite{simo1985consistent, simo1986return, frohne2016efficient},
which is rather straightforward for the Hencky strains used in our
theory~\cite{weber1990finite, simo1992algorithms}.
The diffusion of Li into and out of the aSi anode particle follows
classical diffusion theory.
However, instead of relying on a standard mixing entropy
ansatz for the chemical part of the free energy density,
e.g.,~\cite{di-leo2015diffusion-deformation},
we rely on a measured experimental OCV curve to model the chemical
part~\cite{kolzenberg2022chemo-mechanical, schoof2023simulation,
  schoof2022parallelization},
which has not received a lot of attention in literature.
To the best of the authors' knowledge, a combination of a measured experimental
OCV curve with a (visco-)plastic battery active particle has not been used
before.
In this work, we limit our model to show numerical results with a single-phase
mechanism after
the initial two-phase transformation of aSi in the first half cycle
as described
in~\cite{mcdowell2013in}.
This finding of~\cite{mcdowell2013in} is in contrast
to~\cite{uxa2019on},
where a phase-separation is also measured after the first
half cycle.
As a boundary condition, various charging rates (C-rates) are applied for
the lithium flux.
For solving the coupled system of equations we extend the solution scheme
proposed
in~\cite{castelli2021efficient, castelli2021numerical,
  schoof2022parallelization},
relying on a spatial and temporal adaptive algorithm. We consider
a one-dimensional computational domain with radial symmetry and
a two-dimensional elliptical computational domain
and compare stress and plastic
strain
development as well as concentration distributions
after a various number of half cycles for both rate-dependent and
rate-independent
plasticity models.
In addition, we investigate the computational performance
and numerical efficiency of our implementation scheme.

The remainder of this article is organized
as follows: in \cref{sec:theory} we introduce the theoretical basis for our
work and derive the equations describing chemo-mechanically coupled diffusion
processes in aSi anodes.
\cref{sec:numerics} summarizes the numerical approach taken in this work to
solve the derived system of equations. Subsequently, in \cref{sec:results},
we present results for various investigated cases. We close with a
conclusion and an outlook in~\cref{sec:conclusion}.




\section{Theory}
\label{sec:theory}

In a first step we review and summarize our used constitutive theory
adapted from~\cite{kolzenberg2022chemo-mechanical,
  di-leo2015diffusion-deformation,
  castelli2021efficient, schoof2023simulation}
to couple chemical, elastic and plastic material behavior.
We base our model on a thermodynamically consistent theory for the
chemo-mechanical
coupling during lithiation and delithiation.

\subsection{Finite Deformation}
\label{subsec:finite_deformation}

We start by considering a mapping~$\vect{\Phi}
\colon \IR_{\geq 0} \times \Omega_{0} \to \Omega$,
$\vect{\Phi} \left(t, \vect{X}_{0}\right) \coloneqq \vect{x}
{\,\,= \vect{X}_{0} + \vect{u}(t, \vect{X}_0)}$
from the Lagrangian
domain~$\Omega_{0}$
to the Eulerian domain~$\Omega$,
see for more information~\cite[Section~2]{holzapfel2010nonlinear},
\cite[Section~8.1]{lubliner2006plasticity},
\cite[Chapter~VI]{braess2007finite}
and
\cite{castelli2021efficient, di-leo2014cahn-hilliard-type,
  di-leo2015diffusion-deformation, kolzenberg2022chemo-mechanical,
  schoof2023simulation}.
The
deformation gradient~$\tens{F}
={\ptl \vect{\Phi}}/{\ptl \vect{X}_{0}}
=\tens{I} + \gradL \vect{u}$
with the identity tensor~$\tens{I}$
and displacement~$\vect{u}$
is multiplicatively decomposed into chemical, elastic and plastic
parts~\cite[Section~10.4]{bertram2021elasticity} and
\cite[Section~8.2.2]{lubliner2006plasticity},
\cite{kolzenberg2022chemo-mechanical,
di-leo2015diffusion-deformation}
\begin{align}
  \tens{F}
  (\gradu)
  = \tens{F}_\text{ch} \tens{F}_\text{el} \tens{F}_\text{pl}
  = \tens{F}_\text{rev} \tens{F}_\text{pl}
\end{align}
with various expressions for the volume changes
defined as
\begin{alignat}{3}
  \label{eq:volume_changes}
  J
  &= \det (\tens{F} ) = J_\text{el} J_\text{ch} J_\text{pl}
  = V_\Omega / V_{\Omega_{0}}> 0,
  \q
  &&
  J_\text{ch}
  &&= \det (\tens{F}_\text{ch}) > 0,
  \\
  J_\text{el}
  &= \det (\tens{F}_\text{el}) > 0,
  \q
  &&
  J_\text{pl}
  &&= \det (\tens{F}_\text{pl}) \stackrel{!}{=} 1,
\end{alignat}
respectively.
The chemical and elastic deformations are reversible and summarized
in~$\tens{F}_\text{rev}$.
The polar decomposition of the elastic
part~$\tens{F}_\text{el}$ of the
deformation gradient
is given by its rotational and stretch
part~\cite[Section~2.6]{holzapfel2010nonlinear}:
\begin{align}
  \label{eq:polar_decomposition}
  \tens{F}_\text{el}
  =\tens{R}_\text{el} \tens{U}_\text{el},
  \q
  \tens{R}_\text{el}^\trp \tens{R}_\text{el} =
  \tens{I},
\end{align}
with the right stretch tensor~$\tens{U}_\text{el}$
being unique,
positive definite and
symmetric. With the symmetric elastic right
Cauchy--Green tensor
\begin{align}
  \tens{C}_\text{el}
  = \tens{F}_\text{el}^\trp \tens{F}_\text{el}
  = \tens{U}_\text{el}^2,
\end{align}
the (Lagrangian) logarithmic Hencky strain can be defined as strain measure
with a spectral
decomposition
\begin{align}
  \tens{E}_\text{el}
  = \ln\left(\tens{U}_\text{el}\right)
  = \ln\Big(\sqrt{\tens{F}_\text{el}^\trp \tens{F}_\text{el}}\Big)
  = \ln\Big(\sqrt{\tens{C}_\text{el}}\Big)
  = \sum_{\alpha=1}^{3}
  \myln\left(\sqrt{\eta_{\el,\alpha}}\right)
  \vect{r}_{\el, \alpha} \otimes \vect{r}_{\el, \alpha},
\end{align}
where~$\sqrt{\eta_{\el,\alpha}}$
and~$\vect{r}_{\el, \alpha}$
are the eigenvalues and eigenvectors of~$\tens{U}_\text{el}$, respectively.
In literature, typically the Green--St-Venant~(GSV) strain tensor, often
called
\textit{the} Lagrangian strain
tensor~\cite[Section~8.1]{lubliner2006plasticity},
is used
\begin{align}
  \tens{E}_\text{el,GSV}
  = \halb \big(\tens{F}_\text{el}^\trp \tens{F}_\text{el} -
  \tens{I} \big).
\end{align}
We will later compare results obtained for both strain measures.

We consider an isotropic, volumetric swelling due to the Li concentration
with~$\tens{F}_\text{ch}$ being
defined as
\begin{align}
  \tens{F}_\text{ch}
  (\conb(t, \vect{X}_0))
  =
  \lambda_\ch
  (\conb)
  \tens{I}
  =
  \sqrt[3]{J_\ch}
  (\conb)
  \tens{I},
\end{align}
where~$\lambda_\ch = \sqrt[3]{1 + v_\text{pmv} c_{\max} \conb}$,
$v_\text{pmv}$
is the constant partial molar volume (pmv) of lithium
inside the host material~\cite{castelli2021efficient}
and
$\conb \in \left[0,1\right]$ is the normalized concentration
in the reference configuration regarding the maximal concentration of the
host material~$c_{\max}$
of the reference domain in \si{\mol\per\cubic\meter}.
Thus, the total concentration in reference
configuration is $c = c_{\max} \conb$ in \si{\mol\per\cubic\meter}.
With the chemical part~$\tens{F}_\text{ch}$ of the total deformation
gradient
and the plastic part~$\tens{F}_\text{pl}(t, \vect{X}_0)$,
introduced as an internal variable~\cite{neff2003comparison}
in more detail in~\cref{subsec:plasticity},
the elastic part~$\tens{F}_\text{el}$ of the deformation gradient is
given
by~$\tens{F}_\text{el}
(\conb, \gradu, \Fpl)
= \tens{F}_\text{ch}^{\minus 1} \tens{F} \tens{F}_\text{pl}^{\minus 1}
= \lambda_\ch^{\minus 1} \tens{F} \tens{F}_\text{pl}^{\minus 1}$.

\subsection{Free Energy}
\label{subsec:free_energy}

To obtain a thermodynamically consistent material model, which guarantees a
strictly positive entropy production, we introduce a Helmholtz free
energy~$\psi$,
being in general
a function of the lithium concentration $c$
and
the displacement gradient $\gradL
\vect{u}$~\cite{kolzenberg2022chemo-mechanical, latz2015multiscale,
  latz2011thermodynamic, schammer2021theory, schoof2023simulation}
as well as the internal
variable~$\tens{F}_\pl$~\cite{neff2003comparison},
while we assume the temperature $T$ to be constant.
Therefore, we use an additive decomposition into a chemical part~$\psi_\ch$,
being unaltered by elastic deformations,
and a mechanical
part~$\psi_\el$, depending solely on elastic deformations,
according to~\cite{kolzenberg2022chemo-mechanical, latz2015multiscale,
latz2011thermodynamic, schammer2021theory}
\begin{align}
  \label{eq:free_energy}
  \psi
  (\conb, \gradu, \Fpl)
  = \psi_\ch
  (\conb)
  + \psi_\el
  (\conb, \gradu, \Fpl)
\end{align}
in the Lagrangian frame.

Following~\cite{kolzenberg2022chemo-mechanical, schoof2022parallelization,
  schoof2023simulation},
we define the chemical part by incorporating an experimentally obtained
open-circuit voltage~($\text{OCV}$)
curve~$U_\text{OCV}(\conb)$
in \si{\volt}
\begin{align}
  \label{eq:chemical_part_free_energy}
  \rho_0 \psi_\ch
  (\conb)
  = \minus
  c_\text{max}
  \int_{0}^{\conb} \mathrm{Fa} \, U_\text{OCV}(z)
  \de z,
\end{align}
with the Faraday constant~$\mathrm{Fa}$ and
the mass density~$\rho_0$ of aSi in the reference configuration.
This measurement based model for the contribution to the chemical part of the
free energy density is in contrast with most existing works, as
often the chemical contribution is based on mixture entropy approaches for
ideal gases, see, e.g.,~\cite{di-leo2015diffusion-deformation,
	zhao2011large}. The assumption of Li behaving like an ideal gas could
be considered questionable, motivating the approach to base the chemical
contribution part on a measurement. This leads to a chemical
potential being non-symmetric to $\conb = 0.5$, c.f.~\cref{fig:voltage},
which is not the case for the mixture entropy model.
The mechanical part is given
by
a linear elastic approach
being
quadratic in the elastic Hencky strain~\cite{weber1990finite,
di-leo2015chemo-mechanics, di-leo2014cahn-hilliard-type}:
\begin{align}
  \label{eq:mechanical_part_free_energy}
  \rho_0 \psi_\el
  (\conb, \gradu, \Fpl)
  =\rho_0 \psi_\el
    (\tens{F}_\el)
  = \halb
  \tens{E}_\text{el} \dualreduction \Ctensor \left[ \tens{E}_\text{el} \right].
\end{align}
Here, $\Ctensor$ is the constant, isotropic stiffness tensor of aSi,
\linebreak
$\Ctensor
\left[ \tens{E}_\text{el} \right]
= \lambda \tr \big(\tens{E}_\text{el}\big) \tens{I} +
2 G
\tens{E}_\text{el}$
and
the first and second Lam\'{e} constants are~$\lambda
= 2 G \nu / \left(1-2\nu \right)$
and
$G = E/\big(2 \left(1+\nu\right)\big)$,
depending on the
Young's modulus~$E$
and Poisson's ratio~$\nu$
of the host material.

\subsection{Chemistry}
\label{subsec:chemistry}

Inside the host material, we use a continuity equation to describe the change
in lithium concentration
in the reference configuration
via
\begin{align}
  \ptl_t c
  = \minus \divgL \vect{N}
  \qquad
  \tin (0, \tfinal) \times \Omega_{0},
\end{align}
where
$\vect{N}
(\conb, \gradL{\mu}, \gradu, \Fpl)
\coloneqq \minus m
(\conb, \gradu, \Fpl)
\gradL{\mu}
$
is the lithium flux,
\begin{align}
  m
  =
  D \, \bigl( \ptl_c \mu
  \bigr)^{\minusone}
  \label{eq:mobility}
\end{align}
is the scalar mobility, where $m > 0$, of Li in aSi
for the considered isotropic case
and
$D$ the diffusion coefficient for lithium atoms inside the active
material~\cite{kolzenberg2022chemo-mechanical, schoof2023simulation}.
The chemical potential~$\mu(t, \vect{X}_0)$
is given as
partial derivative of the free energy density with respect to the
total concentration~$c$
using~\cref{eq:free_energy}--\eqref{eq:mechanical_part_free_energy}
\begin{align}
  \mu
  &
  = \ptl_c (\rho_0 \psi) = \ptl_c
  (\rho_0 \psi_\ch) + \ptl_c (\rho_0 \psi_\el)
  \\
  &= \minus \mathrm{Fa} \, U_\text{OCV}
  - \frac{v_\text{pmv}}{3}\lambda_\ch^{\minus 3}
  \tens{I}\dualreduction\Ctensor \left[
  \tens{E}_\text{el} \right]
  =
  \minus \mathrm{Fa} \, U_\text{OCV}
  - \frac{v_\text{pmv}}{3}\lambda_\ch^{\minus 3}
  \tr\big( \Ctensor \left[ \tens{E}_\text{el} \right] \big).
\end{align}
Following~\cite{castelli2021efficient, di-leo2015diffusion-deformation},
we apply a uniform
and constant external flux~$N_\ext$
as constant current
with either positive or negative sign for cycling the host particle in terms of
the
charging rate~($\si{C}$-rate).
The \si{C}-rate characterizes the hours of complete lithiation of the
particle:
$\si{C}\text{-rate} = 1 / t_\text{cycle}$.
So a charging rate of 1~\si{C}
means that the
battery is fully lithiated in 1~\si{\hour},
of 0.5~\si{C} in 2~\si{\hour} and so on~\cite{deng2015li-ion}.
$N_\ext$ is given in the dimension~$1/\si{h}$
after scaling with~$c_{\max} / A_\text{V}$ for the specific
surface~$A_\text{V}= \text{surface} / \text{volume}$ of $\Omega_0$ in
$\si{\square\meter\per\cubic\meter}$~\cite{castelli2021efficient}.
The simulation time~$t$
and the state of charge~($\text{SOC}$)
can be connected via
\begin{align}
  \mathrm{\text{SOC}}
  = \frac{1}{V_{\Omega_{0}}}\int_{\Omega_{0}} \frac{c}{c_\text{max}} \de
  \vect{X}_{0}
  = \frac{c_0}{c_\text{max}}
  + N_\ext t,
\end{align}
with the volume~$V_{\Omega_{0}}$
of $\Omega_{0}$
and a constant initial condition~$c_0
\in \left(0, c_\text{max}\right)$.

\subsection{Mechanics}
\label{subsec:mechanics}

The deformation in the Lagrangian domain is considered by static balance of
linear
momentum~\cite{castelli2021efficient, kolzenberg2022chemo-mechanical,
  schoof2023simulation}
\begin{align}
  \vect{0} & = \divgL \tens{P}
  \qquad
  \tin (0, \tfinal) \times \Omega_{0},
\end{align}
with the first Piola--Kirchhoff stress tensor\linebreak
$\tens{P}
(\conb, \gradu, \Fpl)
= \ptl_\tens{F}
(\rho_0 \psi_\el)
= 2\tens{F}
\ptl_{\tens{C}}
(\rho_0 \psi_\el)
=\tens{F}
\big(\tens{F}_\text{rev}^\trp \tens{F}_\text{rev} \big)^\minusone
\big(\tens{F}_\pl^\minusone\big)^\trp \tens{F}_\pl^\minusone
\Ctensor \left[ \tens{E}_\text{el} \right]$,
compare~\cref{app:derivation_p}.
The Cauchy stress~$\boldsymbol{\sigma}$
in the Eulerian frame is
linked to
$\tens{P}$ via
$\boldsymbol{\sigma} = \tens{P}\tens{F}^{\trp}/\det \left( \tens{F} \right)
$~\cite[Section~3.1]{holzapfel2010nonlinear}.
Furthermore, we introduce the Mandel
stress~$\tens{M}
(\conb, \gradu, \Fpl)
= \tens{C}_\text{rev} \tens{S}_\text{rev}
=J_\text{rev} \tens{F}_\text{rev}^\trp
\boldsymbol{\sigma} \tens{F}_\text{rev}^{\minus \trp}
=J_\el J_\ch \tens{F}_\el^\trp
\boldsymbol{\sigma} \tens{F}_\el^{\minus \trp} $
with the second Piola--Kirchhoff stress
tensor~$\tens{S}_\text{rev} = J_\text{rev} \tens{F}_\text{rev}^\minusone
\boldsymbol{\sigma} \tens{F}_\text{rev}^{\minus \trp}$,
see for further
information~\cite{kolzenberg2022chemo-mechanical}.
Based on the derivations
presented in \cite{di-leo2015diffusion-deformation,
  di-leo2014cahn-hilliard-type, di-leo2015chemo-mechanics},
for $\tens{C}_\el$ and $\dot{\tens{C}}_\el$
being coaxial and isotropic
material behavior,
a hyperelastic law relating the free energy density in the
stress-free
configuration and the Mandel
stress~$\tens{M}$
is retrieved.
In the case considered in this work $\tens{M}$ is linear
in~$\tens{E}_\text{el}$ and given by
\begin{align}
  \tens{M}
  = \Fp{(\rho_0 \psi)}{\tens{E}_\text{el}}
  = \Ctensor \left[ \tens{E}_\text{el} \right]
  &= \lambda \tr\big(\tens{E}_\text{el}\big) \tens{I}
  + 2 G
  \tens{E}_\text{el}.
\end{align}

\subsection{Inelastic Constitutive Theory}
\label{subsec:plasticity}

Following~\cite[Section~2.7]{holzapfel2010nonlinear}
and
\cite{di-leo2015diffusion-deformation},
the evolution equation for the plastic part~$\tens{F}_\text{pl}$ of the
deformation gradient
takes for elastically and plastically isotropic materials
(where a plastic spin
is negligible) the
form
\begin{align}
  \label{eq:time_evolution_plastic_deformation_gradient}
  \dot{\tens{F}}_\pl
  = \tens{D}_\pl \tens{F}_\text{pl}.
\end{align}
As mentioned in the introduction, we consider two inelastic models
and compare their influence on battery performance.
We start with a rate-independent von Mises plasticity with isotropic
hardening, which is formulated for the Mandel stress,
see~\cite[Section~2.3]{simo1998computational}
and \cite{kolzenberg2022chemo-mechanical,
  di-leo2015diffusion-deformation, kossa2009exact, simo1986return,
  frohne2016efficient}.
The yield function reads
\begin{align}
  \label{eq:yield_function}
  F_\text{Y}
  (\conb, \gradu, \Fpl, \varepsilon_\pl^{\text{eq}})
  = \norm{\tens{M}^\devi}{}{} - \sigma_\text{F}(\conb,
  \varepsilon_\pl^{\text{eq}})
  \leq 0,
\end{align}
with the deviatoric stress tensor~$\tens{M}^\devi
= \tens{M} - 1/3\tr\big(\tens{M}\big)\tens{I}$,
the Frobenius norm~$\norm{\something}{}{}$
and the yield stress~$\sigma_\text{F}
(\conb, \varepsilon_\pl^{\text{eq}}) \coloneqq
\sigma_\text{Y}(\conb) +
\gamma^\iso \varepsilon_\pl^{\text{eq}}$.
The yield stress consists of two parts:
a concentration dependent part~$\sigma_\text{Y}(\conb)$,
which will describe a softening behavior,
and a linear isotropic hardening part
with a scalar parameter~$\gamma^\iso > 0$.
$\varepsilon_\pl^{\text{eq}}(t, \vect{X}_0)$ is the
accumulated
equivalent
inelastic strain and will be introduced subsequently.
Ideal plasticity is present for~$\gamma^\iso = 0$.
The softening behavior modeled by~$\sigma_\text{Y}(\conb)$
is inspired
by~\cite{di-leo2015diffusion-deformation, sitinamaluwa2017deformation}
and is incorporated via
\begin{align}
  \sigma_\text{Y}(\conb)
  \coloneqq \sigma_\text{Y,\text{min}}
  \conb + (1 - \conb)
  \sigma_\text{Y,\text{max}}.
\end{align}
Plastic flow is only allowed if~$F_\text{Y}
=
0$.
To describe the plastic flow when this yield point is reached, we base on the
maximum plastic dissipation principle~\cite{kolzenberg2022chemo-mechanical,
  simo1998computational, han1995computational}.
With this postulate from plasticity theory we can define the associated flow
rule constraining the plastic flow to the normal direction of the yield
surface~$\ptl_{\tens{M}} F_\text{Y}$:
\begin{align}
  \tens{D}_\pl
  \big(\conb, \gradu, \Fpl,
  \varepsilon_\pl^{\text{eq}}\big)
  = \dot{\varepsilon}_\pl^\text{eq} \tens{N}_\pl
  = \dot{\varepsilon}_\pl^\text{eq} \frac{\ptl
  F_\text{Y}}{\ptl \tens{M}}
  = \dot{\varepsilon}_\pl^\text{eq}
  \frac{\tens{M}^\devi}{\norm{\tens{M}^\devi}{}{}}.
\end{align}
Here, $\tens{D}_\pl$ is the plastic strain rate
measure~\cite{di-leo2015diffusion-deformation, lubliner1986normality}.
Now, we can define the scalar equivalent plastic strain
\begin{align}
  \varepsilon_\pl^{\text{eq}}
  (t, \something)
  = \int_{0}^{
    t
  } \norm{\tens{D}_\pl}{}{} \de
  z
  = \int_{0}^{
    t
  } \dot{\varepsilon}_\pl^\text{eq} \de
  z
  ,
\end{align}
which is used to describe an increase in yield
stress~$\sigma_\text{F}(\conb,
\varepsilon_\pl^{\text{eq}})$.
Note that we always start with
$\varepsilon_\pl^{\text{eq}}(0, \something) = 0$.
To be consistent with a one-dimensional tensile test, we scale our
concentration dependent yield stress~$\sigma_\text{Y}(\conb)$
with the factor~$\sqrt{2/3}$,
compare~\cite[Section~2.3.1]{simo1998computational}.

The second model
studied
is a viscoplastic material model which was
proposed
by~\citeauthor{di-leo2015diffusion-deformation}~\cite{di-leo2015diffusion-deformation},
where a viscoplastic material behavior is considered without isotropic
hardening, i.e.~$\gamma^\iso = 0$.
This results in a formulation for the equivalent plastic strain
\begin{subequations}
  \label{eq:viscoplatic_epsplver}
  \begin{empheq}[
    left={
      \dot{\varepsilon}_\pl^\text{eq} =
      \empheqlbrace}
    ]{alignat=3}
    &0,\label{eq:viscoplatic_epsplver_a}
    && \Norm{\tens{M}^{\devi}}{}{} \leq \sigma_\text{Y}(\conb),
    \\
    &
    \dot{\varepsilon}_0
    \Bigg( \frac{\Norm{\tens{M}^{\devi}}{}{} -
      \sigma_\text{Y}(\conb)}
    {\sigma_{\text{Y}^{*}}}\Bigg)^\beta,
       \q
    &&
    \label{eq:viscoplatic_epsplver_b}
    \Norm{\tens{M}^{\devi}}{}{} > \sigma_\text{Y}(\conb),
  \end{empheq}
\end{subequations}
where $\sigma_{\text{Y}^{*}}, \dot{\varepsilon}_0$ and $\beta$ are a
positive-valued
stress-dimensioned constant, a reference tensile plastic strain rate and a
measure of the strain rate sensitivity of the material, respectively.

The classical loading and unloading conditions for rate-independent plasticity
can be conveniently expressed
via the Karush--Kuhn--Tucker (KKT)
conditions~\cite[Section~1.2.1]{simo1998computational},
\cite[Section~3.2]{lubliner2006plasticity}
and~\cite{kolzenberg2022chemo-mechanical}
\begin{align}
  F_\text{Y} \leq 0,\q
  \dot{\varepsilon}_\pl^\text{eq} \geq 0, \q
  F_\text{Y} \dot{\varepsilon}_\pl^\text{eq} = 0.
\end{align}
Compared to classical notation of loading and unloading, the process is elastic
if~$F_\text{Y} < 0$
requiring~$\dot{\varepsilon}_\pl^\text{eq} \equiv 0$ and no
plastic deformation
occurs.
The consistency condition for the evolution of inelastic strains in the case of
rate-independent
plasticity reads
\begin{align}
  \text{when } F_\text{Y} = 0: \q
  \dot{\varepsilon}_\pl^\text{eq} \geq 0, \q
  \dot{F}_\text{Y} \leq 0, \q
  \dot{\varepsilon}_\pl^\text{eq} \dot{F}_\text{Y} = 0,
\end{align}
so the plastic strain can increase during loading but not during
unloading.
In the viscoplastic case, the KKT conditions and the consistency condition
are replaced by the evolution equation of the equivalent plastic strains
as given
by~\cref{eq:viscoplatic_epsplver},
compare~\cite[Section~1.7]{simo1998computational}.




\section{Numerical Approach}
\label{sec:numerics}

In the following section we present the numerical treatment of our set of
coupled partial differential equations of~\cref{sec:theory}, i.e.
the problem formulation, the normalization of the model
parameters and the numerical solution procedure including the weak formulation,
space and time discretization and our applied adaptive solution algorithm.

\subsection{Problem Formulation}
\label{subsec:problem_statement}

Before we state our problem formulation we introduce a nondimensionalization
of the model to improve numerical stability. Our cycle time~$t_\text{cycle} = 1
/ \si{C}$-rate
depends on the $\si{C}$-rate, i.e.
the hours for charging or discharging of the particle.
Further, the particle radius~$L_0$ and the maximal concentration $c_\text{max}$
in
the
Lagrangian frame are used as reference parameters. For the yield stress
$\sigma_\text{Y}(c)$ we use the same nondimensionalization as for the Young's
modulus~$E$.
The resulting dimensionless numbers~$\tilde{{E}}$
and the \textit{Fourier number}~$\mathrm{Fo}$
relate the mechanical energy scale to the chemical energy scale and the
diffusion time scale to the process time scale, respectively.
All dimensionless variables are listed in~\cref{tab:normalization}
and will be used for model equations from now on, neglecting the accentuation
for
better readability.

\begin{table}[b]
  \centering
  \caption{Dimensionless variables of the used model equations.}
  \label{tab:normalization}
  \begin{tabular}[t]{@{}llll@{}}
    \toprule
    $\tilde{t} = t / t_\text{cycle}$  \quad \quad
    &
    $\tilde{c} = c / c_{\max}$ \quad \quad
    &
    $\tilde{E} = E / R_\text{gas} T c_{\max}$
    &
    $\widetilde{\rho_0 \psi} = \rho_0 \psi / R_\text{gas} T c_{\max}$
    \\
    $\tilde{\vect{X}}_0 = \vect{X}_0 / L_0$
    &
    $\tilde{v}_\text{pmv} = {v}_\text{pmv} c_{\max}$
    &
    $\mathrm{Fo} = D t_\text{cycle}/ L_0^2$
    &
    $\tilde{U}_\text{OCV}
    = \mathrm{Fa} \, {U}_\text{OCV} / R_\text{gas} T$
    \\
    $\tilde{\vect{u}} = \vect{u} / L_0$
    &
    $\tilde{\mu} = \mu / R_\text{gas} T$
    &
    $\tilde{N}_{\ext} = {N}_{\ext} t_\text{cycle}/ L_0 c_{\max}$
    &
    \\
    \bottomrule
  \end{tabular}
\end{table}

We state our general mathematical problem
formulation~\cite{castelli2021efficient,
  schoof2023simulation}
by solving our set of equations for the concentration~$c$,
the chemical potential~$\mu$
and the displacements~$\vect{u}$,
whereas the quantities~$\tens{F}$,
$\tens{F}_\text{el}$,
$\tens{E}_\text{el}$
and
$\tens{P}$
are calculated in dependency of the
solution
variables.
The Cauchy stresses~$\boldsymbol{\sigma}$ are computed in a postprocessing
step.
The handling of the internal variable~$\Fpl$~\cite{neff2003comparison}
is explained in more detail
in~\cref{subsec:num_sol_approach}.

The dimensionless initial boundary value problem with inequality boundary
conditions is
given
as follows:
let~$\tfinal > 0$
be the final simulation time and
$\Omega_0 \subset \IR^d$
a representative bounded electrode particle in reference configuration with
dimension~$d = 1, 2, 3$.
Find the normalized concentration~$c
\colon[0, \tfinal] \times \overline{\Omega}_0 \rightarrow
[0,1]$,
the chemical potential~$\mu
\colon[0, \tfinal] \times \overline{\Omega}_0 \rightarrow \IR$
and the displacement~$\vect{u}
\colon[0, \tfinal] \times \overline{\Omega}_0 \rightarrow \IR^d$
satisfying
\begin{subequations}
  \label{eq:ibvp}
  \begin{empheq}[left=\empheqlbrace]{alignat=4}
    \partial_t c
    &= \minus \divgL \vect{N}
    (\con, \gradL{\mu}, \gradu, \Fpl)
    && \qquad
    && \tin (0, \tfinal) \times \Omega_{0},
    \label{eq:ibvp_a}
    \\
    \mu
    &= \partial_c
    \big(\rho_0
    \psi
    (\con, \gradu, \Fpl)\big)
    &&
    && \tin (0, \tfinal) \times \Omega_{0},
    \label{eq:ibvp_b}
    \\
    \tens{0}
    &=
    \divgL \tens{P}
    (\con, \gradu, \Fpl)
    &&
    && \tin (0, \tfinal) \times \Omega_{0},
    \label{eq:ibvp_c}
    \\
    F_\text{Y}
    (\con, \gradu, \Fpl, \varepsilon_\pl^{\text{eq}})
    &\leq 0,\q
    \dot{\varepsilon}_\pl^\text{eq} \geq 0, \q
    F_\text{Y}
    \dot{\varepsilon}_\pl^\text{eq} = 0
    &&
    && \tin (0, \tfinal) \times \Omega_{0},
    \label{eq:ibvp_d}
    \\
    \vect{N} \cdot \normalL
    &= N_\ext
    &&
    && \ton (0, \tfinal) \times \partial \Omega_{0},
    \label{eq:ibvp_e}
    \\
    \tens{P} \cdot \normalL
    &= \tens{0}
    &&
    && \ton (0, \tfinal) \times \partial\Omega_{0},
    \label{eq:ibvp_f}
    \\
    c(0, \something)
    &= c_0
    &&
    && \tin \Omega_{0},
    \label{eq:ibvp_g}
    \\
    \tens{F}_\text{pl}(0, \something)
    &= \tens{I}
    &&
    && \tin \Omega_{0},
    \label{eq:ibvp_h}
    \\
    \varepsilon_\pl^{\text{eq}}(0, \something)
    &= 0
    &&
    && \tin \Omega_{0},
    \label{eq:ibvp_i}
  \end{empheq}
\end{subequations}
with a boundary-consistent initial concentration~$c_0$ and boundary
conditions for the displacement excluding rigid body motions.
Note that the
original definition of the chemical
part~$\tens{F}_\text{ch}$ of the
deformation gradient
is done in three dimensions, but all variables and equations are also
mathematically valid in dimensions $d= 1, 2$.
In this case,
the deviatoric part is computed
with the factor~$1/ d$
in contrast to the factor~$1/3$ after~\cref{eq:yield_function}.
In case of $d=2$, for example, we have no displacement and no stresses in the
third direction, so it is purely a two-dimensional computational consideration.
A final remark:
\cref{eq:ibvp_a} and~\cref{eq:ibvp_b} could be condensed to one equation.
The second equation would be treated on integration point level. However, we
choose the splitting method, detailed in~\cite{castelli2021efficient}, to
easily extend the physical behavior to phase-separating materials with a
Cahn--Hilliard approach and stabilize the numerical solution procedure.
Using only the chemical potential as solution variable would lead to an
additional calculation of the concentration with a scalar Newton method.

\subsection{Numerical Solution Procedure}
\label{subsec:num_sol_approach}

This subsection describes the way to obtain a numerical solution and especially
the handling of the KKT condition in~\cref{eq:ibvp_d}: formulating a primal
mixed
variational inequality, using static condensation to obtain a primal
formulation as well as space and time discretization, finally completed with an
adaptive solution algorithm.

\subsubsection{Weak Formulation}
\label{subsubsec:weak_formulation}

In a first step towards the numerical solution we state the weak formulation
of~\cref{eq:ibvp} as primal mixed variational inequality like
in~\cite{frohne2016efficient}.
However, it can also be derived from a minimization
problem~\cite[Section~1.4.2]{simo1998computational}
and~\cite[Section~7.3]{gromann2007numerical}.
We introduce the $L^2$-inner
product for two functions~$f$, $g \in L^2\left( \Omega_{0} \right)$
as~$\left(f, g\right) = \int_{\Omega_{0}} f g \de \vect{X}_{0}$,
for two vector
fields~$\vect{v}$, $\vect{w} \in L^2\left( \Omega_{0}; \IR^d\right)$
as~$\left(\vect{v}, \vect{w}\right) = \int_{\Omega_{0}} \vect{v} \cdot \vect{w}
\de
\vect{X}_{0}$,
and for two tensor
fields~$\tens{S}$, $\tens{T} \in L^2\left( \Omega_{0}; \IR^{d,
  d}\right)$
as~$\left(\tens{S}, \tens{T}\right) = \int_{\Omega_{0}} \tens{S} \dualreduction
\tens{T}
\de \vect{X}_{0}$
and boundary integrals with the respective boundary as subscript.
Now, we define the function
spaces~$V\coloneqq
H^k\left(\Omega_{0}, \IR \right)$,
$\vect{V}^{*} \coloneqq
H_{*}^k\big(\Omega_{0}, \IR^d\big)$
and~$\hat{\vect{V}} \coloneqq
H^k\big(\Omega_{0}, \IR_\text{sym}^{d,d}\big)$
with $k \in \IN$, $k>1$,
being sufficiently large
in order to allow a mathematically well-posed
formulation~\cite{neff2003comparison}.
$\vect{V}^{*}$ includes displacement boundary constraints for the precise
applications case
from~\cref{sec:results}.
Next,
we multiply with test
functions, integrate over~$\Omega_{0}$
and integrate by parts.
Following~\cite{castelli2021efficient,
  han1995computational, frohne2016efficient, lubliner1986normality},
we finally have the
weak
primal mixed variational inequality
formulation:
find solutions~$\lbrace c$, $\mu$, $\vect{u}\rbrace$
and internal variables $\Fpl$ and
$\varepsilon_\pl^{\text{eq}}$
with~$c, \mu \in V$,
$\ptl_t c \in L^2(\Omega_{0}, \IR)$,
$\vect{u} \in \vect{V}^{*}$,
$\Fpl \in \hat{\tens{V}}$
and~$\big(\tens{P},
\varepsilon_\pl^{\text{eq}}\big)
\in
\Big\{
L^2\big(\Omega_{0}, \IR^{d, d}\big),
L^2\big(\Omega_{0},
\IRp \big)
\fdg
F_\text{Y}
\leq 0
\Big\} \eqqcolon \tens{Y}$
such that
\begin{subequations}
  \label{eq:weak_formulation}
  \begin{empheq}[left=\empheqlbrace]{alignat=2}
    \big( \varphi,\partial_t c \big)
    &= \minus
    \Bigl(
    \gradL \varphi,
    m
    (\con, \gradu, \Fpl)
    \gradL\mu
    \Bigr)
    - \left(\varphi, N_\ext \right)_{\ptl \Omega_{0}},
    \vspace{0.2cm}
    \label{eq:weak_formulation_a}
    \\
    0 &= \minus
    (
    \varphi, \mu
    )
    + \Big( \varphi, \ptl_c
    \big(\rho_0
    \psi_\ch(c)
    \big)
    + \ptl_c
    \big(\rho_0
    \psi_\el
    (\con, \gradu, \Fpl)
    \big)
    \Big),
    \vspace{0.2cm}
    \label{eq:weak_formulation_b}
    \\
    \tens{0}
    & =
    \minus
    \Big( \gradL \vect{\xi}, \tens{P}
    (\con, \gradu, \Fpl)
    \Big),
    \label{eq:weak_formulation_c}
    \\
    \tens{0}
    & \leq
    \Big(\tens{D}_\pl,
    \tens{P}
    (\con, \gradu, \Fpl)
    - \tens{P}^{*} \Big)
    +
    \gamma^\iso
    \Big(\dot{\varepsilon}_\pl^\text{eq},
    \varepsilon_\pl^{\text{eq}*} -
    \varepsilon_\pl^{\text{eq}}
    \Big)
    \label{eq:weak_formulation_d}
  \end{empheq}
\end{subequations}
for all test functions~$\varphi \in V$,
$\vect{\xi} \in \vect{V}^{*}$
and~$\big(\tens{P}^{*}, \varepsilon_\pl^{\text{eq}*}\big) \in
\tens{Y}$.

\cref{eq:weak_formulation}
becomes a saddle point problem requiring special techniques for solving the
related
linear
system~\cite{frohne2016efficient, braess2007finite, wriggers2008nonlinear}.
However,
we prefer a version
of~\cref{eq:weak_formulation}
without~\cref{eq:weak_formulation_d}
for easier handling of the numerical studies.
Therefore, we eliminate~\cref{eq:weak_formulation_d}
by using a primal formulation with a
projector onto the set of admissible
stresses~\cite{frohne2016efficient, suttmeier2010on},
also known as static
condensation~\cite{wilson1974static, di-pietro2014arbitrary-order,
  di-pietro2015hybrid, oday2005superposition}.
Following~\cite{frohne2016efficient},
we introduce the
continuous
projector onto the admissible Mandel stress for both
considered
inelastic
constitutive theories. For the rate-independent model, it reads
\begin{subequations}
  \label{eq:projector_rate_independent}
  \begin{align*}
    \tens{M}
    (\con, \gradu, \Fpl, \varepsilon_\pl^{\text{eq}})
    =\tens{P}_\Pi
    (\con, \gradu, \Fpl, \varepsilon_\pl^{\text{eq}})
    \coloneqq
    \hspace{25.5em}\\[-2.9em]
  \end{align*}
  \begin{empheq}[
    left=
    \empheqlbrace
    ]{alignat=3}
    &\tens{M}^\trial,
    && \Norm{\tens{M}^{\trial, \devi}}{}{}
    \leq
    \sigma_\text{Y}(c) + \gamma^\iso
    \varepsilon_\pl^{\text{eq}},
    \\
    &
    \scalebox{0.89}{$
      \left[
      \frac{\gamma^\iso}{2G + \gamma^\iso}
      + \left(1-\frac{\gamma^\iso}{2G + \gamma^\iso}
      \varkappa(c, \varepsilon_\pl^\text{eq})
      \right)
      \frac{\sigma_\text{Y}(c)}{\norm{\tens{M}^{\trial, \devi}}{}{}}
      \right]
      \tens{M}^{\trial, \devi}
      $}
    &&
    \nonumber
    \\
    &
    \quad + \frac{1}{3} \tr \left( \tens{M}^{\trial} \right)
    \tens{I},
    && \Norm{\tens{M}^{\trial, \devi}}{}{}
    >
    \sigma_\text{Y}(c) + \gamma^\iso
    \varepsilon_\pl^{\text{eq}},
  \end{empheq}
\end{subequations}
with~$\varkappa(c,\varepsilon_\pl^\text{eq}) =
1-\frac{2G}{\sigma_\text{Y}(c)}\varepsilon_\pl^{\text{eq}}$
and $\tens{M}^{\trial}$ follows from a purely elastic deformation, denoted as
the
trial part of $\tens{M}$.
The projector for ideal plasticity follows with~$\gamma^\iso=0$, while
the projector for the rate-dependent viscoplastic approach is given
in~\cref{subsubsec:time_discretization}.

Finally, we can reformulate~\cref{eq:weak_formulation} using the projector
formulation for the specific plastic behavior and arrive at the primal
formulation:
for given~$\tens{F}_\pl$
and~$\varepsilon_\pl^\text{eq}$
find solutions~$\lbrace c, \mu, \vect{u} \rbrace$
with~$c, \mu \in V$,
$\ptl_t c \in L^2(\Omega_{0}, \IR)$ and
$\vect{u} \in \vect{V}^{*}$,
such that
\begin{subequations}
  \label{eq:weak_formulation_primal}
  \begin{empheq}[left=\empheqlbrace]{alignat=2}
    \left( \varphi,\partial_t c \right)
    &= \minus
    \Big( m
    (\con, \gradu, \Fpl)
     \gradL\mu,
     \gradL \varphi
    \Big)
    - \left(\varphi, N_\ext \right)_{\ptl \Omega_{0}},
    \vspace{0.2cm}
    \label{eq:weak_formulation_primal_a}
    \\
    0 &= \minus
    (
    \varphi, \mu
    )
    + \Big( \varphi, \ptl_c
    \big(\rho_0
    \psi_\ch(c)
    \big)
    + \ptl_c
    \big(\rho_0
    \psi_\el
    (\con, \gradu, \Fpl)
    \big)
    \Big),
    \vspace{0.2cm}
    \label{eq:weak_formulation_primal_b}
    \\
    \tens{0}
    & =
    \minus
    \Big( \gradL \vect{\xi}, \tens{P}
    \big( c, \gradL \vect{u}, \Fpl, \tens{P}_\Pi (c, \gradu, \Fpl,
    \varepsilon_\pl^\text{eq})
    \big)
    \Big)
    \label{eq:weak_formulation_primal_c}
  \end{empheq}
\end{subequations}
holds for all test functions~$\varphi \in V$,
$\vect{\xi} \in \vect{V}^{*}$
and~$\tens{P}
\big( c, \gradL \vect{u}, \Fpl, \tens{P}_\Pi (c, \gradu, \Fpl,
\varepsilon_\pl^\text{eq})
\big)
=
\tens{F}
\big(\tens{F}_\text{rev}^\trp \tens{F}_\text{rev} \big)^\minusone
\big(\tens{F}_\pl^\minusone\big)^\trp \tens{F}_\pl^\minusone
\tens{P}_\Pi
\big(c, \gradu, \Fpl,
\varepsilon_\pl^\text{eq}
\big)$.
This means that inequality~\cref{eq:ibvp_d}
is condensed
into~\cref{eq:weak_formulation_primal_c}
with~\cref{eq:projector_rate_independent}
or \cref{eq:projector_rate_dependent},
respectively.
Note that by using the projector~$\tens{P}_\Pi$
we have transformed the
plasticity inequality into a (non-smooth, actually Lipschitz continuous)
nonlinear equation.

\subsubsection{Space Discretization}
\label{subsubsec:space_discretization}

Next, we introduce the spatial discretization and therefore choose a
polytop approximation~$\Omega_h$
as computational domain for particle geometry~$\Omega_{0}$.
For the approximation of curved boundaries, an isoparametric Lagrangian finite
element method is chosen~\cite[Chapter~III~\S{}2]{braess2007finite}
on an admissible mesh~${\mathcal{T}_n}$.
For the spatial discrete solution we define the finite dimensional subspaces for
the basis functions with bases
\begin{subequations}
  \begin{alignat}{2}
    &V_h
    &&= {\text{span}} \{ \varphi_i \fdg i=1,\dots, N \}
    \subset {V},
    \\
    &\vect{V}_h^{*}
    &&= {\text{span}} \{ \vect{\xi}_j \fdg j=1,\dots, d N \}
    \subset \vect{V}^{*}.
  \end{alignat}
\end{subequations}
On these finite dimensional subspaces we solve for the
variables~$c_h \colon [0,\tfinal] \rightarrow \{ V_h: c_h \in
[0,1] \}$,
$\mu_h \colon [0, \tfinal] \rightarrow V_h$
and
$\vect{u}_h\colon[0, \tfinal] \rightarrow \vect{V}_h^{*}$
the spatial discrete version of~\cref{eq:weak_formulation_primal}.
Relating the discrete solution variables with the finite basis functions
\begin{subequations}
  \begin{alignat}{2}
    &c_h(t, \vect{X}_{0})
    = \sum_{i=1}^N c_i(t)\varphi_i(\vect{X}_{0}),
    \qq
    &&\mu_h(t, \vect{X}_{0}) = \sum_{j=1}^N
    \mu_j(t)\varphi_j(\vect{X}_{0}),
    \\
    &\vect{u}_h(t, \vect{X}_{0})
    = \sum_{k=1}^{d N} u_k(t)\vect{\xi}_k(\vect{X}_{0}),
  \end{alignat}
\end{subequations}
we gather all time-dependent coefficients in the vector valued function
\begin{align}
  \label{eq:discrete_solution_vector}
  \vect{y} \colon[0, \tfinal] \rightarrow \IR^{(2+d)N}, \quad
  t \mapsto \vect{y}(t) =
  \begin{pmatrix}
    \vect{c}_h(t) \\
    \vect{\mu}_h(t) \\
    \vect{u}_h(t)
  \end{pmatrix}.
\end{align}
The internal variables~$\Fpl$
and~$\varepsilon_\pl^\text{eq}$
are discretized in the same way as the solution variables and are
denoted by~$\tens{F}_{\pl, h}$
and~$\varepsilon_{\pl,h}^\text{eq}$.
Note that both internal variables~$\Fpl$
and~$\varepsilon_\pl^\text{eq}$
are not part of the solution vector$~\vect{y}$.
However, they are handled separately as introduced
in~\cref{subsubsec:time_discretization}.
Both projector formulas of~\cref{eq:projector_rate_independent}
and~\cref{eq:projector_rate_dependent}
are also valid in the discrete case.
This leads to our spatial discrete problem formulated as general nonlinear
differential algebraic equation~$\text{(DAE)}$:
for given~$\tens{F}_{\pl, h}$,
$\varepsilon_{\pl, h}^\text{eq}$
find
$\vect{y} \colon [0, \tfinal] \rightarrow \IR^{(2+d)N}$
satisfying
\begin{align}
  \label{eq:resulsting_dae}
  \hat{\tens{M}}_h \ptl_t \vect{y} - \vect{f}(t, \vect{y},
  \tens{F}_{\pl, h}, \varepsilon_{\pl, h}^\text{eq}
  )
  = \tens{0} \qquad \text{for }
  t \in
  (0, \tfinal], \qquad \vect{y}(0) = \vect{y}^0.
\end{align}
The system matrix $\hat{\tens{M}}_h$ is singular since it has only one
nonzero-block entry given
by~$\tens{M}_h =\left[(\varphi_i, \varphi_j)\right]_{i,j}$
representing the mass matrix of the finite element space~$V_h$.
The vector~$\vect{f}$
consists of matrices and tensors, given
as~$\vect{f}
\colon [0, \tfinal] \times \IR^{(2+d) N}
\times \IR_{\text{sym}}^{d,d} \times \IR_{\geq 0}
\rightarrow
\IR^{(2+d) N }$,
\begin{align}
  \label{eq:right_hand_side_dae}
  \left(t, \vect{y}\right) \mapsto \vect{f}(t, \vect{y},
  \tens{F}_{\pl, h}, \varepsilon_{\pl, h}^\text{eq}
  ) \coloneqq
  \q \q \q \q \q \q \q \q \q \q \q \q \q \q \q \q \q\q \q \q \q
  \nonumber
  \\
  \q \q \q \q \q \q
  \begin{pmatrix}
    \minus \tens{K}_m
    (\con_h, \gradL \vect{u}_h,
    \tens{F}_{\text{pl},h})
    \vect{\mu}_h
    - \vect{N}_\ext
    \\
    \minus \tens{M}_h \vect{\mu}_h
    + \vect{\Psi}_\ch(c_h) +
    \vect{\Psi}_\el
    (\con_h, \gradL \vect{u}_h,
    \tens{F}_{\text{pl},h})
    \\
    \minus
    \tens{P}_h
    \big( c_h, \gradL \vect{u}_h, \tens{F}_{\pl, h}, \tens{P}_\Pi (c_h, \gradL
    \vect{u}_h, \tens{F}_{\pl, h},
    \varepsilon_{\pl,h}^\text{eq})
    \big)
  \end{pmatrix}
\end{align}
with the same indices as above:
the mass matrix~$\tens{M}_h$,
the stiffness matrix~$\tens{K}_m
(\con_h, \gradL \vect{u}_h,
\tens{F}_{\text{pl},h})
=
\left[
\left(m
(\con_h, \gradL \vect{u}_h,
\tens{F}_{\text{pl},h})
\gradL \varphi_i, \gradL \varphi_j\right)
\right]_{i,j}$,
the vectors for the nonlinearities~$\vect{\Psi}_\ch(c_h)
= \left[
\big(\varphi_i, \ptl_c
\big(\rho_0
\psi_\ch(c_h)
\big)
\big)
\right]_{i}$
and~$\vect{\Psi}_\el
(\con_h, \gradL \vect{u}_h,
\tens{F}_{\text{pl},h})
=
\left[
\big(\varphi_i, \ptl_c
\big(\rho_0
\psi_\el
(\con_h, \gradL \vect{u}_h,
\tens{F}_{\text{pl},h})
\big)
\big)
\right]_{i}$,
$\tens{P}_h
\big( c_h, \gradL \vect{u}_h, \tens{F}_{\pl, h}, \tens{P}_\Pi (c_h, \gradL
\vect{u}_h, \tens{F}_{\pl, h},
\varepsilon_{\pl,h}^\text{eq})
\big)
=
\left[
\left( \gradL \vect{\xi}_k, \tens{P}
\big( c_h, \gradL \vect{u}_h, \tens{F}_{\pl, h}, \tens{P}_\Pi (c_h, \gradL
\vect{u}_h, \tens{F}_{\pl, h},
\varepsilon_{\pl,h}^\text{eq})
\big)
\right)
\right]_{k}$
as well as the boundary condition~$\vect{N}_\ext =
\left[
\left(\varphi_i, N_\ext\right)_{\Gamma_\ext}\right]_i$,
respectively.

\subsubsection{Time Discretization}
\label{subsubsec:time_discretization}

Before we
formulate
the space and time discrete problem
for~$\vect{y}$,
we have to consider the
time
evolution
of the internal
variables~$\tens{F}_\text{pl}$
and~$\varepsilon_\pl^\text{eq}$
and derive a time integration scheme for them.
Therefore, we update the time integration separately from the time advancing of
our~\cref{eq:resulsting_dae}.
Applying an implicit exponential map
to~\cref{eq:time_evolution_plastic_deformation_gradient}
to the continuous internal variable~$\Fpl$
leads to~\cite{neff2003comparison}
\begin{alignat}{3}
  \label{eq:exponential_map}
  &
  && \q \tens{F}_{\text{pl}}^\npo
  (\con^\npo, \gradL \vect{u}^\npo,
  \tens{F}_{\text{pl}}^n)
  &&=
  \exp\left(
  \tau_n
  \tens{D}_{\pl}^\npo\right) \tens{F}_{\text{pl}}^n
\end{alignat}
from one time step~$t_n$
to the next~$t_\npo
= t_n + \tau_n
$
with time step size~$\tau_n > 0$.
For the rate-independent plasticity we use the well known return mapping
algorithm~\cite{weber1990finite, simo1992algorithms}
in an explicit form:
\begin{align}
  \norm{\tens{M}^{\npo, \devi}}{}{} - \sigma_\text{F} (c^\npo,
  \varepsilon_\pl^{\text{eq},
    \npo})
  &= \norm{\tens{M}^{\trial, \devi}}{}{} - 2G
  \triangle_n
  \varepsilon_\pl^{\text{eq}} -
  \left(\sigma_\text{Y}(c^\npo) + \gamma^\iso
  \varepsilon_\pl^{\text{eq},
  \npo}
  \right)
  \nonumber
  \\
  &
  \stackrel{!}{=} 0
  \label{eq:return_mapping_a}
  \\
  \Longleftrightarrow \qq
  \varepsilon_\pl^{\text{eq}, \npo}
  & =\frac{\norm{\tens{M}^{\trial, \devi}}{}{} +
    2 G \varepsilon_\pl^{\text{eq},n} -
    \sigma_\text{Y}(c^\npo)
  }{2G + \gamma^\iso},
  \label{eq:return_mapping_b}
\end{align}
where the trial Mandel stress is
$\tens{M}^{\trial, \devi} = \Ctensor[\tens{E}_\el^{\trial, \devi}
(\con^\npo, \grad_0 \vect{u}^\npo, \Fpl^n)
]$
and the difference of the new time step and current time step
of~$\varepsilon_\pl^\text{eq}$
is~$\triangle_n \varepsilon_\pl^\text{eq}
= \varepsilon_\pl^{\text{eq}, \npo} -
\varepsilon_\pl^{\text{eq}, n}$.
This
is straightforward for isotropic linear hardening.
With the solution for~$\varepsilon_\pl^{\text{eq}, \npo}$
from~\cref{eq:return_mapping_b}
and the initial conditions~$\tens{F}_\text{pl}(0, \something) =
\tens{I}$
and~$\varepsilon_\pl^{\text{eq}}(0, \something) = 0$,
all necessary quantities can be updated
and~$\tens{F}_{\text{pl}}^\npo$
can be computed.
For more details regarding the time integration scheme for~$\tens{F}_\text{pl}$
in the
rate-independent case,
see~\cite[Appendix~C.5]{di-leo2015chemo-mechanics}.

Following~\cite{di-leo2015diffusion-deformation, di-leo2014cahn-hilliard-type}
for the rate-dependent viscoplastic approach,
the projector
is
  \begin{subequations}
    \label{eq:projector_rate_dependent}
    \begin{align*}
      \tens{M}
      \big(\con^\npo, \gradL \vect{u}^\npo, \Fpl^n,
      \triangle_n \varepsilon_\pl^{\text{eq}} \big)
      =\tens{P}_\Pi
      \big(\con^\npo, \gradL \vect{u}^\npo, \Fpl^n,
      \triangle_n \varepsilon_\pl^{\text{eq}} \big)
      \coloneqq
      \hspace{7.0em}\\[-2.9em]
    \end{align*}
    \begin{empheq}[
      left=
      \empheqlbrace
      ]{alignat=3}
      &\tens{M}^\trial,
      && \Norm{\tens{M}^{\trial, \devi}}{}{}
      \leq
      \sigma_\text{Y}(c^\npo),
      \\
      &
      \frac{\norm{\tens{M}^{\trial, \devi}}{}{} - 2 G \triangle_n
        \varepsilon_\pl^{\text{eq}}}{\norm{\tens{M}^{\trial,
        \devi}}{}{}}
      \tens{M}^{\trial, \devi} 	+ \frac{1}{3} \tr \left( \tens{M}^{\trial}
      \right)
      \tens{I}, \qq
      && \Norm{\tens{M}^{\trial, \devi}}{}{} > \sigma_\text{Y}(c^\npo).
    \end{empheq}
    \label{eq:viscoplastic_projector}
\end{subequations}
\!\!\!\!\!\! For
this
case, however, no explicit form can be retrieved
for the accumulated plastic strain, due to the nonlinearity
in~\cref{eq:viscoplatic_epsplver_b}.
Therefore, we have to use a scalar Newton--Raphson method for the
current time increment~$\tau_n$
for~$\triangle_n
\varepsilon_\pl^{\text{eq}}$ of the implicit Eulerian scheme:
solve~\cref{eq:viscoplatic_epsplver_b}
for~$\Norm{\tens{M}^{\trial, \devi}}{}{} > \sigma_\text{Y}(
c^\npo)$
using the relation
$\norm{\tens{M}^{\npo, \devi}}{}{}
= \norm{\tens{M}^{\trial, \devi}}{}{} - 2G
\triangle_n
\varepsilon_\pl^{\text{eq}}$
from~\cref{eq:return_mapping_a}
with the residual of~\cref{eq:viscoplatic_epsplver_b}
\begin{align}
  r_{\varepsilon}
  = \dot{\varepsilon}_0
  \bigg(
  \frac{\norm{\tens{M}^{\trial, \devi}}{}{} - 2G
    \triangle_n
    \varepsilon_\pl^{\text{eq}} -
    \sigma_\text{F}(
    c^\npo
    )}
  {\sigma_{\text{Y}^{*}}}
  \bigg)^\beta
  - \frac{
    \triangle_n
    \varepsilon_\pl^{\text{eq}}}{\tau_n}.
  \label{eq:residual_viscoplastic}
\end{align}

All computations from~\cref{eq:exponential_map}
to~\cref{eq:residual_viscoplastic} are also valid for the spatial discrete
variables
but are written for the continuous case for sake of readability.

For the time evolution of the~$\text{DAE}$~\eqref{eq:resulsting_dae},
we apply the family of numerical differentiation formulas (NDFs)
in a variable-step, variable-order algorithm, based on
the \textit{Matlab}'s
ode15s~\cite{reichelt1997matlab, shampine1997matlab, shampine1999solving,
  shampine2003solving},
since our~$\text{DAE}$
has similar properties as stiff ordinary differential
equations~\cite{castelli2021efficient}.
An error control handles the switch in time step sizes~$\tau_n$
and order.
We have chosen the same time scale for this temporal discretization and
the temporal discretization of the internal
variables~$\tens{F}_\pl$
and~$\varepsilon_\pl^{\text{eq}}$.
We arrive at the space and time discrete problem to go on
from one time
step~$t_n$
to the next~$t_\npo$:
for given~$\tens{F}_{\pl,h}^n$
and~$\varepsilon_{\pl,h}^{\text{eq},n}$
find the discrete
solution~$\vect{y}^{\npo} \approx \vect{y}(t_{\npo})$
satisfying
\begin{align}
  \label{eq:space_and_time_discretization}
  \alpha_{k_n} \hat{\tens{M}}_h \left(\vect{y}^\npo - \vect{\chi}^n\right)
  - \tau_n
  \vect{f} \left(t_\npo, \vect{y}^\npo,
  \tens{F}_{\pl, h}^{n},
  \varepsilon_{\pl,h}^{\text{eq}, n}
  \right)
  = \tens{0}
\end{align}
with~$\vect{\chi}^n$
composed of solutions on former time
steps~$\vect{y}^n, \dots, \vect{y}^{n-k}$
and
a constant~\mbox{$\alpha_{k_n}>0$}
dependent on the chosen order~$k_n$
at time~$t_n$~\cite[Section~2.3]{shampine1997matlab}.
The vector~$\vect{f}$ depends explicitly on the time~$t$
due to the time-dependent Neumann boundary
condition~$N_\ext$
due to the change of the sign for charging and
discharging, respectively.
After computing the new discrete solution~$\vect{y}^{n+1}$
we can update~$\tens{F}_{\pl,h}^n$
and~$\varepsilon_{\pl,h}^{\text{eq},n}$
to the new time steps~$\tens{F}_{\pl,h}^{n+1}$
and~$\varepsilon_{\pl,h}^{\text{eq},{n+1}}$ as already
introduced.

\subsubsection{Adaptive Solution Algorithm}
\label{subsubsec:adaptive_algorithm}

Since our~$\text{DAE}$~\eqref{eq:resulsting_dae}
is nonlinear, we apply the Newton--Raphson method and thus need to compute the
Newton update
of the Jacobian in each time step.
For this, we
have to
linearize our equations, especially the projectors
of~\cref{eq:projector_rate_independent}
and~\cref{eq:projector_rate_dependent}.
However, both projectors are not formally differentiable.
Nevertheless, they fulfill the conditions to be slantly
differentiable~\cite{frohne2016efficient, hintermuller2002primal-dual} and we
can work with the formal linearization as an approximation.
For
the
projector
defined in~\cref{eq:projector_rate_independent},
we
follow~\cite[Appendix~C.6]{di-leo2015chemo-mechanics}
and
propose
the consistent algorithmic modulus
of~$\tens{P}_\Pi
(\con, \gradu, \Fpl, \varepsilon_\pl^{\text{eq}})$
around~$\tens{E}_\text{el}^\trial$~\cite{frohne2016efficient,
  simo1985consistent}
as
\begin{subequations}
  \label{eq:linearization_rate_independent}
  \begin{align*}
    \tensorIpi
    \big(\con^\npo, \grad_0 \vect{u}^\npo, \Fpl^n,
      \varepsilon_\pl^{\text{eq},n} \big)
    \coloneqq\hspace{22.0em}\\[-2.7em]
  \end{align*}
  \begin{empheq}[
    left=
    \empheqlbrace
    ]{alignat=3}
    &
    \Ctensor_{K} + \Ctensor_{G},
    &&
    \scalebox{0.9}{$
    \NOrm{\Ctensor \left[\tens{E}_\text{el}^{\trial, \devi}
      \right]}{}{}
    \leq
    \sigma_\text{Y}(c^\npo
    ) + \gamma^\iso \varepsilon_\pl^{\text{eq},n},
    $}
    \\
    &
    \scalebox{1.0}{$
      \Big(1- \frac{\gamma^\iso}{2G + \gamma^\iso}
      \varkappa(c^\npo, \varepsilon_\pl^{\text{eq},n})
      \Big)
      \frac{\sigma_\text{Y}(c^\npo
        )}{\norm{\tens{M}^{\trial, \devi}}{}{}}
      $}
    &&
    \nonumber\\
    &
    \quad
    \scalebox{1.0}{$
      \bigg(\Ctensor_{G} - 2G
    \frac{\tens{M}^{\trial, \devi}\otimes
      \tens{M}^{\trial, \devi}}{\norm{\tens{M}^{\trial, \devi}}{}{2}}
    \bigg)
    $}
    &&
    \nonumber
    \\
    &\quad \scalebox{1.0}{$
      + \frac{\gamma^\iso}{2G + \gamma^\iso} \Ctensor_{G} +
      \Ctensor_{K},
      $}
    &&
    \scalebox{0.9}{$
    \NOrm{\Ctensor \left[\tens{E}_\text{el}^{\trial, \devi}
      \right]}{}{}
    >
    \sigma_\text{Y}(c^\npo
    ) + \gamma^\iso \varepsilon_\pl^{\text{eq},n}$}
  \end{empheq}
\end{subequations}
with
$\Ctensor_{G}
= 2G \left(\II - \frac{1}{3} \tens{I} \otimes
\tens{I}\right)$,
$\Ctensor_{K}
= K \tens{I} \otimes
\tens{I}$, the bulk modulus~$K$
and~$\varkappa$
as in~\cref{subsubsec:weak_formulation}.
The derivation for
the approximation of
the
linearization of~\cref{eq:linearization_rate_independent}
is given in~\cref{app:derivation}
as well as the formulation of the linearization for the projector
of~\cref{eq:projector_rate_dependent}.
Again, both linearizations are written for continuous functions but are
also valid for the discrete variables.

Another possibility is the automatic differentiation ($\text{AD}$) framework,
provided
by~\cite{arndt2021deal-ii}, which
offers the possibility to use an AD framework via
Sacado (a component of Trilinos)~\cite{team2020trilinos}.
This allows us to compute the partial derivative, required for the Newton
scheme, automatically instead assembling the Newton matrix by hand.
We use the tapeless dynamic forward-mode Sacado number type (one time
differentiable).
To compute the Newton update, we use a direct LU-decomposition.
The iteration number can be decreased if an appropriate
initialization is chosen. Therefore, the starting condition for the first time
step is stated in~\cref{subsec:simulation_setup}
while a predictor scheme is used during the further time
integration~\cite{shampine1997matlab}.

Finally, we follow Algorithm~1 in~\cite{castelli2021efficient}
for the space and time adaptive solution algorithm.
Both a temporal error
estimator~\cite{reichelt1997matlab, shampine1997matlab, shampine1999solving,
  shampine2003solving}
and a spatial error estimator are considered for the respective adaptivity.
To measure spatial regularity, a gradient recovery estimator is
applied for the solution variables $c$, $\mu$ and
$\vect{u}$~\cite[Chapter~4]{ainsworth2000posteriori}.
For marking the cells of the discrete triangulation for coarsening and
refinement,
two parameters~$\theta_{\mathrm{c}}$
and~$\theta_{\mathrm{r}}$
are applied with a maximum strategy~\cite{banas2008adaptive}.
Altogether, we use a mixed error control with the
parameters~$\RelTol_t$,
$\AbsTol_t$,
$\RelTol_x$
and $\AbsTol_x$.
For further details we refer to~\cite{castelli2021efficient}.




\section{Numerical Studies}
\label{sec:results}
This section deals with the investigation of the presented model
from~\cref{sec:theory}
with the numerical tools of~\cref{sec:numerics}.
For this purpose we introduce the simulation setup
in~\cref{subsec:simulation_setup}
and discuss the numerical results in~\cref{subsec:numerical_results}
with 1D and 2D simulations, which are a 3D spherical symmetric particle reduced
to the 1D unit interval
and in addition a 2D quarter ellipse reduced from a 3D elliptical nanowire,
respectively.
The latter one is chosen to reveal the
influence of asymmetric half-axis length
on plastic deformation.

\subsection{Simulation Setup}
\label{subsec:simulation_setup}

\begin{table}[t]
  \centering
  \caption{Model parameters for numerical experiments
    \cite{kolzenberg2022chemo-mechanical, schoof2022parallelization,
      di-leo2015diffusion-deformation}.}
  \label{tab:parameters}
  \begin{tabular}[t]{@{}llccc@{}}
    \toprule
    \textbf{Description} & \textbf{Symbol} &
    \multicolumn{1}{c}{\textbf{Value}} & \textbf{Unit} &
    \multicolumn{1}{c}{\textbf{Dimensionless}} \\

    \midrule

    Universal gas constant & $R_\text{gas}$ & $8.314$ &
    \si{\joule\per\mol\per\kelvin}
    & $1$ \\[\defaultaddspace]

    Faraday constant & $\mathrm{Fa}$ & $96485$ & \si{\joule\per\volt\per\mol} &
    $1$
    \\[\defaultaddspace]

    Operation temperature & $T$ & $298.15$ & \si{\kelvin} & $1$
    \\[\defaultaddspace]

    \midrule

    \multicolumn{5}{c}{Silicon} \\

    \midrule

    Particle length scale & $L_0$ & $\num{50e-9}$ &
    \si{\meter} & 1
    \\[\defaultaddspace]

    Cycle time & $t_\text{cycle}$ & $\num{3600}$ &
    \si{\second} & 1
    \\[\defaultaddspace]

    Diffusion coefficient & $D$ & $\num{1e-17}$ &
    \si{\square\meter\per\second} &
    $14.4$ \\[\defaultaddspace]

    OCV curve & $U_\text{OCV}$ & $\text{Equation~}\eqref{eq:ocv}$ &
    \si{\volt}
    & $\mathrm{Fa}/R_\text{gas} T \cdot \eqref{eq:ocv}$
    \\[\defaultaddspace]

    Young's modulus & $E$ & $\num{90.13e9}$ & \si{\pascal} & $116.74$
    \\[\defaultaddspace]

    Maximal yield stress & $\sigma_\text{Y,max}$ & $\num{8e8}$ &
    \si{\pascal} & $0.85$
    \\[\defaultaddspace]

    Minimal yield stress & $\sigma_\text{Y,min}$ & $\num{2e8}$ &
    \si{\pascal} & $0.21$
    \\[\defaultaddspace]

    Stress constant & $\sigma_{\text{Y}^{*}}$ & $\num{2e8}$ &
    \si{\pascal} & $0.21$
    \\[\defaultaddspace]

    Isotropic hardening parameter & $\gamma^\iso$ & $\num{1,0e9}$ &
    \si{\pascal} & $0.77$
    \\[\defaultaddspace]

    Tensile plastic strain rate & $\dot{\varepsilon}_0$ & $\num{2.3e-3}$ &
    \si{\second} & $8.28$
    \\[\defaultaddspace]

    Partial molar volume & $v_\text{pmv}$ & $\num{10.96e-6}$ &
    \si{\cubic\meter\per\mol} & $3.41$
    \\[\defaultaddspace]

    Maximal concentration & $c_\text{max}$ & $\num{311.47e3}$ &
    \si{\mol\per\cubic\meter} & $1$
    \\[\defaultaddspace]

    Initial concentration & $c_0$ & $\num{6.23e3}$ &
    \si{\mol\per\cubic\meter} & $\num{2e-2}$
    \\[\defaultaddspace]

    Poisson's ratio & $\nu$ & $0.22$ & \si{-} & $0.22$
    \\[\defaultaddspace]

    Strain measurement & $\beta$ & $\num{2.94}$ &
    \si{-} & $\num{2.94}$
    \\[\defaultaddspace]

    Rate constant current density & $k_0$ & $\num{0.4207}$ &
    \si{\ampere\per\square\meter} & $\num{1.0079}$
    \\
    \bottomrule
  \end{tabular}
\end{table}

As mentioned in the introduction,
amorphous
silicon is worth investigating due to its
larger energy density and is
therefore chosen as host material.
The used model parameters are listed in~\cref{tab:parameters}.
Most parameters are taken
from~\cite{kolzenberg2022chemo-mechanical, schoof2022parallelization},
however, for the plastic deformation we pick and
adapt the parameters from~\cite{di-leo2015diffusion-deformation}
such that the yield stress is in the range of~\cite{zhang2019stress}.
In particular, the ratio between~$\sigma_\text{Y,max}$
and~$\sigma_\text{Y,min}$ is maintained.
If not otherwise stated, we
charge with $1~\si{C}$ and discharge with $\minus 1~\si{C}$.
Following~\cite{kolzenberg2022chemo-mechanical},
we charge between~$U_{\max} = 0.5~\si{V}$
and $U_{\min} = 0.05~\si{V}$,
corresponding to an initial concentration of~\mbox{$c_0=0.02$}
and a duration of $0.9~\si{h}$ for one half cycle which is one
lithiation of the host particle.
The OCV curve $U_\text{OCV}(c)$ for silicon is taken
from~\cite{chan2007high-performance}
and defined as~$U_\text{OCV} \colon (0, 1) \rightarrow \IR_{> 0}$
with
\begin{align}
  \label{eq:ocv}
  U_\text{OCV}(c) & \coloneqq
  \frac{\minus0.2453 \, c^3-0.005270 \, c^2+0.2477 \, c+0.006457}
  {c+0.002493}.
\end{align}
The curve is depicted in~\cref{app:comparison_voltage_curve}
in~\cref{fig:voltage}.

\subsubsection{Geometrical Setup}
\label{subsubsec:geometries}

\begin{figure}[t]
  \centering
  \begin{subfigure}[htbp]{0.48\textwidth}
    \vspace{0.28cm}
    {\includegraphics[width=\textwidth,page=1]
      {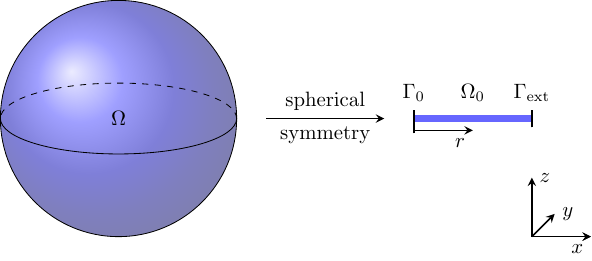}}
    \subcaption{Dimension reduction for the 1D spherical symmetric case.}
    \label{fig:geometry_1d}
  \end{subfigure}%
  \hspace{0.04\textwidth}%
  \begin{subfigure}[htbp]{0.48\textwidth}
    {\includegraphics[width=\textwidth,page=2]
      {\figpath/schematics/particle_obstacle_computational_domain}}
    \subcaption{2D quarter ellipse with
    different lengths of the half-axes.}
    \label{fig:geometry_2d}
  \end{subfigure}
  \caption{Computational domains in 1D in~\subref{fig:geometry_1d}
    and 2D in~\subref{fig:geometry_2d}
    used for the numerical simulations.}
  \label{fig:geometries}
\end{figure}

We proceed by presenting two computational domains and present the boundary
conditions due to new artificial boundaries.
We choose a representative 3D spherical particle and reduce the computational
domain to the 1D unit
interval~$\Omega_0=\left(0,1\right)$,
with a new artificial boundary~$\Gamma_0$, compare~\cref{fig:geometry_1d},
in the particle center with a no flux condition and zero displacement:
\begin{align}
  \vect{N} \cdot \normalL = 0, \qquad u = {0} \qquad
  \ton (0, \tfinal) \times \Gamma_0.
\end{align}
To ensure the radial symmetry, we adapt the quadrature weight to~$\de
\vect{X}_{0} = 4 \pi r^2 \de r$
in the discrete finite element formulation.
In the 1D domain it is consistent to assume that the fields vary solely along
the
radius $r$.
For the calculation of the displacement gradient, we refer
to~\cite[Appendix~B.2.1]{castelli2021numerical}.
As stated above, the initial concentration is $c_0=0.02$,
which leads to a one-dimensional stress-free radial
displacement~$u_0 = r \left(\lambda_\ch(c_0) - 1\right)$.
It follows that the initial chemical potential
is~$\mu_{0} =
\ptl_c
\big(\rho_0\psi_\ch(c_0)\big)$.

In~\cref{fig:geometry_2d}, the 2D simulation case is shown in terms of a
quarter ellipse.
We choose a factor of 0.6 for the short
half-axis compared to the longer half-axis.
We create this geometry by considering a 3D nanowire with no stresses and
no changes
in~$z$-direction as well as symmetry around the $x$- and $y$-axes.
Here, further artificial boundaries on~$\Gamma_{0,x}$
and on~$\Gamma_{0,y}$
with no flux conditions and only radial displacement:
\begin{subequations}
  \begin{alignat}{4}
    &\vect{N} \cdot \normalL &&= 0, \qquad u_y &&= {0},\qquad &&
    \ton (0, \tfinal) \times \Gamma_{0,x},\\
    &\vect{N} \cdot \normalL &&= 0, \qquad u_x &&= {0}, &&\ton (0,
    \tfinal) \times \Gamma_{0,y},
  \end{alignat}
\end{subequations}
have to be introduced.
We make use of an isoparametric mapping for the representation of the curved
boundary on~$\Gamma_\ext$.
Again, we choose a constant initial concentration~\mbox{$c_0 = 0.02$}
and a chemical
potential~$\mu_{0} = \ptl_c
\big(\rho_0\psi_\ch(c_0)\big)$. However, we use
for
the displacement the condition~$\vect{u}_0 = \vect{0}$.

\subsubsection{Implementation Details}
\label{subsubsec:implementation_details}
All numerical simulations are executed with an isoparametric fourth-order
Lagrangian finite element method and all integrals are evaluated through a
Gau\ss--Legendre quadrature formula with six quadrature points in space
direction.
Our code implementation is based on the finite element
library~\texttt{deal.II}~\cite{arndt2021deal-ii}
implemented in~\texttt{C++}.
Further, we use the interface to the Trilinos
library~\cite[Version~12.8.1]{team2020trilinos}
and the UMFPACK package~\cite[Version~5.7.8]{davis2004algorithm}
for the LU-decomposition for solving the linear equation systems.
A desktop computer with \SI{64}{\giga\byte}~RAM,
Intel~i5-9500~CPU,
GCC compiler version 10.5
and the operating system Ubuntu 20.04.6 LTS
is used as working machine.
Furthermore, OpenMP Version 4.5 is applied for shared memory parallelization
for assembling the Newton matrix, residuals and spatial estimates and
message passing interface (MPI)
parallelization with four MPI-jobs for the 2D simulations with
Open~MPI~4.0.3.
Unless otherwise stated, we choose for the space and time adaptive algorithm
tolerances of~$\RelTol_t = \RelTol_x = \num{1e-5}$,
$\AbsTol_t = \AbsTol_x = \num{1e-8}$,
an initial refinement of seven and of four for the 1D case and the 2D
case, respectively,
an initial time step size~$\tau_0 = \num{1e-6}$
and a maximal time step size $\tau_{\max} = \num{1e-2}$.
For the marking parameters of local mesh coarsening and refinement,
$\theta_{\mathrm{c}} = 0.05$ and $\theta_{\mathrm{r}}=0.5$
are set.
A minimal refinement level of five is applied for the 1D simulations, in order
to achieve
a parameterization as general as possible for all different 1D simulation.
Due to limited regularity of the nonlinearity of the plastic deformation, we
limit
the maximal order
of the adaptive temporal adaptivity by two.

\subsection{Numerical Results}
\label{subsec:numerical_results}
In this section we consider the numerical results of the 1D spherical symmetric
particle and the 2D quarter ellipse computational domain.
We analyze the computed fields such as stresses and concentrations as well as
the computational	performance of our presented model and implementation scheme.
Further, we compare the computational times for using the derived linearization
of the projector formulation and an automatic differentiation~(AD) technique.

\subsubsection{1D Spherical Symmetry}
\label{subsubsec:1d_spherical_symmetry}

In a first step, we analyze the effect of plastic deformation on the
chemo-physical behavior of the 1D domain depicted in \cref{fig:geometry_1d}.
Detailed studies for purely elastic behavior can be found in, e.g.,
\cite{schoof2022parallelization, schoof2023simulation} and are included in this
study for comparison.

\begin{figure}[t]
  \centering
  {\includegraphics[width=\textwidth,page=1]
    {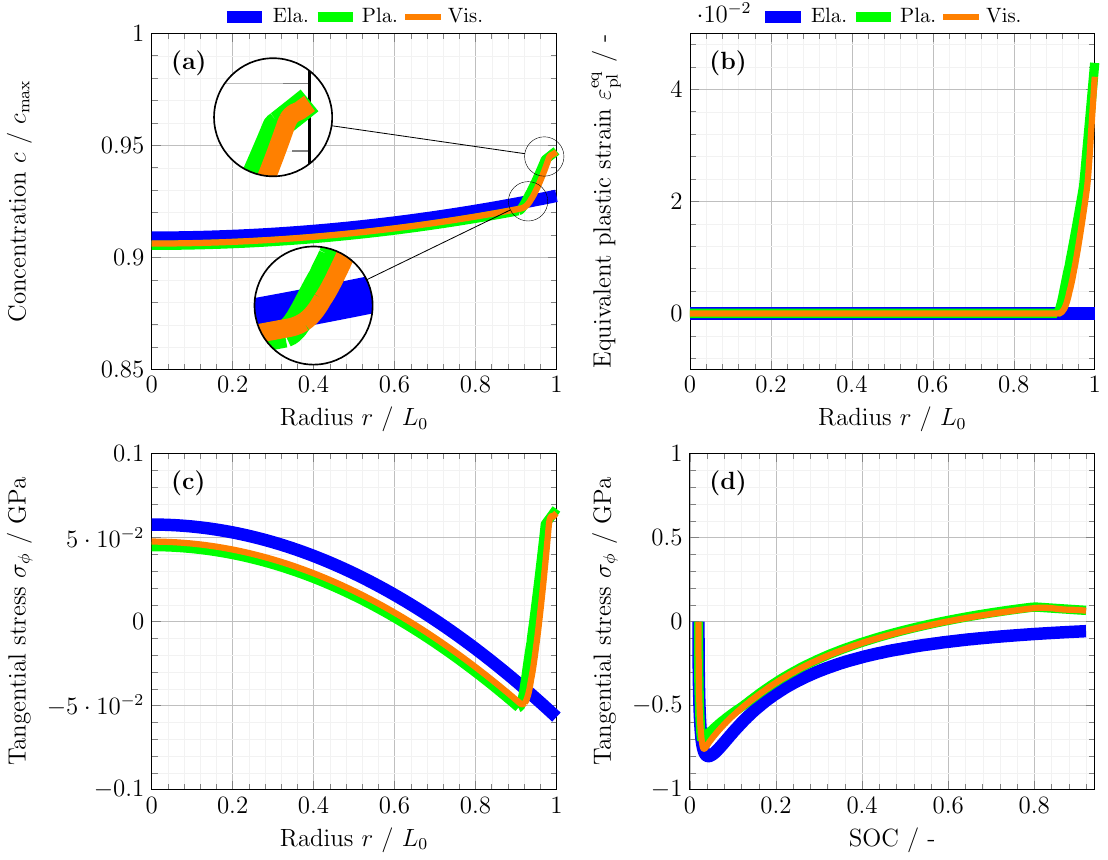}}
  \caption{
    Numerical results for the elastic (Ela.), plastic (Pla.) and
    viscoplastic (Vis.) approaches
    of the 1D radial symmetric case at~$\text{SOC} = 0.92$
    over the particle radius~$r$:
    concentration~$c$ in~(a),
    equivalent plastic strain~$\varepsilon_\pl^{\text{eq}}$
    in~(b) and
    tangential Cauchy
    stress~$\sigma_\phi$ in~(c)
    as well as
    tangential Cauchy stress~$\sigma_\phi$ at the particle
    surface~$r=1.0$
    over $\text{SOC}$ in~(d).
  }
  \label{fig:1d_one_half_cycle}
\end{figure}

\begin{figure}[t]
  \centering
  {\includegraphics[width=\textwidth,page=1]
    {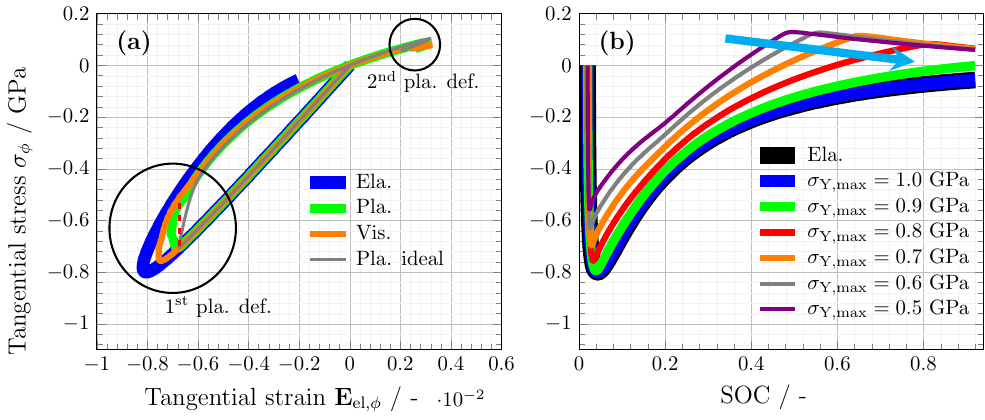}}
  \caption{Influence of the different plasticity
    approaches in~(a) and of the maximal yield stress of
    the viscoplasticity approach in~(b) for the tangential Cauchy
    stress~$\sigma_\phi$
    at the particle surface~$r=1.0$.}
  \label{fig:comparison}
\end{figure}

\textbf{Physical results after one half cycle.}
In~\cref{fig:1d_one_half_cycle},
we compare the numerical results
for the concentration, the plastic strains as well as the
tangential Cauchy
stresses
between the elastic, plastic and viscoplastic model
after one lithiation of the host particle, that means one half cycle
at~$\text{SOC} = 0.92$ over the particle radius~$r$.
The changes of the concentration profiles due to plastic deformation are
displayed in~\cref{fig:1d_one_half_cycle}(a).
It is clearly visible that the
concentration gradient
increases in both the plastic and viscoplastic case
in the vicinity of the particle surface, whereas lower
concentration values occur in the particle center compared to the elastic case.
This is due to the limited
maximal
stress in the plastic models.
The lower stresses inside the particle lead to a lower mobility in the particle
interior, c.f.~\cref{eq:mobility}.
This gradient in the Li mobility leads to the observed pile up of Li atoms at
the
particle
surface.
Comparing the plastic and viscoplastic case, a smoother transition from elastic
to plastic is visible and therefore a little shift of the concentration values
to the particle surface, see the lower magnifying glass
in~\cref{fig:1d_one_half_cycle}(a).
This is commonly
referred to as
viscoplastic regularization
\cite[Section~1.7.1.4.]{simo1998computational}.
The second magnifier shows an area, where the slope of the concentration profile
changes close to the particle surface, revealing a second plastic deformation
process	during lithiation, which is also observable at the change in slope of
the
equivalent plastic strain, c.f.~\cref{fig:1d_one_half_cycle}(b).
This becomes more apparent, when the tangential
Cauchy
stresses are
investigated as a function of the $\text{SOC}$,
c.f.~\cref{fig:1d_one_half_cycle}(d)
and~\cref{fig:comparison}(b).

\cref{fig:1d_one_half_cycle}(b)
shows the equivalent plastic
strain~$\varepsilon_\pl^{\text{eq}}$,
revealing plastic deformations near the particle surface.
During lithiation, the particle deforms plastically twice. The first
plastic deformation process occurs in the initial stages of charging at
low~$\text{SOC} \leq 0.13$,
c.f.~\cref{fig:1d_one_half_cycle}(d),
and leads to equivalent
plastic strains of~$3.4\,\%$.
Upon further lithiation, this process repeats at SOC~$\approx 0.8$,
c.f.~\cref{fig:1d_one_half_cycle}(d)
and the magnitude of the equivalent plastic strain increases to $4\,\%$.
For the tangential Cauchy stress~$\sigma_\phi$, c.f.
\cref{fig:1d_one_half_cycle}(c) and (d),
there is a change of the stress direction from compressive stress to tensile
stress at the particle surface for the plastic cases compared to the elastic
case.
This change in sign for the tangential stresses occurs in the area, where the
particle undergoes plastic deformation beforehand.
The heterogeneous plastic deformation thus leads to an eigenstrain, that results
in tensile stresses near the particle surface, which cannot be observed in the
elastic case.
This means that for an almost fully lithiated particle, there is a significant
shift in the stress development and the plastic deformation leads to tensile
stresses close to the particle surface.
This crucial change in the stress profile
at a small spatial area of the particle is important to recognize for the
battery life time.
\cref{fig:1d_one_half_cycle}(d) displays the tangential Cauchy stress at
the particle surface versus the SOC over the complete half cycle.
Directly after the start compressive Cauchy stresses occur where the elastic
approach increases to larger values compared to the plastic approaches, which
are limited due to the onset of plastic deformations.
The viscoplastic model shows larger negative tangential stresses, i.e. an
overstress above the yield stress, when compared to the rate-independent
models, also allowing larger elastic strains, c.f.~\cref{fig:comparison}(a).
After reaching the maximal values the stresses reduce in all cases.
However, the plastic approaches predict tensile stresses
around~$\text{SOC}=0.6$.
At~$\text{SOC} \approx 0.8$, the particle deforms plastically a second time,
reducing the maximal tensile stresses in tangential direction.
Striking is the difference between the elastic case revealing only compressive
tangential Cauchy stresses compared to the plastic approaches featuring also
tensile tangential Cauchy stresses.
The radial Cauchy stress is not plotted since we have stress-free boundary
condition at the particle surface.
These findings are qualitatively comparable to numerical results
from~\cite{mcdowell2016mechanics},
Fig.~4(c),
and~\cite{zhao2011large},
Fig.~5(d),
compare also the numerical
results at~$\text{SOC} = 0.1$
in~\cref{fig:comparison_stresses_plastic_literature}
of~\cref{app:comparison_stress}.
We want to point out that both the plastic and the viscoplastic model lead to
almost equivalent numerical results at the end of the first half cycle.
This is, however, not the case for results after multiple half cycles, which is
outlined below.
Before we proceed by investigating the influence of several material
parameters and multiple half cycles,
we also compare
the Green--St-Venant strain tensor or Lagrangian strain with our used
logarithmic strain tensor.
Both approaches predict almost identical results,
see~\cref{app:comparison_strain},
as also observed in the numerical results
in~\cite{zhang2018sodium}.

\begin{figure}[b]
  \centering
  {\includegraphics[width=\textwidth,page=1]
    {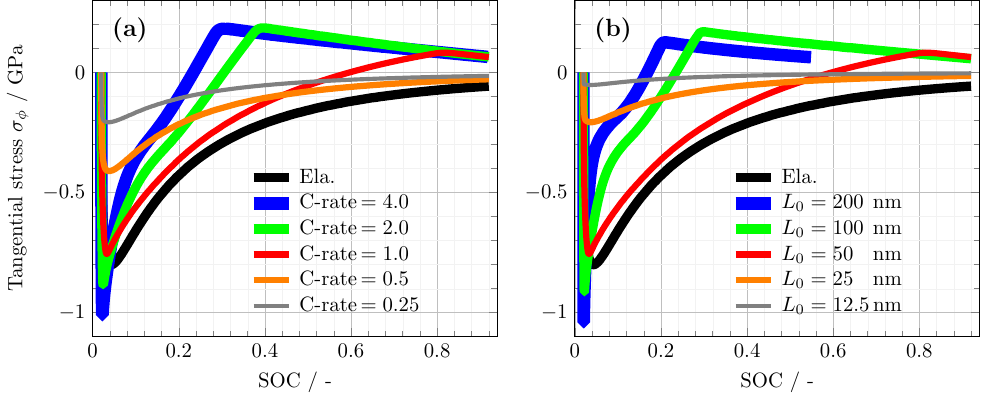}}
  \caption{Tangential Cauchy stress~$\sigma_\phi$
    over SOC at the particle surface~$r=1.0$
    for varying~\si{C}-rate in~(a)
    and particle size in~(b).}
  \label{fig:parameter_study}
\end{figure}

\textbf{Parameter studies.}
In a next step we take a closer look at the
stress-strain curves of the different mechanical approaches and
analyze the influence of the maximal yield stress~$\sigma_\text{Y,max}$ on the
tangential Cauchy stress~$\sigma_\phi$
at the particle surface~$r=1.0$
in~\cref{fig:comparison}.
Furthermore, we compare the dependency on the $\si{C}$-rate and the particle
size in~\cref{fig:parameter_study}.
The AD concept is used for all parameter studies, see the
next subsection \textit{numerical efficiency} for more details.
The comparison
in~\cref{fig:comparison}(a)
shows the effects of the different mechanical approaches:
elastic, plastic and viscoplastic deformation
as well as
ideal plastic deformation with~$\gamma^\iso = 0$.
For the ideal plastic case, we choose a uniform grid with ten refinements and
a backward differentiation formula (BDF) time stepping of order two to
increase numerical performance.
Comparing the elastic and all plastic cases, the maximal compressive Cauchy
stress of the elastic solution is not reached in the plastic cases, however,
the stresses decrease rapidly after reaching the yield stress.
Again, the viscoplastic overstress becomes apparent,
c.f.~\cref{fig:1d_one_half_cycle}(d).
Further, the influence of the concentration dependent yield
stress is clearly visible, since with a constant yield stress, the curve would
just move straight upwards (red dotted reference line) instead featuring a
shift to the right.
The tangential stress decreases rapidly after an initial plastic deformation,
see also~\cref{fig:comparison}(b),
so that no further plastic deformation occurs.
In contrast to the elastic model, tangential tensile stresses occur for large
$\text{SOC}$
values in all plasticity models, which lead to another onset of plastic
deformation.
For larger $\text{SOC}$
values, only the plasticity approaches show tensile tangential Cauchy stresses
with a second plastic deformation at the end of the first lithiation.
As stated below, the viscoplastic model including the viscoplastic
regularization leads to better numerical properties.
In addition, the experiments of
\citeauthor{pharr2014variation}~\cite{pharr2014variation}
indeed indicate some sort of rate-dependent inelastic behavior, which, in the
opinion of
the authors, makes the usage of a viscoplastic model more plausible.
Therefore, we continue our investigation with the viscoplastic model unless
otherwise stated.
%
\cref{fig:comparison}(b) compares the influence of a varying maximal
yield stress~$\sigma_\text{Y,max}$
on the tangential stress at the particle surface over the~$\text{SOC}$.
For smaller values of~$\sigma_\text{Y,max}$,
the particle starts yielding at smaller tangential stresses, leading to an
overall decrease in the observed minimal tangential stresses~$\sigma_\phi$
at small $\text{SOC}$.
In contrast, the earlier the plastic deformation occurs, the larger the
tangential tensile stresses at higher $\text{SOC}$.
In addition, the decrease in yield stress with increasing concentrations can be
observed, as the tangential tensile stresses lead to further plastic
deformation at lower stress levels,
indicated with the cyan arrow in~\cref{fig:comparison}(b).
However, no plastic deformation is visible for a maximal yield stress
of~$1.0~\si{GPa}$, which shows a purely elastic response.

Next, we analyze the dependency of the plastic deformation on
different~$\si{C}$-rates and particle sizes in~\cref{fig:parameter_study}(a)
and (b).
Again, the tangential Cauchy stress~$\sigma_\phi$ is plotted over
the~$\text{SOC}$.
For fast charging batteries, high $\si{C}$-rates
are desirable for a comfortable user experience.
However, \cref{fig:parameter_study}(a)
shows that for higher $\si{C}$-rates higher tangential Cauchy stresses arise
which lead especially for higher~$\text{SOC}$
values to a large area of plastic deformation.
The smaller the~$\si{C}$-rate, the lower
the occurring stresses and the
smaller the plastically deformed regions in the anode particle.
For the parameter set considered in this study, decreasing the $\si{C}$-rate
by $50\,\%$
leads to purely elastic deformations.
Similar results are observed, when an increasing particle diameter is
considered, as in~\cref{fig:parameter_study}(b):
the larger the particles, the greater the stresses and the area of plastic
deformation.
This results from larger concentration gradients for larger particles since the
lithium needs longer to diffuse to the particle center.
The simulation for particle size~$L_0 = 200~\si{nm}$
terminates at~$\text{SOC} \approx 0.55$,
since at this simulation time a concentration of one is reached at the particle
surface.
For the largest particle radius, the heterogeneity in the
concentration profile is even more pronounced. This also explains the
apparent lower yield stress when the $L_0 = 200~\si{nm}$ and $L_0 = 100~\si{nm}$
curves are
compared. In the larger particle, the concentration at the boundary is
larger at smaller SOCs, due to the heterogeneity in the mobility, as
outlined above.
In conclusion, \cref{fig:parameter_study}
shows
that small particles and
low~$\si{C}$-rates
are preferred in order to avoid high stresses and irreversible plastic
deformations.

\begin{figure}[b]
  \centering
  {\includegraphics[width=\textwidth,page=1]
    {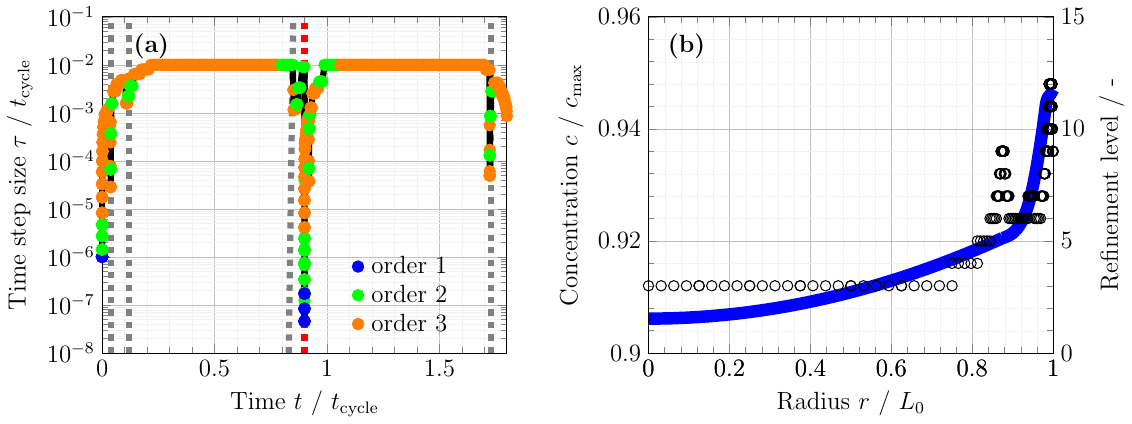}}
  \caption{Advantages of the adaptive solution algorithm:
    time step size~$\tau_n$
    over simulation time~$t$
    for two half cycles (one lithiation and one delithiation)
    in~(a)
    and
    concentration~$c$
    as well as refinement level of the spatial discretization over the
    particle
    radius~$r$
    at~$\text{SOC} = 0.92$
    in~(b).
  }
  \label{fig:numerical_efficiency}
\end{figure}

\textbf{Numerical efficiency.}
To show the capabilities of the adaptive solution algorithm
of~\cref{subsubsec:adaptive_algorithm},
we consider in~\cref{fig:numerical_efficiency}(a)
the time step size~$\tau_n$
and the used order of the NDF
multi-step
procedure
over two half cycles, i.e., one lithiation and one delithiation step
with~$\tfinal= 1.8~\si{\hour}$
in total, and
in~\cref{fig:numerical_efficiency}(b) the refinement level for the spatial
refinement after one lithiation at~$\text{SOC}=0.92$,
comparing the lithium concentration distribution over the particle
radius~$r$.
Here, we use a maximal order for the time adaptivity of three
for~\cref{fig:numerical_efficiency}(a)
and a minimal refinement level of three in~\cref{fig:numerical_efficiency}(b),
respectively.
During the first lithiation, there are three changes in the time step size and
used time order:
after starting with order one and switching to order two and three, the first
plastic deformation arises at around~$t = 0.04~\si{\hour}$
at the first vertical gray reference line.
Here, the time step sizes decrease and the used order goes down to two.
After recovering to larger step sizes and orders again,
the plastic deformation has ended at the second gray line
at~$t = 0.12~\si{\hour}$
and the particle deforms elastically again.
Then, a large time range with the maximal time step size~$\tau_{\max}$ and
maximal order is passed through.
Shortly before the end of the first half cycle,
the step sizes and order decrease again (third gray line).
In this instant, plastic deformation occurs again, due to the tensile stresses
at the particle surface.
Changing the external lithium flux direction from charging to discharging is not
trivial for the adaptive algorithm:
the order decreases to one and the time step sizes drops over five orders of
magnitude at the red reference line.
After recovering from this event, one further plastification occurs
shortly after the red reference line before the maximal time step
size and the maximal order are reached again.
At the time of~$t=1.7~\si{h}$,
indicated with the fourth gray
reference line, the next plastic deformation happens accompanied with a
reduction of time step sizes and lower order.
In total, there are two further plastifications during delithiation.
In~\cref{fig:numerical_efficiency}(b),
the focus is on the spatial adaptivity. We see the concentration
distribution~$c$
over the particle radius~$r$
and the refinement level of the cells of our triangulation.
It is clearly visible, that areas with larger concentration gradients have a
higher refinement level due to the used gradient recovery error estimator.

\begin{table}[b]
  \captionsetup{width=1.0\textwidth}
  \caption{Comparison of
    number of time steps,
    average number of Newton iterations per time step (Av. Newton),
    assembling time
    for the Newton matrix
    and
    the right hand side of the linear
    system,
    between
    the derived linearization
    and
    AD techniques
    for the 1D and the 2D simulation case
    as well as other modern numerical techniques for the
    2D case.}
  \label{tab:ad}
  \begin{tabular}[t]{@{}p{4.9cm}ccc@{}}
    \toprule
    \textbf{Method} & \textbf{Time steps} & \textbf{Av.
    Newton} &
    \textbf{Assembling time} \\
    \midrule
    Linearization 1D & 227 &
    1.29
    & \num{3.96e0} $\si{\sec}$
    \\[\defaultaddspace]
    AD 1D & 229 &
    1.27
    & \num{1.80e0} $\si{\sec}$
    \\[\defaultaddspace]
    Linearization 2D & 314 &
    1.18
    & \num{4.03 e3}
    $\si{\sec}$
    \\[\defaultaddspace]
    AD 2D & 338 &
    1.26
    &
    \num{3.93 e2}
    $\si{\sec}$
    \\
    AD 2D without MPI & 338 &
    1.26
    &
    \num{1.50 e3}
    $\si{\sec}$
    \\
    AD 2D without MPI, uniform in space & 302 &
    1.17
    &
    \num{3.53e3}
    $\si{\sec}$
    \\
    AD 2D without MPI, uniform in space and time
    (extrapolated time)
    &
    \num{4.50e6}
    &
    1.05
    &
    \num{4.50 e7}
    $\si{\sec}$
    \\
    \bottomrule
  \end{tabular}
\end{table}

In~\cref{tab:ad} we consider the average number of Newton iterations per time
step,
the assembling time for the Newton matrix
and
the right hand side of the linear equation
system
for one half cycle.
We compare the results for the derived linearization and the AD technique in
the 1D and the 2D simulation setup.
More numerical results for the latter case are given
in~\cref{subsubsec:2d_quarter_ellipse}.
Additionally, we compare for the AD 2D simulation case also these
simulations:
\begin{itemize}
    \item
    without MPI parallelization, but with space and time adaptivity
    \item
    without MPI parallelization and constant spatial discretization, but
    with time
    adaptivity
    \item
    without MPI parallelization and constant spatial and temporal
    discretization
\end{itemize}
For the second case,
we choose the spatial discretization for the whole duration of the simulation
to be just below the largest number of degrees of freedom~(DOFs) in the
adaptive simulation case.
For the last situation with uniform temporal
discretization, we choose \num{2e-7}~\si{\hour} and BDF with order one for the
time integration scheme.
\num{2e-7}~\si{\hour} is close, but just above the lowest time step size of the
space and
time adaptive simulation case after the initial time step size.
These values are necessary to correctly record all relevant physical phenomena.
Note that we count all Newton steps, also the Newton steps when the
rate for reduction of
the new Newton update is not fast enough or the tolerances of the adaptive
algorithm are not fulfilled.
Considering~\cref{fig:numerical_efficiency}(a),
the majority of time steps in one lithiation cycle (from $t=0$ to $t=0.9$)
are purely elastic and result in larger time step sizes.
Due to the variable-step, variable-order time integration scheme, the
maximal number of Newton iterations is limited.
If the maximal
number of Newton iterations is reached, the rate for the Newton update is
too small or the tolerances of the adaptive algorithm are not achieved,
the
time step size is reduced.
A smaller time step size leads to a lower number of Newton iterations.
These points explain the low average number of Newton steps for our
particular simulation setting.

The number of time steps and the average number of Newton
iterations per time step are similar for the 1D and 2D simulations,
respectively.
However,
there are significant differences in the assembling times:
the linearization is twice as slow as
AD for the 1D simulation and even ten times slower for the 2D
simulation.
This is a remarkable acceleration in the assembling times when solving for the
Newton
update.
Not using the MPI parallelization leads to assembling time roughly
four times slower, as expected.
Considering the constant grid, we see a deceleration of over a factor of two,
but
also a slight reduction of the number of time steps.
For the last case, we computed 100
time steps and extrapolated the total assembling times.
The average number of
Newton steps for these 100 steps is lower than in the other cases, also due to
the
fact that within these first 100 time steps no plastic deformation occurs.
Considering a
simulation
without MPI parallelization, but with constant spatial and temporal
discretization, would
lead to a around $\num{e5}$ higher assembling time.
In addition, there would be another factor of ten without using the AD
technique, so in total,
we have
an speed up
by a factor of
roughly estimated~$\num{e6}$.
All in
all,~\cref{tab:ad} shows the huge acceleration for the assembling during the
simulation with
modern
numerical
techniques.

\begin{figure}[t]
  \centering
  {\includegraphics[width=\textwidth,page=1]
    {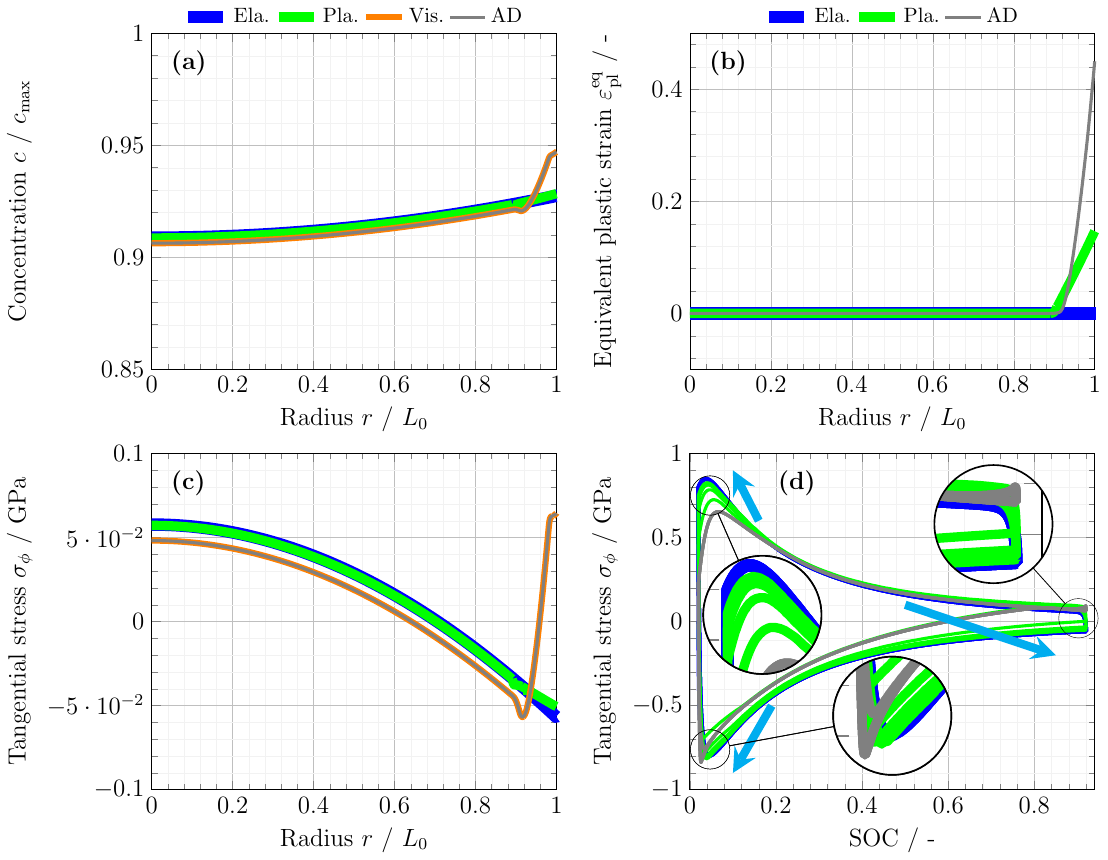}}
  \caption{
    Numerical results for the elastic (Ela.), plastic (Pla.),
    viscoplastic (Vis.) and viscoplastic with AD (AD) approaches
    of the 1D radial symmetric case at~$t =
    8.12~\si{h}$
    over the particle radius~$r$:
    concentration~$c$ in~(a),
    equivalent plastic strain~$\varepsilon_\pl^{\text{eq}}$
    in~(b) and
    tangential Cauchy
    stress~$\sigma_\phi$
    in~(c)
    as well as
    tangential Cauchy stress~$\sigma_\phi$
    over~$\text{SOC}$
    in~(d).
  }
  \label{fig:1d_nine_half_cycle}
\end{figure}

\textbf{Physical results after nine half cycles.}
Before proceeding to investigate the 2D computational domain, we use the
capabilities of the
spatially and temporally adaptive algorithm together with the AD technique, to
study multiple charging and discharging cycles.
To be more precise, we consider five charging
and four discharging cycles.
\cref{fig:1d_nine_half_cycle}
shows that the viscoplastic case with the derived
linearization or use of AD deliver identical results.
The differences between the plastic and viscoplastic deformation are
remarkable after nine half cycles, especially for the concentration
in~\cref{fig:1d_nine_half_cycle}(a).
The plastic approach is almost identical to the elastic results, whereas the
viscoplastic results are more similar to~\cref{fig:1d_one_half_cycle}.
The similarity between the elastic and plastic results can be explained by the
isotropic hardening, which leads to a large increase in the yield stress.
Thus, after nine half cycles, occurring stresses inside the particle do not
reach
the yield stress, preventing the particle to deform plastically any further.
This in turn leads to a more homogeneous distribution of stresses and thus no
pile
up of Li atoms close to the particle surface.
The lower magnifier reveals also for the viscoplastic approach some difference
for
higher cycling numbers:
a small wave occurs due to the change in the sign of the lithium flux for
charging
or discharging the active particle.
This is important to notice because the small wave is already more pronounced
for
the tangential Cauchy stress in~\cref{fig:1d_nine_half_cycle}(c)
and is a candidate for possible further inhomogeneities.
It should also be noted the higher magnitude in~\cref{fig:1d_nine_half_cycle}(b)
compared to~\cref{fig:1d_one_half_cycle}(b),
indicating large extent of plastic deformation which could eventually lead to
damage of the particle.
\cref{fig:1d_nine_half_cycle}(d) shows the development of the tangential
Cauchy stress at the particle surface over the~$\text{SOC}$.
For the elastic case there is no difference recognizable after the first half
cycle and all further half cycles are identical.
Also for the viscoplastic case,
here computed with AD,
only a small difference is visible at the lower left corner,
see the magnifier in the bottom center.
The second and all further lithiation cycles feature a higher compressive
stress,
even slightly higher than the elastic case,
than in the first half cycle.
This can be explained by the constant initial data.
The magnifying glass at the right top shows the effect of the change from
lithiation to delithiation with a short increase and followed decrease
indicating an additional plastic deformation.
A larger difference
between the elastic and viscoplastic approach at lower~$\text{SOC}$
is observable at the end of the discharging cycle at the
left top corner with the magnifier at the left center.
The plastic ansatz with isotropic hardening, however, shows the biggest
difference in each cycling.
As already discussed,
the tangential Cauchy stresses tend from cycle to cycle
to the values of the elastic case indicated with the cyan arrows.
It should be noted that the small wave of~\cref{fig:1d_nine_half_cycle}(a) and
(c)
is not visible in this representation.

\subsubsection{2D Quarter Ellipse}
\label{subsubsec:2d_quarter_ellipse}

\begin{figure}[t]
  \centering
  {\includegraphics[width=0.799\textwidth,page=1]
    {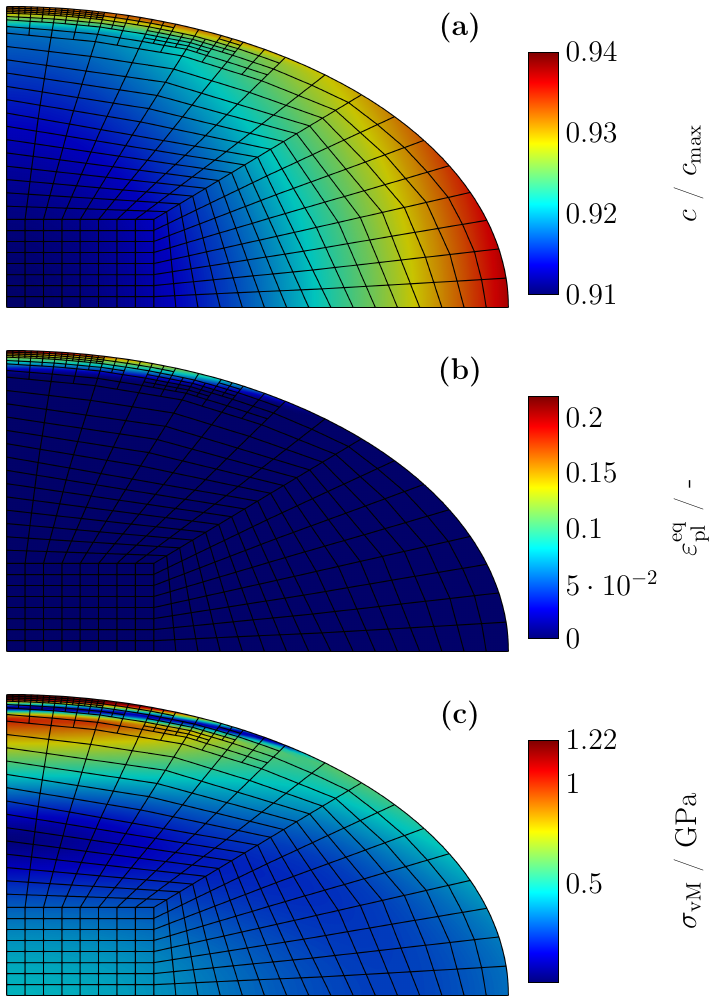}}
  \vspace{0.1cm}
  \caption{
    2D quarter ellipse
    after one lithiation at $\text{SOC}=0.92$
    for concentration~$c$
    in~(a),
    equivalent plastic strain~$\varepsilon_\pl^{\text{eq}}$
    in~(b)
    and von Mises stress~$\sigma_\text{vM}$
    in~(c) using AD techniques.}
  \label{fig:2d_quarter_ellipse}
\end{figure}
We now evaluate numerical results for the 2D quarter ellipse simulation
setup. Here, we adapt the initial time step size to~$\tau_0=\num{1e-8}$
and a minimal refinement level to three.
Further, we set~$\theta_{\mathrm{c}} = 0.005$
and~$\tau_{\max} = \num{5e-3}$.
With the capabilities of the adaptive space and time
solution algorithm as well as the presented AD technique, we consider now the
effects of half-axes of different lengths on plastic
deformation and
concentration distribution.

\textbf{Physical results after one half cycle.}
\cref{fig:2d_quarter_ellipse}
shows the numerical results at the end of one half cycle
with~$\tfinal=0.9~\si{\hour}$
and~$\text{SOC}=0.92$.
The number of
DOFs
are in a range
of~$[320, 121472]$
and the total computation time is less than 15 minutes.
In all subfigures, the adaptive mesh is visible in the background, indicating
higher refinement levels at the external surface near the smaller half-axis.
In~\cref{fig:2d_quarter_ellipse}(a),
the concentration profile is displayed in
a range of~$[0.91, 0.94]$.
It is noticeable that the gradient in the concentration is stronger at the
particle surface on the short half-axis.
Here, a larger area with lower concentration values in blue can be located with
a small but steeply growing area near the particle surface.
This property is comparable to the 1D simulation results with a higher
increase near the particle surface in~\cref{fig:1d_one_half_cycle}(a)
or in~\cref{fig:1d_nine_half_cycle}(a).

Having a look on the scalar equivalent plastic
strain~$\varepsilon_\pl^{\text{eq}}$
in~\cref{fig:2d_quarter_ellipse}(b)
the area of plastic deformation can be confirmed:
plastic deformation occurs at the particle surface of the smaller half-axis.
Large equivalent plastic strains of up to $22\,\%$ occur.
\cref{fig:2d_quarter_ellipse}(c) shows the von Mises stress in the general
plane state
\begin{align}
  \label{eq:von_mises_stress}
  \sigma_{\text{vM}} = \sqrt{\sigma_{11}^2 + \sigma_{22}^2 -
    \sigma_{11}\sigma_{22} + 3\sigma_{12}^2}.
\end{align}
A significant inhomogeneity of the von Mises stress distribution in the
particle is identifiable.
Whereas the stress values at the larger half-axis are lower than
$0.45~\si{GPa}$, the von Mises stresses are more than
1.5 times higher at the
smaller half-axis and feature two peaks near the particle surface:
one due to tensile stresses because of the plastic deformation next to the
surface and one due to
compressive stresses a little further inside of the particle, both again like
in the 1D case in~\cref{fig:1d_one_half_cycle}(c)
or in~\cref{fig:1d_nine_half_cycle}(c).
The strong change of stress gradients near the surface of the smaller
half-axis may appear unusual. However, as outlined,
the change of the gradients results from the change in sign of the tangential
stresses from tensile at the particle surface to compressive in the particle
interior close to the surface, c.f.~\cref{fig:1d_one_half_cycle}(c).
As the von Mises stress, c.f.~\cref{eq:von_mises_stress}, is the absolute
value of the stress deviator, this change results in the two maxima
close together.
This inhomogeneity of the stress distribution
between compressive and tensile stresses
is especially important for
further investigations on particle damage or fracture.

However, the finding of the higher stresses and plastic deformation at the
smaller half-axis are in contrast to the numerical outputs
in~\cref{fig:parameter_study}(b).
There, larger particle size results in larger stresses, so it would be
reasonable to feature higher stresses at the longer half-axis.
Though, the longer half-axis is responsible for lower concentration values in
the particle center.
On the smaller half-axis, the diffusion would create higher concentrations in
the particle center.
Thus, the influence of the lower concentration values at the longer half-axis
is responsible for higher concentration gradients at the smaller half-axis
resulting in higher stress and plastic deformation.
In this context it is important to recall that a constant external
lithium flux is used at the particle surface, which favors this behavior.




\section{Conclusion and Outlook}
\label{sec:conclusion}
We conclude our work with a summary and an outlook.
\subsection{Conclusion}
Within this work, a large deformation, chemo-elasto-plastic model to predict the
diffusion-deformation behavior of aSi anode particles within Li-Ion batteries
is presented. The used plasticity model is in close analogy to the model
presented in \cite{di-leo2015diffusion-deformation}. As an extension, we
specify our theory for both rate-dependent and rate-independent plastic
deformations, see~\cref{subsec:plasticity},
and compare numerical results of both
theories after multiple charging and discharging cycles.
The chemical model relies on \cite{schoof2023simulation}, where a
measured
experimental
OCV curve
is used to model the chemical contribution to the
Helmholtz
free energy.
The derived model is implemented in an efficient finite element scheme, first
presented in~\cite{castelli2021efficient}, relying on space and time adaptive
algorithms and parallelization \cite{schoof2022parallelization}
as well as modern numerical techniques like AD.

For both plastic models, a return mapping algorithm \cite{simo1986return} is
used in the context of static condensation, allowing the evaluation of the
plastic deformations at the finite element integration point level instead of
treating them as additional degrees of
freedom~\cite{kolzenberg2022chemo-mechanical}.
We present
in~\cref{subsubsec:weak_formulation}
a projection onto the set of admissible stresses inspired
by~\cite{frohne2016efficient}, which can be given explicitly in the case of
the rate-independent model with linear hardening,
whereas in the viscoplastic case the projector is
implicit, due to the non-linearity introduced by the
evolution of the equivalent plastic strains.
The linearization of these projectors, necessary for computing Newton updates
in the nonlinear solution scheme, is approximated
in~\cref{subsubsec:adaptive_algorithm}
inspired
by~\cite{di-leo2015diffusion-deformation}.
In addition,
we apply AD techniques
to circumvent the necessity of tedious analytic computations and assembling
operations and compare the numeric performance of both approaches.
We incorporate our DAE system into an existing efficient space and
time adaptive solver, presented in~\cref{sec:numerics}.

The numerical results
in~\cref{subsec:numerical_results}
for a given set of parameters show a heterogeneous
concentration distribution within the particle, when plastic models are
considered. We attribute this to lower stresses inside the particle,
which in turn affect the mobility of Li atoms and lead to a pile up at the
particle's surface. In our studies, plastic deformations are limited to the
outer
parts of the particle and result in an eigenstrain
that leads to tensile
stresses
at the particle surface during charging, which is not observed in
the elastic model. Both the concentration pile up and the plastic deformation
can be mitigated by decreasing the particle size and or decreasing the
charging rate.

Comparing the plastic and viscoplastic model
in~\cref{fig:1d_nine_half_cycle}, differences
become visible after multiple charging and discharging cycles.
On the one hand, the plastic
model hardens isotropic and after several cycles the occurring stresses do
not reach the increased yield limit, preventing further plastification.
On the other hand, the viscoplastic model does show an increased amount of
plastic
deformations, which results in a different stress and
concentration
distribution.
As a result of~\cref{sec:results},
the differences are more pronounced the more cycles are
considered,
which seems not to be investigated before and affects battery performance and
lifetime.
This could be achieved by the use of several modern numerical techniques.

Investigating the performance of the
approximation of the
projector's linearization in comparison
with the AD scheme
in~\cref{tab:ad},
we conclude that the AD scheme is way more efficient,
due to the fact that the assembling of the Newton matrix is done simultaneously
with the residual.
This leads to decreased computation times.
In addition, the use of AD circumvents
the tedious analytical derivation as well as the
fault prone implementation of the linearization.
In total, with the application of efficient tools like parallelization,
adaptive spatial and temporal discretization and AD on our higher order
finite element approach,
we
have reduced our
computational effort for the assembling
significantly.
Adding the plastic part of the
deformation gradient as solution variable to the
set
of equations like
in~\cite{kolzenberg2022chemo-mechanical} would increase the computational time
even more, especially for the purely two-dimensional computational
consideration.

When studying a two-dimensional problem
in~\cref{subsubsec:2d_quarter_ellipse},
an asymmetry in the
concentration and plastic strain distribution is observed, which we attribute
to the ellipsoidal particle shape and the resulting asymmetric concentration
distribution.

To conclude, our study indicates that small particle sizes with a spherical
shape result in a smaller build up of stresses and is therefore desirable.
Moreover, semi-axes with different lengths should be
avoided to obtain homogeneously distributed mechanical stresses.
Regarding material behavior, the considered isotropic hardening mechanism
is favorable, as less plastic strains build up after multiple charging and
discharging cycles.
In addition, this would prevent the sharp concentration gradients observed in
our study and thus less interference with battery performance is to be expected.
A further combination with the viscoplastic approach could be a possible
extension.

\subsection{Outlook}
To further extend the theory originating from this work, numerous paths can be
taken.
We consider the following ones to be especially interesting:
\begin{itemize}
  \item
  Due to the usage of the chemical potential as solution variable
  and the corresponding code
  structure, first introduced for the phase-separating material
  lithium iron phosphate in~\cite{castelli2021efficient},
  we can add a Cahn--Hilliard approach for the phase-separation
  phenomena of aSi which is observed
  by~\cite{uxa2019on}.
  \item Previous research pointed towards the usage of nanowires as electrode
  geometry.
  The reason lies in the observation
  that lower electrode sizes decrease
  occurring
  stresses and plastic deformations. This is in agreement with the results of
  decreasing
  stresses for smaller particle radii. Still, a study of various nanowire
  shaped
  electrode geometries is a reasonable application for the derived model and
  the proposed
  implementation scheme.
  \item Fracture and SEI formation typically occur in aSi anode
  particles~\cite{kolzenberg2022chemo-mechanical}.
  Both have been considered previously,
  see, e.g.~\cite{kolzenberg2022chemo-mechanical, di-leo2014cahn-hilliard-type,
  schomburg2023characterization},
  and should be included in future extensions of this work, to better capture
  the physical
  behavior of the anode particles.
  \item A paper published recently by several of the
  authors~\cite{schoof2023simulation},
  introduced an efficient scheme to compute particle deformation behavior when
  contact between particles or walls
  is introduced. A study of the elasto-plastic material behavior derived in
  this work, in
  combination with the obstacle, could provide further insight into the
  mechanics of
  anode materials during charging and discharging.
  \item
  An integration of the Butler--Volmer condition into the external
  boundary condition for the lithium flux is of further interest to consider
  additional surface boundary reactions.
  \item Finally, a more robust and scalable solver can be developed to
  provide additional speedup compared to the currently used LU-decomposition
  when using MPI parallelization~\cite{schoof2022parallelization}.
  Further consideration can also be given to the matrix-free concept of
  \texttt{deal.ii} \cite{arndt2021deal-ii, castelli2019efficient} and other
  solvers like PARDISO~\cite{schenk2011pardiso}.
\end{itemize}


  \section*{Declaration of competing interest}

  \noindent
  The authors declare that they have no known competing financial interests or
  personal relationships that could have appeared to influence the work in this
  paper.

  \section*{CrediT authorship contribution statement}
  \noindent
  \textbf{R. Schoof:}
  Conceptualization, Data curation, Formal Analysis, Investigation, Methodology,
  Software, Validation, Visualization, Writing -- original draft
  \noindent
  \textbf{J. Niermann:}
  Data curation, Formal Analysis, Methodology, Software, Validation, Writing --
  review \&
  editing
  \noindent
  \textbf{A. Dyck:}
  Conceptualization, Formal Analysis, Investigation, Methodology, Writing --
  original draft
  \noindent
  \textbf{T. B\"ohlke:}
  Funding acquisition, Project administration, Resources, Supervision, Writing
  --
  review \& editing
  \noindent
  \textbf{W. D\"orfler:}
  Funding acquisition, Project administration, Resources, Supervision, Writing
  --
  review \& editing

  \section*{Acknowledgment}

  \noindent
  The authors thank G.~F.~Castelli for the software basis and
  L.~von Kolzenberg,
  L.~K\"obbing
  and C. Wieners
  for intensive and constructive discussions about modeling
  silicon particles.
  R.S. acknowledges financial support by the German Research
  Foundation~(DFG) through the Research Training Group 2218
  SiMET~--~Simulation of Mechano-Electro-Thermal processes in Lithium-ion
  Batteries, project number 281041241.
  T.B. acknowledges partial support by the German Research Foundation (DFG)
  within the Priority Programme SPP2013 ’Targeted Use of Forming Induced
  Residual Stresses in Metal Components’ (Bo 1466/14-2). The support by the
  German Research Foundation (DFG) is gratefully acknowledged.

  \section*{ORCID}

  \noindent
  R. Schoof: \href{https://orcid.org/0000-0001-6848-3844}{0000-0001-6848-3844},
  J. Niermann \href{https://orcid.org/0009-0002-0422-0521}{0009-0002-0422-0521},
  A. Dyck \href{https://orcid.org/0000-0001-7239-9653}{0000-0001-7239-9653},
  T. B\"ohlke \href{https://orcid.org/0000-0001-6884-0530}{0000-0001-6884-0530},
  W. D\"orfler:
  \href{https://orcid.org/0000-0003-1558-9236}{0000-0003-1558-9236}


  \bibliography{SchoofEtAl2023Modeling}


\begin{thebibliography}{65}
\ifx \bisbn   \undefined \def \bisbn  #1{ISBN #1}\fi
\ifx \binits  \undefined \def \binits#1{#1}\fi
\ifx \bauthor  \undefined \def \bauthor#1{#1}\fi
\ifx \batitle  \undefined \def \batitle#1{#1}\fi
\ifx \bjtitle  \undefined \def \bjtitle#1{#1}\fi
\ifx \bvolume  \undefined \def \bvolume#1{\textbf{#1}}\fi
\ifx \byear  \undefined \def \byear#1{#1}\fi
\ifx \bissue  \undefined \def \bissue#1{#1}\fi
\ifx \bfpage  \undefined \def \bfpage#1{#1}\fi
\ifx \blpage  \undefined \def \blpage #1{#1}\fi
\ifx \burl  \undefined \def \burl#1{\textsf{#1}}\fi
\ifx \doiurl  \undefined \def \doiurl#1{\url{https://doi.org/#1}}\fi
\ifx \betal  \undefined \def \betal{\textit{et al.}}\fi
\ifx \binstitute  \undefined \def \binstitute#1{#1}\fi
\ifx \binstitutionaled  \undefined \def \binstitutionaled#1{#1}\fi
\ifx \bctitle  \undefined \def \bctitle#1{#1}\fi
\ifx \beditor  \undefined \def \beditor#1{#1}\fi
\ifx \bpublisher  \undefined \def \bpublisher#1{#1}\fi
\ifx \bbtitle  \undefined \def \bbtitle#1{#1}\fi
\ifx \bedition  \undefined \def \bedition#1{#1}\fi
\ifx \bseriesno  \undefined \def \bseriesno#1{#1}\fi
\ifx \blocation  \undefined \def \blocation#1{#1}\fi
\ifx \bsertitle  \undefined \def \bsertitle#1{#1}\fi
\ifx \bsnm \undefined \def \bsnm#1{#1}\fi
\ifx \bsuffix \undefined \def \bsuffix#1{#1}\fi
\ifx \bparticle \undefined \def \bparticle#1{#1}\fi
\ifx \barticle \undefined \def \barticle#1{#1}\fi
\bibcommenthead
\ifx \bconfdate \undefined \def \bconfdate #1{#1}\fi
\ifx \botherref \undefined \def \botherref #1{#1}\fi
\ifx \url \undefined \def \url#1{\textsf{#1}}\fi
\ifx \bchapter \undefined \def \bchapter#1{#1}\fi
\ifx \bbook \undefined \def \bbook#1{#1}\fi
\ifx \bcomment \undefined \def \bcomment#1{#1}\fi
\ifx \oauthor \undefined \def \oauthor#1{#1}\fi
\ifx \citeauthoryear \undefined \def \citeauthoryear#1{#1}\fi
\ifx \endbibitem  \undefined \def \endbibitem {}\fi
\ifx \bconflocation  \undefined \def \bconflocation#1{#1}\fi
\ifx \arxivurl  \undefined \def \arxivurl#1{\textsf{#1}}\fi
\csname PreBibitemsHook\endcsname

\bibitem[\protect\citeauthoryear{Zhao et~al.}{2019}]{zhao2019review}
\begin{barticle}
\bauthor{\bsnm{Zhao}, \binits{Y.}},
\bauthor{\bsnm{Stein}, \binits{P.}},
\bauthor{\bsnm{Bai}, \binits{Y.}},
\bauthor{\bsnm{Al-Siraj}, \binits{M.}},
\bauthor{\bsnm{Yang}, \binits{Y.}},
\bauthor{\bsnm{Xu}, \binits{B.-X.}}:
\batitle{A review on modeling of electro-chemo-mechanics in lithium-ion
  batteries}.
\bjtitle{J. Power Sources}
\bvolume{413},
\bfpage{259}--\blpage{283}
(\byear{2019})
\doiurl{10.1016/j.jpowsour.2018.12.011}
\end{barticle}
\endbibitem

\bibitem[\protect\citeauthoryear{Tomaszewska
  et~al.}{2019}]{tomaszewska2019lithium-ion}
\begin{barticle}
\bauthor{\bsnm{Tomaszewska}, \binits{A.}},
\bauthor{\bsnm{Chu}, \binits{Z.}},
\bauthor{\bsnm{Feng}, \binits{X.}},
\bauthor{\bsnm{O'Kane}, \binits{S.}},
\bauthor{\bsnm{Liu}, \binits{X.}},
\bauthor{\bsnm{Chen}, \binits{J.}},
\bauthor{\bsnm{Ji}, \binits{C.}},
\bauthor{\bsnm{Endler}, \binits{E.}},
\bauthor{\bsnm{Li}, \binits{R.}},
\bauthor{\bsnm{Liu}, \binits{L.}},
\bauthor{\bsnm{Li}, \binits{Y.}},
\bauthor{\bsnm{Zheng}, \binits{S.}},
\bauthor{\bsnm{Vetterlein}, \binits{S.}},
\bauthor{\bsnm{Gao}, \binits{M.}},
\bauthor{\bsnm{Du}, \binits{J.}},
\bauthor{\bsnm{Parkes}, \binits{M.}},
\bauthor{\bsnm{Ouyang}, \binits{M.}},
\bauthor{\bsnm{Marinescu}, \binits{M.}},
\bauthor{\bsnm{Offer}, \binits{G.}},
\bauthor{\bsnm{Wu}, \binits{B.}}:
\batitle{Lithium-ion battery fast charging: {A} review}.
\bjtitle{eTransportation}
\bvolume{1},
\bfpage{100011}
(\byear{2019})
\doiurl{10.1016/j.etran.2019.100011}
\end{barticle}
\endbibitem

\bibitem[\protect\citeauthoryear{{de Vasconcelos}
  et~al.}{2022}]{vasconcelos2022chemomechanics}
\begin{barticle}
\bauthor{\bsnm{{de Vasconcelos}}, \binits{L.S.}},
\bauthor{\bsnm{R.}, \binits{X.}},
\bauthor{\bsnm{Xu}, \binits{Z.}},
\bauthor{\bsnm{Zhang}, \binits{J.}},
\bauthor{\bsnm{Sharma}, \binits{N.}},
\bauthor{\bsnm{Shah}, \binits{S.R.}},
\bauthor{\bsnm{Han}, \binits{J.}},
\bauthor{\bsnm{He}, \binits{X.}},
\bauthor{\bsnm{Wu}, \binits{X.}},
\bauthor{\bsnm{Sun}, \binits{H.}},
\bauthor{\bsnm{Hu}, \binits{S.}},
\bauthor{\bsnm{Perrin}, \binits{M.}},
\bauthor{\bsnm{Wang}, \binits{X.}},
\bauthor{\bsnm{Liu}, \binits{Y.}},
\bauthor{\bsnm{Lin}, \binits{F.}},
\bauthor{\bsnm{Cui}, \binits{Y.}},
\bauthor{\bsnm{Zhao}, \binits{K.}}:
\batitle{Chemomechanics of rechargeable batteries: Status, theories, and
  perspectives}.
\bjtitle{Chem. Rev.}
\bvolume{122}(\bissue{15}),
\bfpage{13043}--\blpage{13107}
(\byear{2022})
\doiurl{10.1021/acs.chemrev.2c00002}
\end{barticle}
\endbibitem

\bibitem[\protect\citeauthoryear{Uxa et~al.}{2019}]{uxa2019on}
\begin{barticle}
\bauthor{\bsnm{Uxa}, \binits{D.}},
\bauthor{\bsnm{Jerliu}, \binits{B.}},
\bauthor{\bsnm{H\"{u}ger}, \binits{E.}},
\bauthor{\bsnm{D\"{o}rrer}, \binits{L.}},
\bauthor{\bsnm{Horisberger}, \binits{M.}},
\bauthor{\bsnm{Stahn}, \binits{J.}},
\bauthor{\bsnm{Schmidt}, \binits{H.}}:
\batitle{On the lithiation mechanism of amorphous silicon electrodes in li-ion
  batteries}.
\bjtitle{J. Phys. Chem. C}
\bvolume{123}(\bissue{36}),
\bfpage{22027}--\blpage{22039}
(\byear{2019})
\doiurl{10.1021/acs.jpcc.9b06011}
\end{barticle}
\endbibitem

\bibitem[\protect\citeauthoryear{Zhang}{2011}]{zhang2011review}
\begin{barticle}
\bauthor{\bsnm{Zhang}, \binits{W.-J.}}:
\batitle{A review of the electrochemical performance of alloy anodes for
  lithium-ion batteries}.
\bjtitle{J. Power Sources}
\bvolume{196}(\bissue{1}),
\bfpage{13}--\blpage{24}
(\byear{2011})
\doiurl{10.1016/j.jpowsour.2010.07.020}
\end{barticle}
\endbibitem

\bibitem[\protect\citeauthoryear{Zhao et~al.}{2011}]{zhao2011large}
\begin{barticle}
\bauthor{\bsnm{Zhao}, \binits{K.}},
\bauthor{\bsnm{Pharr}, \binits{M.}},
\bauthor{\bsnm{Cai}, \binits{S.}},
\bauthor{\bsnm{Vlassak}, \binits{J.J.}},
\bauthor{\bsnm{Suo}, \binits{Z.}}:
\batitle{Large plastic deformation in high-capacity lithium-ion batteries
  caused by charge and discharge}.
\bjtitle{J. Am. Ceram. Soc.}
\bvolume{94},
\bfpage{226}--\blpage{235}
(\byear{2011})
\doiurl{10.1111/j.1551-2916.2011.04432.x}
\end{barticle}
\endbibitem

\bibitem[\protect\citeauthoryear{Di~Leo
  et~al.}{2015}]{di-leo2015diffusion-deformation}
\begin{barticle}
\bauthor{\bsnm{Di~Leo}, \binits{C.V.}},
\bauthor{\bsnm{Rejovitzky}, \binits{E.}},
\bauthor{\bsnm{Anand}, \binits{L.}}:
\batitle{Diffusion-deformation theory for amorphous silicon anodes: The role of
  plastic deformation on electrochemical performance}.
\bjtitle{Int. J. Solids Struct.}
\bvolume{67-68},
\bfpage{283}--\blpage{296}
(\byear{2015})
\doiurl{10.1016/j.ijsolstr.2015.04.028}
\end{barticle}
\endbibitem

\bibitem[\protect\citeauthoryear{Poluektov
  et~al.}{2018}]{poluektov2018modelling}
\begin{barticle}
\bauthor{\bsnm{Poluektov}, \binits{M.}},
\bauthor{\bsnm{Freidin}, \binits{A.B.}},
\bauthor{\bsnm{Figiel}, \binits{L.}}:
\batitle{Modelling stress-affected chemical reactions in non-linear
  viscoelastic solids with application to lithiation reaction in spherical {S}i
  particles}.
\bjtitle{Internat. J. Engrg. Sci.}
\bvolume{128},
\bfpage{44}--\blpage{62}
(\byear{2018})
\doiurl{10.1016/j.ijengsci.2018.03.007}
\end{barticle}
\endbibitem

\bibitem[\protect\citeauthoryear{Vadhva et~al.}{2022}]{vadhva2022towards}
\begin{barticle}
\bauthor{\bsnm{Vadhva}, \binits{P.}},
\bauthor{\bsnm{Boyce}, \binits{A.M.}},
\bauthor{\bsnm{Hales}, \binits{A.}},
\bauthor{\bsnm{Pang}, \binits{M.-C.}},
\bauthor{\bsnm{Patel}, \binits{A.N.}},
\bauthor{\bsnm{Shearing}, \binits{P.R.}},
\bauthor{\bsnm{Offer}, \binits{G.}},
\bauthor{\bsnm{Rettie}, \binits{A.J.E.}}:
\batitle{Towards optimised cell design of thin film silicon-based solid-state
  batteries via modelling and experimental characterisation}.
\bjtitle{J. Electrochem. Soc.}
\bvolume{169}(\bissue{10}),
\bfpage{100525}
(\byear{2022})
\doiurl{10.1149/1945-7111/ac9552}
\end{barticle}
\endbibitem

\bibitem[\protect\citeauthoryear{Zhao et~al.}{2011}]{zhao2011lithium-assisted}
\begin{barticle}
\bauthor{\bsnm{Zhao}, \binits{K.}},
\bauthor{\bsnm{Wang}, \binits{W.L.}},
\bauthor{\bsnm{Gregoire}, \binits{J.}},
\bauthor{\bsnm{Pharr}, \binits{M.}},
\bauthor{\bsnm{Suo}, \binits{Z.}},
\bauthor{\bsnm{Vlassak}, \binits{J.J.}},
\bauthor{\bsnm{Kaxiras}, \binits{E.}}:
\batitle{Lithium-assisted plastic deformation of silicon electrodes in
  lithium-ion batteries: A first-principles theoretical study}.
\bjtitle{Nano Lett.}
\bvolume{11}(\bissue{7}),
\bfpage{2962}--\blpage{2967}
(\byear{2011})
\doiurl{10.1021/nl201501s}
\end{barticle}
\endbibitem

\bibitem[\protect\citeauthoryear{Castelli et~al.}{2021}]{castelli2021efficient}
\begin{barticle}
\bauthor{\bsnm{Castelli}, \binits{G.F.}},
\bauthor{\bsnm{{von Kolzenberg}}, \binits{L.}},
\bauthor{\bsnm{Horstmann}, \binits{B.}},
\bauthor{\bsnm{Latz}, \binits{A.}},
\bauthor{\bsnm{D{\"o}rfler}, \binits{W.}}:
\batitle{Efficient simulation of chemical-mechanical coupling in battery active
  particles}.
\bjtitle{Energy Technol.}
\bvolume{9}(\bissue{6}),
\bfpage{2000835}
(\byear{2021})
\doiurl{10.1002/ente.202000835}
\end{barticle}
\endbibitem

\bibitem[\protect\citeauthoryear{Schoof et~al.}{2023}]{schoof2023simulation}
\begin{barticle}
\bauthor{\bsnm{Schoof}, \binits{R.}},
\bauthor{\bsnm{Castelli}, \binits{G.F.}},
\bauthor{\bsnm{D\"orfler}, \binits{W.}}:
\batitle{Simulation of the deformation for cycling chemo-mechanically coupled
  battery active particles with mechanical constraints}.
\bjtitle{Comput. Math. Appl.}
\bvolume{149},
\bfpage{135}--\blpage{149}
(\byear{2023})
\doiurl{10.1016/j.camwa.2023.08.027}
\end{barticle}
\endbibitem

\bibitem[\protect\citeauthoryear{Bertram}{2021}]{bertram2021elasticity}
\begin{bbook}
\bauthor{\bsnm{Bertram}, \binits{A.}}:
\bbtitle{Elasticity and Plasticity of Large Deformations: Including Gradient
  Materials},
\bedition{4}th edn.
\bpublisher{Springer},
\blocation{Cham,}
(\byear{2021}).
\doiurl{10.1007/978-3-030-72328-6}
\end{bbook}
\endbibitem

\bibitem[\protect\citeauthoryear{Lubliner}{2006}]{lubliner2006plasticity}
\begin{bbook}
\bauthor{\bsnm{Lubliner}, \binits{J.}}:
\bbtitle{Plasticity Theory}.
\bpublisher{Pearson Education, Inc.},
\blocation{New York}
(\byear{2006}).
\doiurl{10.1115/1.2899459}
\end{bbook}
\endbibitem

\bibitem[\protect\citeauthoryear{{von Kolzenberg}
  et~al.}{2022}]{kolzenberg2022chemo-mechanical}
\begin{barticle}
\bauthor{\bsnm{{von Kolzenberg}}, \binits{L.}},
\bauthor{\bsnm{Latz}, \binits{A.}},
\bauthor{\bsnm{Horstmann}, \binits{B.}}:
\batitle{Chemo-mechanical model of {SEI} growth on silicon electrode
  particles}.
\bjtitle{Batter. Supercaps}
\bvolume{5}(\bissue{2}),
\bfpage{202100216}
(\byear{2022})
\doiurl{10.1002/batt.202100216}
\end{barticle}
\endbibitem

\bibitem[\protect\citeauthoryear{Gritton et~al.}{2017}]{gritton2017using}
\begin{barticle}
\bauthor{\bsnm{Gritton}, \binits{C.}},
\bauthor{\bsnm{Guilkey}, \binits{J.}},
\bauthor{\bsnm{Hooper}, \binits{J.}},
\bauthor{\bsnm{Bedrov}, \binits{D.}},
\bauthor{\bsnm{Kirby}, \binits{R.M.}},
\bauthor{\bsnm{Berzins}, \binits{M.}}:
\batitle{Using the material point method to model chemical/mechanical coupling
  in the deformation of a silicon anode}.
\bjtitle{Modelling and Simulation in Materials Science and Engineering}
\bvolume{25}(\bissue{4}),
\bfpage{045005}
(\byear{2017})
\doiurl{10.1088/1361-651x/aa6830}
\end{barticle}
\endbibitem

\bibitem[\protect\citeauthoryear{Basu et~al.}{2019}]{basu2019structural}
\begin{barticle}
\bauthor{\bsnm{Basu}, \binits{S.}},
\bauthor{\bsnm{Koratkar}, \binits{N.}},
\bauthor{\bsnm{Shi}, \binits{Y.}}:
\batitle{Structural transformation and embrittlement during lithiation and
  delithiation cycles in an amorphous silicon electrode}.
\bjtitle{Acta Mater.}
\bvolume{175},
\bfpage{11}--\blpage{20}
(\byear{2019})
\doiurl{10.1016/j.actamat.2019.05.055}
\end{barticle}
\endbibitem

\bibitem[\protect\citeauthoryear{Pharr et~al.}{2014}]{pharr2014variation}
\begin{barticle}
\bauthor{\bsnm{Pharr}, \binits{M.}},
\bauthor{\bsnm{Suo}, \binits{Z.}},
\bauthor{\bsnm{Vlassak}, \binits{J.J.}}:
\batitle{Variation of stress with charging rate due to strain-rate sensitivity
  of silicon electrodes of li-ion batteries}.
\bjtitle{J. Power Sources}
\bvolume{270},
\bfpage{569}--\blpage{575}
(\byear{2014})
\doiurl{10.1016/j.jpowsour.2014.07.153}
\end{barticle}
\endbibitem

\bibitem[\protect\citeauthoryear{Sitinamaluwa
  et~al.}{2017}]{sitinamaluwa2017deformation}
\begin{barticle}
\bauthor{\bsnm{Sitinamaluwa}, \binits{H.}},
\bauthor{\bsnm{Nerkar}, \binits{J.}},
\bauthor{\bsnm{Wang}, \binits{M.}},
\bauthor{\bsnm{Zhang}, \binits{S.}},
\bauthor{\bsnm{Yan}, \binits{C.}}:
\batitle{Deformation and failure mechanisms of electrochemically lithiated
  silicon thin films}.
\bjtitle{{RSC} Adv.}
\bvolume{7}(\bissue{22}),
\bfpage{13487}--\blpage{13497}
(\byear{2017})
\doiurl{10.1039/c7ra01399j}
\end{barticle}
\endbibitem

\bibitem[\protect\citeauthoryear{Di~Leo}{2015}]{di-leo2015chemo-mechanics}
\begin{botherref}
\oauthor{\bsnm{Di~Leo}, \binits{C.V.}}:
Chemo-mechanics of lithium-ion battery electrodes.
PhD thesis,
Massachusetts Institute of Technology (MIT)
(2015)
\end{botherref}
\endbibitem

\bibitem[\protect\citeauthoryear{Schoof
  et~al.}{2022}]{schoof2022parallelization}
\begin{bchapter}
\bauthor{\bsnm{Schoof}, \binits{R.}},
\bauthor{\bsnm{Castelli}, \binits{G.F.}},
\bauthor{\bsnm{D\"orfler}, \binits{W.}}:
\bctitle{Parallelization of a finite element solver for chemo-mechanical
  coupled anode and cathode particles in lithium-ion batteries}.
In: \beditor{\bsnm{Kvamsdal}, \binits{T.}},
\beditor{\bsnm{Mathisen}, \binits{K.M.}},
\beditor{\bsnm{Lie}, \binits{K.-A.}},
\beditor{\bsnm{Larson}, \binits{M.G.}} (eds.)
\bbtitle{8th European Congress on Computational Methods in Applied Sciences and
  Engineering ({ECCOMAS} Congress 2022)}.
\bpublisher{{CIMNE}},
\blocation{Barcelona}
(\byear{2022}).
\doiurl{10.23967/eccomas.2022.106}
\end{bchapter}
\endbibitem

\bibitem[\protect\citeauthoryear{Wilson}{1974}]{wilson1974static}
\begin{barticle}
\bauthor{\bsnm{Wilson}, \binits{E.L.}}:
\batitle{The static condensation algorithm}.
\bjtitle{International Journal for Numerical Methods in Engineering}
\bvolume{8}(\bissue{1}),
\bfpage{198}--\blpage{203}
(\byear{1974})
\doiurl{10.1002/nme.1620080115}
\end{barticle}
\endbibitem

\bibitem[\protect\citeauthoryear{Di~Pietro
  et~al.}{2014}]{di-pietro2014arbitrary-order}
\begin{barticle}
\bauthor{\bsnm{Di~Pietro}, \binits{D.A.}},
\bauthor{\bsnm{Ern}, \binits{A.}},
\bauthor{\bsnm{Lemaire}, \binits{S.}}:
\batitle{An arbitrary-order and compact-stencil discretization of diffusion on
  general meshes based on local reconstruction operators}.
\bjtitle{Comput. Methods Appl. Math.}
\bvolume{14}(\bissue{4}),
\bfpage{461}--\blpage{472}
(\byear{2014})
\doiurl{10.1515/cmam-2014-0018}
\end{barticle}
\endbibitem

\bibitem[\protect\citeauthoryear{Di~Pietro and Ern}{2015}]{di-pietro2015hybrid}
\begin{barticle}
\bauthor{\bsnm{Di~Pietro}, \binits{D.A.}},
\bauthor{\bsnm{Ern}, \binits{A.}}:
\batitle{A hybrid high-order locking-free method for linear elasticity on
  general meshes}.
\bjtitle{Comput. Methods Appl. Mech. Engrg.}
\bvolume{283},
\bfpage{1}--\blpage{21}
(\byear{2015})
\doiurl{10.1016/j.cma.2014.09.009}
\end{barticle}
\endbibitem

\bibitem[\protect\citeauthoryear{O'Day and
  Curtin}{2005}]{oday2005superposition}
\begin{barticle}
\bauthor{\bsnm{O'Day}, \binits{M.P.}},
\bauthor{\bsnm{Curtin}, \binits{W.A.}}:
\batitle{A superposition framework for discrete dislocation plasticity}.
\bjtitle{J. Appl. Mech.}
\bvolume{71}(\bissue{6}),
\bfpage{805}--\blpage{815}
(\byear{2005})
\doiurl{10.1115/1.1794167}
\end{barticle}
\endbibitem

\bibitem[\protect\citeauthoryear{Frohne et~al.}{2016}]{frohne2016efficient}
\begin{barticle}
\bauthor{\bsnm{Frohne}, \binits{J.}},
\bauthor{\bsnm{Heister}, \binits{T.}},
\bauthor{\bsnm{Bangerth}, \binits{W.}}:
\batitle{Efficient numerical methods for the large-scale, parallel solution of
  elastoplastic contact problems}.
\bjtitle{Internat. J. Numer. Methods Engrg.}
\bvolume{105}(\bissue{6}),
\bfpage{416}--\blpage{439}
(\byear{2016})
\doiurl{10.1002/nme.4977}
\end{barticle}
\endbibitem

\bibitem[\protect\citeauthoryear{Barrera et~al.}{2016}]{barrera2016modelling}
\begin{barticle}
\bauthor{\bsnm{Barrera}, \binits{O.}},
\bauthor{\bsnm{Tarleton}, \binits{E.}},
\bauthor{\bsnm{Tang}, \binits{H.}},
\bauthor{\bsnm{Cocks}, \binits{A.}}:
\batitle{Modelling the coupling between hydrogen diffusion and the mechanical
  behaviour of metals}.
\bjtitle{Computational Materials Science}
\bvolume{122},
\bfpage{219}--\blpage{228}
(\byear{2016})
\end{barticle}
\endbibitem

\bibitem[\protect\citeauthoryear{Di~Leo
  et~al.}{2014}]{di-leo2014cahn-hilliard-type}
\begin{barticle}
\bauthor{\bsnm{Di~Leo}, \binits{C.V.}},
\bauthor{\bsnm{Rejovitzky}, \binits{E.}},
\bauthor{\bsnm{Anand}, \binits{L.}}:
\batitle{A {C}ahn--{H}illiard-type phase-field theory for species diffusion
  coupled with large elastic deformations: {Application} to phase-separating
  {Li}-ion electrode materials}.
\bjtitle{J. Mech. Phys. Solids}
\bvolume{70},
\bfpage{1}--\blpage{29}
(\byear{2014})
\doiurl{10.1016/j.jmps.2014.05.001}
\end{barticle}
\endbibitem

\bibitem[\protect\citeauthoryear{Simo and Hughes}{1998}]{simo1998computational}
\begin{bbook}
\bauthor{\bsnm{Simo}, \binits{J.C.}},
\bauthor{\bsnm{Hughes}, \binits{T.J.R.}}:
\bbtitle{Computational Inelasticity}.
\bsertitle{Interdisciplinary applied mathematics}.
\bpublisher{Springer},
\blocation{New York}
(\byear{1998})
\end{bbook}
\endbibitem

\bibitem[\protect\citeauthoryear{Simo and Taylor}{1985}]{simo1985consistent}
\begin{barticle}
\bauthor{\bsnm{Simo}, \binits{J.C.}},
\bauthor{\bsnm{Taylor}, \binits{R.L.}}:
\batitle{Consistent tangent operators for rate-independent elastoplasticity}.
\bjtitle{Comput. Methods Appl. Mech. Eng.}
\bvolume{48}(\bissue{1}),
\bfpage{101}--\blpage{118}
(\byear{1985})
\doiurl{10.1016/0045-7825(85)90070-2}
\end{barticle}
\endbibitem

\bibitem[\protect\citeauthoryear{Simo and Taylor}{1986}]{simo1986return}
\begin{barticle}
\bauthor{\bsnm{Simo}, \binits{J.C.}},
\bauthor{\bsnm{Taylor}, \binits{R.L.}}:
\batitle{A return mapping algorithm for plane stress elastoplasticity}.
\bjtitle{Internat. J. Numer. Methods Engrg.}
\bvolume{22}(\bissue{3}),
\bfpage{649}--\blpage{670}
(\byear{1986})
\doiurl{10.1002/nme.1620220310}
\end{barticle}
\endbibitem

\bibitem[\protect\citeauthoryear{Weber and Anand}{1990}]{weber1990finite}
\begin{barticle}
\bauthor{\bsnm{Weber}, \binits{G.}},
\bauthor{\bsnm{Anand}, \binits{L.}}:
\batitle{Finite deformation constitutive equations and a time integration
  procedure for isotropic, hyperelastic-viscoplastic solids}.
\bjtitle{Comput. Methods Appl. Mech. Eng.}
\bvolume{79}(\bissue{2}),
\bfpage{173}--\blpage{202}
(\byear{1990})
\doiurl{10.1016/0045-7825(90)90131-5}
\end{barticle}
\endbibitem

\bibitem[\protect\citeauthoryear{Simo}{1992}]{simo1992algorithms}
\begin{barticle}
\bauthor{\bsnm{Simo}, \binits{J.C.}}:
\batitle{Algorithms for static and dynamic multiplicative plasticity that
  preserve the classical return mapping schemes of the infinitesimal theory}.
\bjtitle{Comput. Methods Appl. Mech. Engrg.}
\bvolume{99}(\bissue{1}),
\bfpage{61}--\blpage{112}
(\byear{1992})
\doiurl{10.1016/0045-7825(92)90123-2}
\end{barticle}
\endbibitem

\bibitem[\protect\citeauthoryear{McDowell et~al.}{2013}]{mcdowell2013in}
\begin{barticle}
\bauthor{\bsnm{McDowell}, \binits{M.T.}},
\bauthor{\bsnm{Lee}, \binits{S.W.}},
\bauthor{\bsnm{Harris}, \binits{J.T.}},
\bauthor{\bsnm{Korgel}, \binits{B.A.}},
\bauthor{\bsnm{Wang}, \binits{C.}},
\bauthor{\bsnm{Nix}, \binits{W.D.}},
\bauthor{\bsnm{Cui}, \binits{Y.}}:
\batitle{In situ tem of two-phase lithiation of amorphous silicon nanospheres}.
\bjtitle{Nano Lett.}
\bvolume{13}(\bissue{2}),
\bfpage{758}--\blpage{764}
(\byear{2013})
\doiurl{10.1021/nl3044508}
\end{barticle}
\endbibitem

\bibitem[\protect\citeauthoryear{Castelli}{2021}]{castelli2021numerical}
\begin{botherref}
\oauthor{\bsnm{Castelli}, \binits{G.F.}}:
Numerical investigation of {C}ahn--{H}illiard-type phase-field models for
  battery active particles.
PhD thesis,
Karlsruhe Institute of Technology (KIT)
(2021).
\doiurl{10.5445/IR/1000141249}
\end{botherref}
\endbibitem

\bibitem[\protect\citeauthoryear{Holzapfel}{2010}]{holzapfel2010nonlinear}
\begin{bbook}
\bauthor{\bsnm{Holzapfel}, \binits{G.A.}}:
\bbtitle{Nonlinear Solid Mechanics}.
\bpublisher{John Wiley \& Sons, Ltd.},
\blocation{Chichester}
(\byear{2010})
\end{bbook}
\endbibitem

\bibitem[\protect\citeauthoryear{Braess}{2007}]{braess2007finite}
\begin{bbook}
\bauthor{\bsnm{Braess}, \binits{D.}}:
\bbtitle{Finite Elements},
\bedition{3}rd edn.
\bpublisher{Cambridge University Press},
\blocation{Cambridge}
(\byear{2007}).
\doiurl{10.1007/978-3-540-72450-6}
\end{bbook}
\endbibitem

\bibitem[\protect\citeauthoryear{Neff and Wieners}{2003}]{neff2003comparison}
\begin{barticle}
\bauthor{\bsnm{Neff}, \binits{P.}},
\bauthor{\bsnm{Wieners}, \binits{C.}}:
\batitle{Comparison of models for finite plasticity: a numerical study}.
\bjtitle{Comput. Vis. Sci.}
\bvolume{6}(\bissue{1}),
\bfpage{23}--\blpage{35}
(\byear{2003})
\doiurl{10.1007/s00791-003-0104-1}
\end{barticle}
\endbibitem

\bibitem[\protect\citeauthoryear{Latz and Zausch}{2015}]{latz2015multiscale}
\begin{barticle}
\bauthor{\bsnm{Latz}, \binits{A.}},
\bauthor{\bsnm{Zausch}, \binits{J.}}:
\batitle{Multiscale modeling of lithium ion batteries: thermal aspects}.
\bjtitle{Beilstein J. Nanotechnol.}
\bvolume{6},
\bfpage{987}--\blpage{1007}
(\byear{2015})
\doiurl{10.3762/bjnano.6.102}
\end{barticle}
\endbibitem

\bibitem[\protect\citeauthoryear{Latz and Zausch}{2011}]{latz2011thermodynamic}
\begin{barticle}
\bauthor{\bsnm{Latz}, \binits{A.}},
\bauthor{\bsnm{Zausch}, \binits{J.}}:
\batitle{Thermodynamic consistent transport theory of {Li}-ion batteries}.
\bjtitle{J. Power Sources}
\bvolume{196}(\bissue{6}),
\bfpage{3296}--\blpage{3302}
(\byear{2011})
\doiurl{10.1016/j.jpowsour.2010.11.088}
\end{barticle}
\endbibitem

\bibitem[\protect\citeauthoryear{Schammer et~al.}{2021}]{schammer2021theory}
\begin{barticle}
\bauthor{\bsnm{Schammer}, \binits{M.}},
\bauthor{\bsnm{Horstmann}, \binits{B.}},
\bauthor{\bsnm{Latz}, \binits{A.}}:
\batitle{Theory of transport in highly concentrated electrolytes}.
\bjtitle{J. Electrochem. Soc.}
\bvolume{168}(\bissue{2}),
\bfpage{026511}
(\byear{2021})
\doiurl{10.1149/1945-7111/abdddf}
\end{barticle}
\endbibitem

\bibitem[\protect\citeauthoryear{Deng}{2015}]{deng2015li-ion}
\begin{barticle}
\bauthor{\bsnm{Deng}, \binits{D.}}:
\batitle{Li-ion batteries: basics, progress, and challenges}.
\bjtitle{Energy Sci. Eng.}
\bvolume{3}(\bissue{5}),
\bfpage{385}--\blpage{418}
(\byear{2015})
\doiurl{10.1002/ese3.95}
\end{barticle}
\endbibitem

\bibitem[\protect\citeauthoryear{Kossa and Szab{\'o}}{2009}]{kossa2009exact}
\begin{barticle}
\bauthor{\bsnm{Kossa}, \binits{A.}},
\bauthor{\bsnm{Szab{\'o}}, \binits{L.}}:
\batitle{Exact integration of the von mises elastoplasticity model with
  combined linear isotropic-kinematic hardening}.
\bjtitle{Int. J. Plast.}
\bvolume{25}(\bissue{6}),
\bfpage{1083}--\blpage{1106}
(\byear{2009})
\doiurl{10.1016/j.ijplas.2008.08.003}
\end{barticle}
\endbibitem

\bibitem[\protect\citeauthoryear{Han et~al.}{1995}]{han1995computational}
\begin{barticle}
\bauthor{\bsnm{Han}, \binits{W.}},
\bauthor{\bsnm{Reddy}, \binits{B.D.}},
\bauthor{\bsnm{Oden}, \binits{J.T.}}:
\batitle{Computational plasticity: the variational basis and numerical
  analysis}.
\bjtitle{Comp. Mech. Advances}
\bvolume{2},
\bfpage{283}--\blpage{400}
(\byear{1995})
\end{barticle}
\endbibitem

\bibitem[\protect\citeauthoryear{Lubliner}{1986}]{lubliner1986normality}
\begin{barticle}
\bauthor{\bsnm{Lubliner}, \binits{J.}}:
\batitle{Normality rules in large-deformation plasticity}.
\bjtitle{Mechanics of Materials}
\bvolume{5}(\bissue{1}),
\bfpage{29}--\blpage{34}
(\byear{1986})
\end{barticle}
\endbibitem

\bibitem[\protect\citeauthoryear{Gro{\ss}mann}{2007}]{gromann2007numerical}
\begin{bbook}
\bauthor{\bsnm{Gro{\ss}mann}, \binits{C.V.}}:
\bbtitle{Numerical Treatment of Partial Differential Equations}.
\bsertitle{Universitext}.
\bpublisher{Springer},
\blocation{Berlin}
(\byear{2007})
\end{bbook}
\endbibitem

\bibitem[\protect\citeauthoryear{Wriggers}{2008}]{wriggers2008nonlinear}
\begin{bbook}
\bauthor{\bsnm{Wriggers}, \binits{P.}}:
\bbtitle{Nonlinear Finite Element Methods}.
\bpublisher{Springer},
\blocation{Berlin}
(\byear{2008})
\end{bbook}
\endbibitem

\bibitem[\protect\citeauthoryear{Suttmeier}{2010}]{suttmeier2010on}
\begin{barticle}
\bauthor{\bsnm{Suttmeier}, \binits{F.-T.}}:
\batitle{On plasticity with hardening: an adaptive finite element
  discretisation}.
\bjtitle{Int. Math. Forum}
\bvolume{5}(\bissue{49-52}),
\bfpage{2591}--\blpage{2601}
(\byear{2010})
\end{barticle}
\endbibitem

\bibitem[\protect\citeauthoryear{Reichelt et~al.}{1997}]{reichelt1997matlab}
\begin{botherref}
\oauthor{\bsnm{Reichelt}, \binits{M.W.}},
\oauthor{\bsnm{Shampine}, \binits{L.F.}},
\oauthor{\bsnm{Kierzenka}, \binits{J.}}:
MATLAB \texttt{ode15s}.
Copyright 1984--2020 {The MathWorks, Inc.}
(1997).
\url{https://www.mathworks.com}
\end{botherref}
\endbibitem

\bibitem[\protect\citeauthoryear{Shampine and
  Reichelt}{1997}]{shampine1997matlab}
\begin{barticle}
\bauthor{\bsnm{Shampine}, \binits{L.F.}},
\bauthor{\bsnm{Reichelt}, \binits{M.W.}}:
\batitle{The {MATLAB} {ODE} suite}.
\bjtitle{SIAM J. Sci. Comput.}
\bvolume{18}(\bissue{1}),
\bfpage{1}--\blpage{22}
(\byear{1997})
\doiurl{10.1137/S1064827594276424}
\end{barticle}
\endbibitem

\bibitem[\protect\citeauthoryear{Shampine et~al.}{1999}]{shampine1999solving}
\begin{barticle}
\bauthor{\bsnm{Shampine}, \binits{L.F.}},
\bauthor{\bsnm{Reichelt}, \binits{M.W.}},
\bauthor{\bsnm{Kierzenka}, \binits{J.A.}}:
\batitle{Solving index-{$1$} {DAE}s in {MATLAB} and {S}imulink}.
\bjtitle{SIAM Rev.}
\bvolume{41}(\bissue{3}),
\bfpage{538}--\blpage{552}
(\byear{1999})
\doiurl{10.1137/S003614459933425X}
\end{barticle}
\endbibitem

\bibitem[\protect\citeauthoryear{Shampine et~al.}{2003}]{shampine2003solving}
\begin{bbook}
\bauthor{\bsnm{Shampine}, \binits{L.F.}},
\bauthor{\bsnm{Gladwell}, \binits{I.}},
\bauthor{\bsnm{Thompson}, \binits{S.}}:
\bbtitle{Solving {ODE}s with {MATLAB}}.
\bpublisher{Cambridge University Press},
\blocation{Cambridge}
(\byear{2003}).
\doiurl{10.1017/CBO9780511615542}
\end{bbook}
\endbibitem

\bibitem[\protect\citeauthoryear{Hinterm\"{u}ller
  et~al.}{2002}]{hintermuller2002primal-dual}
\begin{barticle}
\bauthor{\bsnm{Hinterm\"{u}ller}, \binits{M.}},
\bauthor{\bsnm{Ito}, \binits{K.}},
\bauthor{\bsnm{Kunisch}, \binits{K.}}:
\batitle{The primal-dual active set strategy as a semismooth {N}ewton method}.
\bjtitle{SIAM J. Optim.}
\bvolume{13}(\bissue{3}),
\bfpage{865}--\blpage{8882003}
(\byear{2002})
\doiurl{10.1137/S1052623401383558}
\end{barticle}
\endbibitem

\bibitem[\protect\citeauthoryear{Arndt et~al.}{2021}]{arndt2021deal-ii}
\begin{barticle}
\bauthor{\bsnm{Arndt}, \binits{D.}},
\bauthor{\bsnm{Bangerth}, \binits{W.}},
\bauthor{\bsnm{Blais}, \binits{B.}},
\bauthor{\bsnm{Fehling}, \binits{M.}},
\bauthor{\bsnm{Gassm{\"o}ller}, \binits{R.}},
\bauthor{\bsnm{Heister}, \binits{T.}},
\bauthor{\bsnm{Heltai}, \binits{L.}},
\bauthor{\bsnm{K{\"o}cher}, \binits{U.}},
\bauthor{\bsnm{Kronbichler}, \binits{M.}},
\bauthor{\bsnm{Maier}, \binits{M.}},
\bauthor{\bsnm{Munch}, \binits{P.}},
\bauthor{\bsnm{Pelteret}, \binits{J.-P.}},
\bauthor{\bsnm{Proell}, \binits{S.}},
\bauthor{\bsnm{Simon}, \binits{K.}},
\bauthor{\bsnm{Turcksin}, \binits{B.}},
\bauthor{\bsnm{Wells}, \binits{D.}},
\bauthor{\bsnm{Zhang}, \binits{J.}}:
\batitle{The \texttt{deal.II} library, version 9.3}.
\bjtitle{J. Numer. Math.}
\bvolume{29}(\bissue{3}),
\bfpage{171}--\blpage{186}
(\byear{2021})
\doiurl{10.1515/jnma-2021-0081}
\end{barticle}
\endbibitem

\bibitem[\protect\citeauthoryear{{T}rilinos~{P}roject
  {T}eam}{2020}]{team2020trilinos}
\begin{botherref}
\oauthor{\bsnm{{T}eam}, \binits{T.}}:
The {T}rilinos {P}roject {W}ebsite.
(2020).
\url{https://trilinos.github.io}
\end{botherref}
\endbibitem

\bibitem[\protect\citeauthoryear{Ainsworth and
  Oden}{2000}]{ainsworth2000posteriori}
\begin{bbook}
\bauthor{\bsnm{Ainsworth}, \binits{M.}},
\bauthor{\bsnm{Oden}, \binits{J.T.}}:
\bbtitle{A Posteriori Error Estimation in Finite Element Analysis}.
\bsertitle{Pure and Applied Mathematics}.
\bpublisher{John Wiley \& Sons, Inc.},
\blocation{New York}
(\byear{2000})
\end{bbook}
\endbibitem

\bibitem[\protect\citeauthoryear{Ba\v{n}as and
  N{\"u}rnberg}{2008}]{banas2008adaptive}
\begin{barticle}
\bauthor{\bsnm{Ba\v{n}as}, \binits{L.}},
\bauthor{\bsnm{N{\"u}rnberg}, \binits{R.}}:
\batitle{Adaptive finite element methods for {C}ahn--{H}illiard equations}.
\bjtitle{J. Comput. Appl. Math.}
\bvolume{218}(\bissue{1}),
\bfpage{2}--\blpage{11}
(\byear{2008})
\doiurl{10.1016/j.cam.2007.04.030}
\end{barticle}
\endbibitem

\bibitem[\protect\citeauthoryear{Zhang et~al.}{2019}]{zhang2019stress}
\begin{barticle}
\bauthor{\bsnm{Zhang}, \binits{K.}},
\bauthor{\bsnm{Li}, \binits{Y.}},
\bauthor{\bsnm{Wang}, \binits{F.}},
\bauthor{\bsnm{Zheng}, \binits{B.}},
\bauthor{\bsnm{Yang}, \binits{F.}}:
\batitle{Stress effect on self-limiting lithiation in silicon-nanowire
  electrode}.
\bjtitle{Appl. Phys. Express}
\bvolume{12}(\bissue{4}),
\bfpage{045004}
(\byear{2019})
\doiurl{10.7567/1882-0786/ab0ce8}
\end{barticle}
\endbibitem

\bibitem[\protect\citeauthoryear{Chan et~al.}{2007}]{chan2007high-performance}
\begin{barticle}
\bauthor{\bsnm{Chan}, \binits{C.K.}},
\bauthor{\bsnm{Peng}, \binits{H.}},
\bauthor{\bsnm{Liu}, \binits{G.}},
\bauthor{\bsnm{McIlwrath}, \binits{K.}},
\bauthor{\bsnm{Zhang}, \binits{X.F.}},
\bauthor{\bsnm{Huggins}, \binits{R.A.}},
\bauthor{\bsnm{Cui}, \binits{Y.}}:
\batitle{High-performance lithium battery anodes using silicon nanowires}.
\bjtitle{Nat. Nanotechnol.}
\bvolume{3}(\bissue{1}),
\bfpage{31}--\blpage{35}
(\byear{2007})
\doiurl{10.1038/nnano.2007.411}
\end{barticle}
\endbibitem

\bibitem[\protect\citeauthoryear{Davis}{2004}]{davis2004algorithm}
\begin{barticle}
\bauthor{\bsnm{Davis}, \binits{T.A.}}:
\batitle{Algorithm 832: {UMFPACK} {V}4.3---an unsymmetric-pattern multifrontal
  method}.
\bjtitle{ACM Trans. Math. Software}
\bvolume{30}(\bissue{2}),
\bfpage{196}--\blpage{199}
(\byear{2004})
\doiurl{10.1145/992200.992206}
\end{barticle}
\endbibitem

\bibitem[\protect\citeauthoryear{McDowell et~al.}{2016}]{mcdowell2016mechanics}
\begin{barticle}
\bauthor{\bsnm{McDowell}, \binits{M.T.}},
\bauthor{\bsnm{Xia}, \binits{S.}},
\bauthor{\bsnm{Zhu}, \binits{T.}}:
\batitle{The mechanics of large-volume-change transformations in high-capacity
  battery materials}.
\bjtitle{Extreme Mech. Lett.}
\bvolume{9},
\bfpage{480}--\blpage{494}
(\byear{2016})
\doiurl{10.1016/j.eml.2016.03.004}
\end{barticle}
\endbibitem

\bibitem[\protect\citeauthoryear{Zhang and Kamlah}{2018}]{zhang2018sodium}
\begin{barticle}
\bauthor{\bsnm{Zhang}, \binits{T.}},
\bauthor{\bsnm{Kamlah}, \binits{M.}}:
\batitle{Sodium ion batteries particles: {Phase}-field modeling with coupling
  of {C}ahn--{H}illiard equation and finite deformation elasticity}.
\bjtitle{J. Electrochem. Soc.}
\bvolume{165}(\bissue{10}),
\bfpage{1997}--\blpage{2007}
(\byear{2018})
\doiurl{10.1149/2.0141810jes}
\end{barticle}
\endbibitem

\bibitem[\protect\citeauthoryear{Schomburg
  et~al.}{2023}]{schomburg2023characterization}
\begin{barticle}
\bauthor{\bsnm{Schomburg}, \binits{F.}},
\bauthor{\bsnm{Drees}, \binits{R.}},
\bauthor{\bsnm{Kurrat}, \binits{M.}},
\bauthor{\bsnm{Danzer}, \binits{M.A.}},
\bauthor{\bsnm{R\"oder}, \binits{F.}}:
\batitle{Characterization of the solid--electrolyte interphase growth during
  cell formation based on differential voltage analysis}.
\bjtitle{Energy Technol.}
\bvolume{11}(\bissue{5}),
\bfpage{2200688}
(\byear{2023})
\doiurl{10.1002/ente.202200688}
\end{barticle}
\endbibitem

\bibitem[\protect\citeauthoryear{Castelli and
  D\"orfler}{2019}]{castelli2019efficient}
\begin{bchapter}
\bauthor{\bsnm{Castelli}, \binits{G.}},
\bauthor{\bsnm{D\"orfler}, \binits{W.}}:
\bctitle{An efficient matrix-free finite element solver for the
  {C}ahn--{H}iliard equation}.
In: \beditor{\bsnm{Gleim}, \binits{T.}},
\beditor{\bsnm{Lange}, \binits{S.}} (eds.)
\bbtitle{Proceedings of 8th {GACM} Colloquium on Computational Mechanics},
pp. \bfpage{441}--\blpage{444}.
\bpublisher{Kassel University Press},
\blocation{Kassel}
(\byear{2019}).
\doiurl{10.19211/KUP978737650939}
\end{bchapter}
\endbibitem

\bibitem[\protect\citeauthoryear{Schenk and
  G{\"a}rtner}{2011}]{schenk2011pardiso}
\begin{bbook}
\bauthor{\bsnm{Schenk}, \binits{O.}},
\bauthor{\bsnm{G{\"a}rtner}, \binits{K.}}:
In: \beditor{\bsnm{Padua}, \binits{D.}} (ed.)
\bbtitle{{PARDISO}},
pp. \bfpage{1458}--\blpage{1464}.
\bpublisher{Springer},
\blocation{Boston, MA}
(\byear{2011}).
\doiurl{10.1007/978-0-387-09766-4_90}
\end{bbook}
\endbibitem

\end{thebibliography}




\appendix

\addcontentsline{toc}{section}{Appendices}
\section*{Appendices}

\renewcommand{\thesection}{A}
\section{Abbreviations and Symbols}
\label{app:abbreviations_and_symbols}
\hspace*{-0.60cm}
\vspace{-0.4cm}
\begin{longtable}[l]{p{2.27cm}p{5.5cm}}
  \multicolumn{2}{l}{\textbf{Abbreviations}} \\
  \ \\
  AD & automatic differentiation \\
  aSi & amorphous silicon \\
  BDF & backward differentiation formula \\
  C-rate & charging rate \\
  DAE & differential algebraic equation \\
  DOF & degree of freedom \\
  GSV & Green--St-Venant \\
  KKT & Karush--Kuhn--Tucker \\
  MPI & message passing interface \\
  NDF & numerical differentiation formula \\
  OCV & open-circuit voltage \\
  pmv & partial molar volume \\
  SOC & state of charge
\end{longtable}

\hspace*{-0.60cm}
\vspace{-0.4cm}
\begin{longtable}[l]{p{2.27cm}p{10.5cm}}
  \textbf{Symbol} & \textbf{Description} \\
  \ \\
  \multicolumn{2}{l}{Latin symbols}
  \vspace{0.2cm}\\
  $c$ & concentration\\
  $\tens{C}_\el$ & right Cauchy-Green tensor \\
  $\Ctensor$ & fourth-order stiffness tensor \\
  $\tens{D}_{\pl}$ & plastic strain rate\\
  $E$ & Young's modulus\\
  $\tens{E}_\el$ & elastic strain tensor \\
  $F_\text{Y}$ & yield function \\
  $\tens{F} = \gradL\vect{\Phi}$ & deformation gradient \\
  $\tens{F} = \tens{F}_\ch\tens{F}_\el\tens{F}_\pl$ & multiplicative
  decomposition of $\tens{F}$ \\
  $\tens{F}_\ch$ & chemical part of the deformation gradient \\
  $\tens{F}_\el$ & elastic part of the deformation gradient \\
  $\tens{F}_\pl$ & plastic part of the deformation gradient \\
  $\mathrm{Fa}$ & Faraday constant \\
  $G$ & shear modulus / second Lam\'{e} constant \\
  $\tens{I}$ & identity tensor \\
  $\tensorIpi$ & approximation for the linearization of $\tens{P}_\Pi$\\
  $J = J_\ch J_\el J_\pl$ & multiplicative
  decomposition of volume change \\
  $K$ & bulk modulus \\
  ${m}$ & scalar valued mobility \\
  $\tens{M}$ & Mandel stress tensor \\
  $\vect{n}$, $\vect{n}_0$ & normal vector on $\Omega$, $\Omega_0$ \\
  $N$ & number of nodes of ${V}_h$ \\
  $\vect{N}$ & lithium flux \\
  ${N}_\text{ext}$ & external lithium flux \\
  $\tens{P}$ & first Piola--Kirchhoff stress tensor \\
  $\tens{P}_\Pi$ & projector on admissible stresses \\
  $\tens{R}_\el$ & rotational part of polar decomposition of $\tens{F}_\el$ \\
  $\tens{S}_\el$ & deformation part of polar decomposition of $\tens{F}_\el$ \\
  $t$ & time \\
  $t_n$ & time step \\
  $t_\text{cycle}$ & cycle time \\
  $U_\text{OCV}$ & OCV curve \\
  $\vect{u} = \vect{x} - \vect{X}_0$ & time dependent displacement vector \\
  $\vect{u}_h$ & discrete displacement vector or algebraic representation \\
  $v_\text{pmv}$ & partial molar volume of lithium \\
  ${V}$ & scalar valued function space \\
  $\vect{V}$ & vector valued function space \\
  $\vect{V}^*$ & subset of $\vect{V}$ \\
  $\vect{x} = \vect{\Phi} \left(t, \vect{X_0}\right) $ & position in Eulerian
  domain \\
  $\vect{X}_0$ & initial placement \\
  & \\
  %
  \multicolumn{2}{l}{Greek symbols} \\
  $\alpha_{k_n} > 0$ & coefficient for adaptive time
  discretization\\
  $\varepsilon_\pl^\text{eq}$ & equivalent plastic strain \\
  $\gamma^\iso$ & isotropic hardening parameter \\
  $\lambda$ & first Lam\'{e} constant \\
  $\lambda_\ch$ & factor of concentration induced deformation gradient \\
  $\mu$ & chemical potential \\
  $\nu$ & Poisson's ratio \\
  $\Omega$ & Eulerian domain \\
  $\Omega_0$ & Lagrangian domain \\
  $\varphi$ & scalar valued test function \\
  $\psi=\psi_\ch+\psi_\el$ & total free energy \\
  $\psi_\ch$ & chemical part of free energy \\
  $\psi_\el$ & elastic part of free energy \\
  $\rho_0$ & density \\
  ${\sigma}_\text{F}$ & yield stress function\\
  ${\sigma_\text{Y}}(c)$ & concentration dependent yield stress \\
  $\boldsymbol{\sigma}$ & Cauchy stress tensor\\
  $\tau_n$ & time step size at time $t_n$\\
  $\vect{\xi}_j$ & vector valued test function of node $j$\\
  ${\xi}_j$ & scalar basis function: nonzero entry of $\vect{{\xi}_j}$ of node
  $j$ 	\\
  & \\
  %
  \multicolumn{2}{l}{Mathematical symbols}	\\
  $\partial \Omega_0$ & boundary of $\Omega_0$ \\
  $\gradL$ & gradient vector in Lagrangian domain \\
  $\square \!:\! \tilde{\square}$ & reduction of two dimensions of two tensors
  $\square$ and $\tilde{\square}$\\
  ${\mathds{A}} [ \tens{B}]$
  &
  reduction of the last two dimensions of the fourth order
  tensor~${\mathds{A}}$
  and the second order tensor $\tens{B}$
  \\
  $\partial_\square$ & partial derivative with respect to $\square$ \\
  $\triangle_n \square$ & difference of next and current time step:
  $\square^\npo - \square^n$
  or
  $\square_\npo - \square_n$\\
  $\norm{\square}{}{}$ & Frobenius norm for tensors \\
  & \\
  %
  %
  \multicolumn{2}{l}{Indices} \\
  ${\square}^{\devi}$ & deviatoric part of $\square$\\
  ${\square}^{\trial}$ & trial part of $\square$\\
  ${\square}_{0}$ & considering variable in Lagrangian domain or initial
  condition
  \\
  ${\square}_\text{ch}$ & chemical part of~$\square$ \\
  ${\square}_\text{el}$ & elastic part of~$\square$ \\
  ${\square}_{h}$ & finite dimensional function of $\square$ or algebraic
  representation of $\square$ with respect to basis function \\
  ${\square}_{\text{max}}$ & maximal part of $\square$\\
  ${\square}_{\text{min}}$ & minimal part of $\square$\\
  ${\square}_\text{pl}$ & plastic part of~$\square$ \\
\end{longtable}


\renewcommand{\thesection}{B}
\section{Derivation of First Piola--Kirchhoff Tensor}
\label{app:derivation_p}
The first Piola--Kirchhoff stress tensor is defined as the derivative of the
free energy density by the deformation gradient~$\tens{F}$,
compare~\cite[Section~6.1]{holzapfel2010nonlinear},
as
\begin{align}
  \label{eq:derivation_first_piola}
  \tens{P}
  &
  = \ptl_\tens{F} (\rho_0 \psi)
  = \ptl_\tens{F} (\rho_0 \psi_\el)
  =
  \ptl_\tens{F} \tens{C}
  \big[
  \ptl_{\tens{C}}(\rho_0 \psi_\el)
  \big]
  \\
  &
  = 2\tens{F}
  \ptl_{\tens{C}}
  (\rho_0 \psi_\el)
  = 2\tens{F} \ptl_{\tens{C}} \tens{E}_\el
  \big[
  \Ctensor \big[
  \tens{E}_\el\big]
  \big].
\end{align}
The chain rule for the second equality
in~\cref{eq:derivation_first_piola}
can be computed as stated,
since
$\rho_0 \psi_\el =
\halb
\tens{E}_\text{el} \dualreduction \Ctensor \left[ \tens{E}_\text{el} \right]
$,
$\ptl_{\tens{C}}(\rho_0 \psi_\el)$
is a symmetric tensor.
Further computation with a symmetric tensor~$\tens{S}$
and $n,m = 1, \dots, d$ leads to:
\begin{align}
  \frac{\ptl}{\ptl{F}_{mn}}
  \tens{C}
  \big[
  \tens{S}
  \big]
  &
  =
  \frac{\ptl
  }
  {\ptl{F}_{mn}}
  (
  \tens{C} \dualreduction \tens{S}
  )
  \qq \qq \vspace{-0.1cm}
  =
  \frac{\ptl
  }
  {\ptl{F}_{mn}}
  (
  \tens{F}^\trp \tens{F} \dualreduction \tens{S}
  )
  \\
  &
  =
  \frac{\ptl
  }
  {\ptl{F}_{mn}}
  \tr \,
  (
  (\tens{F}^\trp \tens{F})^\trp \tens{S}
  )
  \q
  =
  \frac{\ptl
  }
  {\ptl{F}_{mn}}
  \tr \,
  (
  \tens{F}^\trp \tens{F} \tens{S}
  )
  \\
  &
  =
  \frac{\ptl
  }
  {\ptl{F}_{mn}}
  \sum_{j}
  \sum_{l}
  (
  \tens{F}^\trp \tens{F}
  )_{jl}
  {S}_{lj}
  \\
  &
  =
  \frac{\ptl
  }
  {\ptl{F}_{mn}}
  \sum_{j}
  \sum_{l}
  \sum_{k}
  {F}_{kj}
  {F}_{kl}
  {S}_{lj}
  \\
  &
  =
  \sum_{j}
  \sum_{l}
  \sum_{k}
  \big(
  \delta_{km}
  \delta_{jn}
  {F}_{kl}
  {S}_{lj}
  \big)
  +
  \big(
  {F}_{kj}
  \delta_{km}
  \delta_{ln}
  {S}_{lj}
  \big)
  \\
  &
  =
  \sum_{l}
  \big(
  {F}_{ml}
  {S}_{ln}
  \big)
  +
  \sum_{j}
  \big(
  {F}_{mj}
  {S}_{nj}
  \big)
  \\
  &
  =
  \sum_{l}
  \big(
  {F}_{ml}
  {S}_{ln}
  \big)
  +
  \sum_{j}
  \big(
  {F}_{mj}
  (
  \tens{S}^\trp
  )_{jn}
  \big)
  \\
  &
  =
  (
  \tens{F}
  \tens{S}
  )_{mn}
  +
  (
  \tens{F}
  \tens{S}^\trp
  )_{mn}
  =
  2(
  \tens{F} \tens{S}
  )_{mn}.
\end{align}
Using the definition of the logarithmic elastic strain, it follows:
\begin{align}
  \tens{E}_\el
  &
  = \ln
  \left( \tens{U}_\el \right)
  = \ln \left( \sqrt{\tens{C}_\el} \right)
  = \ln   \left( \sqrt{ \tens{F}_\el^{\trp}  \tens{F}_\el } \right)
  =
  \ln \left( \lambda_\ch^{-1} \sqrt{ \left(
    \tens{F}_\pl^{-1} \right)^{\trp} \tens{F}^{\trp} \tens{F}
    \tens{F}_\pl^{-1}  } \right)
  \\
  &
  = \ln \left( \lambda_\ch^{\minusone}
  \sqrt{ \left(
    \tens{F}_\pl^{\minusone} \right)^{\trp} \tens{C}
    \tens{F}_\pl^{\minusone}
  }
  \right).
\end{align}
Now, we have
\begin{align}
  \ptl_{\tens{C}} \tens{E}_\el
  \big[
  \Ctensor \big[
  \tens{E}_\el\big]
  \big]
  &
  =
  \left\{
  \left(
  \lambda_\ch^{\minusone} \sqrt{ \left(
    \tens{F}_\pl^{\minusone} \right)^{\trp} \tens{C}
    \tens{F}_\pl^{\minusone}
  }
  \right)^{\minusone}
  \halb
  \lambda_\ch^{\minusone}
  \left( \left(
  \tens{F}_\pl^{\minusone} \right)^{\trp} \tens{C}
  \tens{F}_\pl^{\minusone} \right)^{\minus 1/2}
  \right.
  \\
  &
  \q
  \left.
  \left( \left(
  \tens{F}_\pl^{\minusone} \right)^{\trp}
  \tens{I}
  \tens{F}_\pl^{\minusone}
  \right)
  \right\}
  \Ctensor
  \big[\tens{E}_\el\big]
  \\
  &
  =
  \halb
  \left\{
  \left( \left(
  \tens{F}_\pl^{\minusone} \right)^{\trp} \tens{C}
  \tens{F}_\pl^{\minusone}   \right)^{\minusone}
  \left( \left(
  \tens{F}_\pl^{\minusone} \right)^{\trp}
  \tens{F}_\pl^{\minusone}   \right) \right\}
  \Ctensor
  \big[\tens{E}_\el\big].
\end{align}
In total, the first Piola--Kirchhoff tensor~$\tens{P}$ is computed via
\begin{align}
  \tens{P}
  &
  =
  2\tens{F}
  \halb
  \left\{
  \left( \left(
  \tens{F}_\pl^{\minusone} \right)^{\trp} \tens{C}
  \tens{F}_\pl^{\minusone}   \right)^{\minusone}
  \left( \left(
  \tens{F}_\pl^{\minusone} \right)^{\trp}
  \tens{F}_\pl^{\minusone}   \right) \right\}
  \Ctensor
  \big[\tens{E}_\el\big]
  \\
  &
  =
  \tens{F}
  \left\{
  \left( \left(
  \tens{F}_\pl^{\minusone} \right)^{\trp}
  \tens{F}_\pl^\trp \tens{F}_\el^\trp \tens{F}_\ch^\trp
  \tens{F}_\ch \tens{F}_\el \tens{F}_\pl
  \tens{F}_\pl^{\minusone}   \right)^{\minusone}
  \left( \left(
  \tens{F}_\pl^{\minusone} \right)^{\trp}
  \tens{F}_\pl^{\minusone}   \right) \right\}
  \Ctensor
  \big[\tens{E}_\el\big]
  \\
  &
  =
  \tens{F}
  \left(
  \tens{F}_\text{rev}^{\trp}
  \tens{F}_\text{rev}
  \right)^{\minusone}
  \left( \left(
  \tens{F}_\pl^{\minusone} \right)^{\trp}
  \tens{F}_\pl^{\minusone}   \right)
  \Ctensor
  \big[\tens{E}_\el\big].
\end{align}

\renewcommand{\thesection}{C}
\section{Linearization of Projector}
\label{app:derivation}

We estimate the formulation of the exact Jacobian for the rate-independent
plastic approach as
in~\citet[Appendix~C.6]{di-leo2015chemo-mechanics}
with
the formal linearization of the projector~$\tens{P}_\Pi$
of~\cref{eq:projector_rate_independent}
around $\tens{E}_\el^\trial$~\cite{simo1985consistent, frohne2016efficient}
\begin{align}
  \tensorIpi
  \big(\con^\npo, \gradL \vect{u}^\npo, \Fpl^n,
  \varepsilon_\pl^{\text{eq},n} \big)
  = \Fp{\tens{M}^\npo}{\tens{E}_\el^\trial}.
\end{align}
Then, the consistent algorithmic modulus is given as
\begin{equation}
  \tensorIpi
  =
  \left\{\begin{array}{ll} \Ctensor_{\el},
    &
    \NOrm{\Ctensor \left[\tens{E}_\el^{\trial, \devi}
      \right]}{}{}
    \leq
    \sigma_\text{Y}(c^\npo
    ) + \gamma^\iso \varepsilon_\pl^{\text{eq},n},
    \\[1em]
    \Ctensor_{\pl},
    &
    \NOrm{\Ctensor \left[\tens{E}_\el^{\trial, \devi}
      \right]}{}{}
    >
    \sigma_\text{Y}(c^\npo
    ) + \gamma^\iso
    \varepsilon_\pl^{\text{eq},n}.\end{array}\right.
\end{equation}
With the expressions for $ \tens{M}^\npo$
in the projector~$\tens{P}_\Pi$
of~\cref{eq:projector_rate_independent}
for the elastic and plastic cases,
the elastic tangent~$\Ctensor_{\el}$
and plastic tangent~$\Ctensor_{\pl}$
can be derived, respectively:\\
\textbf{Elastic tangent.}
It immediately reads
\begin{alignat}{2}
  & \Ctensor_{\el}
  = \Fp{\tens{M}^\npo}{\tens{E}_\el^\trial}
  = \Fp{\tens{M}^\trial}{\tens{E}_\el^\trial}
  &&
  = \Fp{
    \Big(
    2 G \tens{E}_\el^\trial
    + (K - \frac{2}{3}G) \tr \left({\tens{E}_\el^\trial}\right)
    \tens{I}
    \Big)
  }
  {\tens{E}_\el^\trial}
  \\
  & &&= 2G \mathds{Id} + (K - \frac{2}{3} G) \tens{I}
  \otimes \tens{I}
  \\
  & &&= 2G \left(\mathds{Id} - \frac{1}{3} \tens{I}
  \otimes \tens{I}\right) +
  K
  \tens{I} \otimes \tens{I}
  \\
  & &&= \Ctensor_{G}+\Ctensor_{K}.
\end{alignat}

\noindent \textbf{Plastic tangent.}
With the abbreviations
$a = \frac{\gamma^{\iso}}{2G+\gamma^{\iso}}$,
$b =
\Big(1-
\frac{\gamma^{\iso}}{2G+\gamma^{\iso}}
\big(1-\frac{2G}{\sigma_\text{Y}(c^\npo
  )}\varepsilon_\pl^{\text{eq},n} \big)
\Big)$
it follows
\begin{equation}
  \begin{aligned}
    \Ctensor_{\pl}
    &
    =
    \Fp{\tens{M}^\npo}{\tens{E}_\el^\trial}
    =
    \Fp{
      \big(
      \tens{M}^\trial - 2 G \tau_n \tens{D}_\pl^\npo \big)
    }{\tens{E}_\el^\trial}
    =
    \Fp{
      \big(
      \tens{M}^\trial - 2 G \triangle_n \varepsilon_\pl^\text{eq}
      \frac{\tens{M}^{\trial,\devi}}{\norm{\tens{M}^{\trial, \devi}}{}{}}\big)
    }
    {\tens{E}_\el^\trial}
    \\
    &=
    \Fp{
      \Bigg( \drittel
      \tr \left({\tens{M}^\trial} \right) \tens{I}
      + \left(
      a
      +b
      \frac{\sigma_\text{Y}(c^\npo
        )}{\norm{\tens{M}^{\trial, \devi}}{}{}}
      \right)
      \tens{M}^{\trial, \devi}
      \Bigg)
    }
    {
      \tens{E}_\el^\trial
    }
    \\
    &=
    \Fp{
      \left(a+b\frac{\sigma_\text{Y}(c^\npo
        )}{\norm{\tens{M}^{\trial,
            \devi}}{}{}}\right)
    }
    {\tens{E}_\el^{\trial}}
    \tens{M}^{\trial, \devi}
    \\
    & \quad+
    \left(
    a+b\frac{\sigma_\text{Y}(c^\npo
      )}{\norm{\tens{M}^{\trial, \devi}}{}{}}
    \right)
    \Fp{
      \tens{M}^{\trial, \devi}
    }
    {\tens{E}_\el^{\trial}}
    + \drittel
    \Fp{
      \tr\big({\tens{M}^{\trial}}\big)
    }
    {\tens{E}_\el^{\trial}}\tens{I}
    \\
    &= b\sigma_\text{Y}(c^\npo
    )
    \Fp{
      \left(\tens{M}^{\trial, \devi}:\tens{M}^{\trial, \devi}\right)^{-\halb}
    }
    {\tens{E}_{\el}^{\trial}}
    \tens{M}^{\trial, \devi}
    \\
    & \quad
    +
    \left(a+b\frac{\sigma_\text{Y}(c^\npo
      )}{\norm{\tens{M}^{\trial,
          \devi}}{}{}}\right)
    \Fp{\tens{M}^{\trial, \devi}}{\tens{E}_\el^{\trial}}
    + \drittel
    \Fp{
      \tr\big({\tens{M}^{\trial}}\big)}
    {\tens{E}_\el^{\trial, \devi}}\tens{I}.
  \end{aligned}
\end{equation}
Introducing a modified~$\tilde{\Ctensor}_{\el}$ gives
\begin{equation}
  \tilde{\Ctensor}_{\el}
  =
  a
  \Fp{\tens{M}^{\trial, \devi}}
  {\tens{E}_\el^{\trial, \devi}}
  +
  \drittel
  \Fp{\tr\left({\tens{M}^{\trial}}\right)}
  {\tens{E}_\el^{\trial}}\tens{I}
  =
  a2G(\mathds{Id} - \drittel\tens{I} \otimes
  \tens{I})
  +K\tens{I} \otimes \tens{I}
  =
  a\Ctensor_{G}+\Ctensor_{K} .
\end{equation}
Further, it states
\begin{equation}
  \begin{aligned}
    &
    \Fp{
      \Big(\tens{M}^{\trial, \devi}:\tens{M}^{\trial, \devi}
      \Big)^{-\halb}
    }
    {\tens{E}_\el^{\trial}}
    \\
    &
    \qquad \qquad
    = -\halb
    \Big(\tens{M}^{\trial, \devi}
    \dualreduction
    \tens{M}^{\trial, \devi}\Big)^{-\frac{3}{2}}
    \Bigg(
    \Fp{\tens{M}^{\trial, \devi}}{\tens{E}_\el^{\trial}}
    \big[\tens{M}^{\trial, \devi}\big]
    +
    \Fp{\tens{M}^{\trial, \devi}}{\tens{E}_\el^{\trial}}
    \big[\tens{M}^{\trial, \devi}\big]
    \Bigg)
    \\
    &
    \qquad \qquad
    = -\Big(
    \tens{M}^{\trial, \devi}
    \dualreduction
    \tens{M}^{\trial, \devi}
    \Big)^{-\frac{3}{2}}
    \Fp{\tens{M}^{\trial, \devi}}{\tens{E}_\el^{\trial}}
    \big[\tens{M}^{\trial, \devi}\big]
    \\
    &
    \qquad \qquad
    = -2G
    \left(\tens{M}^{\trial, \devi}
    \dualreduction
    \tens{M}^{\trial,\devi}
    \right)^{-\frac{3}{2}}\tens{M}^{\trial, \devi}.
  \end{aligned}
\end{equation}
Finally, the plastic tangent can be specified by
\begin{alignat}{2}
  &
  \Ctensor_{\pl}
  &&
  = \tilde{\Ctensor}_{\el}
  - b \sigma_\text{Y}(c^\npo
  ) 2 G
  \Big(\tens{M}^{\trial,\devi}
  \dualreduction
  \tens{M}^{\trial, \devi}\Big)^{-\frac{3}{2}}
  \tens{M}^{\trial, \devi} \otimes \tens{M}^{\trial, \devi}
  \\
  & && \quad
  +b \frac{\sigma_\text{Y}(c^\npo
    )}
  {\norm{\tens{M}^{\trial, \devi}}{}{}}
  \Fp{\tens{M}^{\trial, \devi}}{\tens{E}_\el^{\trial}}
  \\
  & &&
  = \tilde{\Ctensor}_{\el}
  + \frac{b\sigma_\text{Y}(c^\npo
    )}
  {\norm{\tens{M}^{\trial, \devi}}{}{}}
  \Bigg(
  \Ctensor_{G}-2G
  \frac{\tens{M}^{\trial, \devi}\otimes\tens{M}^{\trial, \devi}}
  {\norm{\tens{M}^{\trial, \devi}}{}{2}}
  \Bigg)
  \\
  & &&
  = \tilde{\Ctensor}_{\el}
  + \frac{
    \Big(1-\frac{\gamma^{\iso}}
    {2G+\gamma^{\iso}}
    \big(1-\frac{2G}{\sigma_\text{Y}(c^\npo
      )} \varepsilon_\pl^{\text{eq},n}\big)
    \Big)
    \sigma_\text{Y}(c^\npo
    )}
  {\norm{\tens{M}^{\trial, \devi}}{}{}}
  \Bigg(\Ctensor_{G}
  \Bigg.
  \\
  & &&
  \quad \quad \quad \quad \quad \quad \quad
  \quad \quad \quad \quad \quad \quad \quad
  \quad \quad \quad \quad \quad \quad \quad
  \Bigg.
  -2G
  \frac{\tens{M}^{\trial, \devi}\otimes\tens{M}^{\trial, \devi}}
  {\norm{\tens{M}^{\trial, \devi}}{}{2}}
  \Bigg).
\end{alignat}
With the deviatoric trial Mandel stress
$\tens{M}^{\trial, \devi} = \Ctensor[\tens{E}_\el^{\trial, \devi}
]
$,
the
linearization of $\tens{P}_\Pi$
around $\tens{E}_\el^{\trial}$
can be expressed by~\cref{eq:linearization_rate_independent}.
A similar calculation leads for the rate-dependent case to
\begin{subequations}
  \label{eq:viscoplastic_linearization}
  \begin{align*}
    \tensorIpi
    \big(\con^\npo, \gradL \vect{u}^\npo, \Fpl^n,
    \triangle_n \varepsilon_\pl^{\text{eq}} \big)
    = \hspace{22.0em}\\[-2.7em]
  \end{align*}
  \begin{empheq}[
    left=
    \empheqlbrace
    ]{alignat=3}
    &
    \Ctensor_{K} + \Ctensor_{G},
    &&
    \NORm{\Ctensor \left[\tens{E}_\el^{\trial, \devi}
      \right]}{}{}
    \leq
    \sigma_\text{Y}(c^\npo
    ),
    \\
    &
    \frac{\minus 2G
      \triangle_n
      \varepsilon_\pl^{\text{eq}}}{\norm{\tens{M}^{\trial,
          \devi}}{}{}}
    \Bigg(
    \Ctensor_{G} - 2G
    \frac{\tens{M}^{\trial, \devi}\otimes
      \tens{M}^{\trial, \devi}}{\norm{\tens{M}^{\trial, \devi}}{}{2}}
    \Bigg)
    &&
    \nonumber
    \\
    &
    + \Ctensor_{G} + \Ctensor_{K},
    &&
    \NORm{\Ctensor \left[\tens{E}_\el^{\trial, \devi}
      \right]}{}{}
    >
    \sigma_\text{Y}(c^\npo
    ).
  \end{empheq}
\end{subequations}

\renewcommand{\thesection}{D}
\section{Comparison of the Plastic Impact on Stress Development}
\label{app:comparison_stress}

\begin{figure}[b]
  \centering
  {\includegraphics[width=0.5\textwidth,page=1]
    {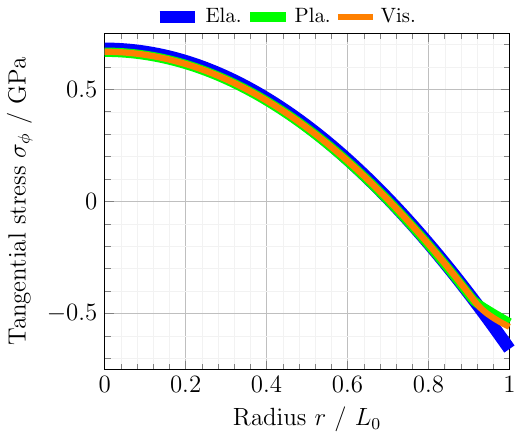}}
  \caption{
    Numerical results
    for the elastic (Ela.), plastic (Pla.) and
    viscoplastic (Vis.) approaches of
    the 1D radial symmetric case at~$\text{SOC} = 0.1$
    for the elastic, plastic and viscoplastic approach of the tangential Cauchy
    stress~$\sigma_\varphi$
    over the particle radius~$r$.
  }
  \label{fig:comparison_stresses_plastic_literature}
\end{figure}

\cref{fig:comparison_stresses_plastic_literature} shows the 1D radial symmetric
results of the tangential Cauchy stress~$\sigma_\varphi$
over the particle radius~$r$
for the elastic, plastic and viscoplastic case
at $\text{SOC}=0.1$.
The area at the particle surface is plastically deformed for the non elastic
approaches and is qualitatively comparable to numerical results in Fig.~4(c)
of~\cite{mcdowell2016mechanics}
and Fig.~5(d) of~\cite{zhao2011large}.
With a reduction of the yield stress~$\sigma_\text{Y}(c)$,
the plastically deformed range would be larger.

\begin{figure}[t]
  \centering
  \includegraphics[scale=0.75,
  page=1]{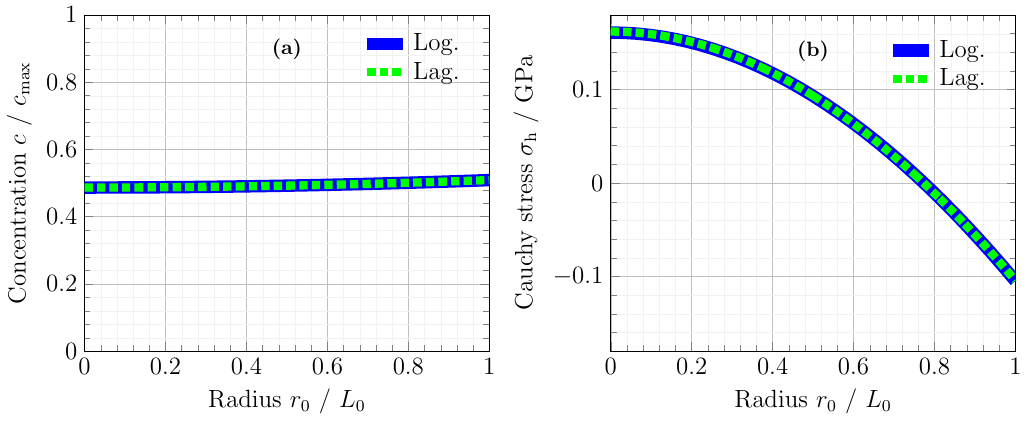}
  \caption{Comparison of concentration~$c$
    in~(a) and hydrostatic Cauchy
    stress~$\sigma_{\text{h}} = 1/3 \big(\sigma_\text{r} + 2
    \sigma_\varphi\big)$
    in~(b)
    over particle radius for logarithmic
    strain (Log.)
    vs. Green--St-Venant strain / Lagrangian strain
    (Lag.)
    at $\text{SOC}=0.5$.}
  \label{fig:comparison_log_lag_c_sigma_h}
\end{figure}

\renewcommand{\thesection}{E}
\section{Comparison of Logarithmic Strain vs. Green--St-Venant Strain}
\label{app:comparison_strain}

Comparing the logarithmic strain approach and the
Green--St-Venant strain
approach as discussed in~\cref{subsec:finite_deformation}
and~\cref{subsubsec:1d_spherical_symmetry},
both approaches result in equal
outputs, displayed in~\cref{fig:comparison_log_lag_c_sigma_h}.
This corresponds to the findings in~\cite{zhang2018sodium}.

\begin{figure}[b]
  \centering
  {\includegraphics[width=0.5\textwidth,page=1]
    {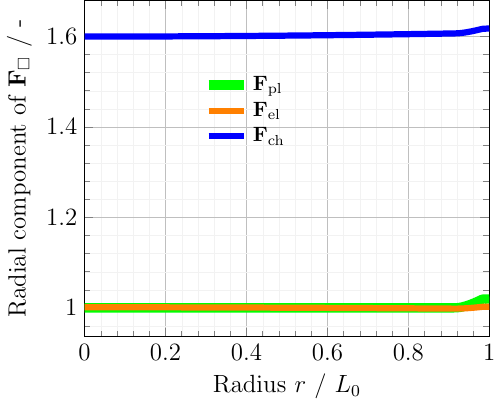}}
  \caption{Comparison of the radial component of the plastic, elastic and
    chemical part of the deformation gradient tensor~$\tens{F}_\pl$,
    $\tens{F}_\el$
    and~$\tens{F}_\ch$
    over the particle radius~$r$
    at~$\text{SOC}=0.92$.
  }
  \label{fig:comparison_deformation_gradient}
\end{figure}

\renewcommand{\thesection}{F}
\section{Comparison of Radial Deformation Gradient}
\label{app:comparison_deformation_gradient}

In this section, we take a closer look at the different components of the three
different parts of the deformation gradient~$\tens{F}_\pl$,
$\tens{F}_\el$
and~$\tens{F}_\ch$.
At the start of the simulation, no displacement gradient is present and the
total deformation gradient tensor is the identity tensor.
In~\cref{fig:comparison_deformation_gradient},
we consider the radial part of the each deformation gradient tensor,
respectively.
It can be clearly seen that the chemical part makes the largest contribution to
the deformation, followed at a large distance by the plastic and the elastic
parts which stay close to one.
This finding justifies the used approach of a linear elastic theory.

\renewcommand{\thesection}{G}
\section{Comparison of Voltage Curve}
\label{app:comparison_voltage_curve}

\begin{figure}[t]
  \centering
  {\includegraphics[width=0.5\textwidth,page=1]
    {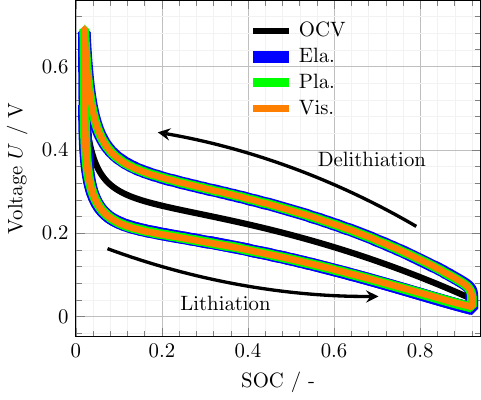}}
  \caption{
    Voltage curve~$U$
    of three half cycles for the elastic, plastic and viscoplastic approach
    over the~$\text{SOC}$ for three half cycles.
  }
  \label{fig:voltage}
\end{figure}

Following~\cite{di-leo2015diffusion-deformation}, we compute the voltage~$U$
with a Butler--Volmer condition at the particle surface to consider also
electrochemical surface influence.
Then,
the voltage~$U$
is given by
\begin{align}
  U
  =
  U_0
  - \frac{\mu_\text{surf}}{\mathrm{Fa}}
  -\frac{2 R_\text{gas} T}{\mathrm{Fa}} \sinh^\minusone
  \bigg( \frac{ \mathrm{Fa} \, N_\ext}{2 j_0(c_\text{surf}) }\bigg),
\end{align}
with some reference potential~$U_0$
depending on the counter-electrode and
a concentration-dependent exchange current
density~$j_0(c_\text{surf}) = k_0 \sqrt{c_\text{surf} (1-c_\text{surf})}$.
Here,~$k_0$ is a rate constant current density and is specific to the
material, $c_\text{surf}$ the non-dimensional surface concentration and
$\mu_\text{surf}$ the surface chemical potential.
For the moment, we set~$U_0$ to zero.
There is almost no difference on the resulting voltage~$U$
between the elastic, plastic and viscoplastic approaches,
compare~\cref{fig:voltage}.
However, the computed voltage values differs from the underlying voltage
curve~$U_\text{OCV}$ due to mechanical effects and the Butler--Volmer interface
condition.
In total, the numerical results are qualitatively comparable with the results
in~Fig.~7(a)
in~\cite{di-leo2015diffusion-deformation}.
We believe that the similarities of the three model approaches emerges from the
applied constant external lithium flux.

\end{document}